\documentclass[aip,jcp,11pt,reprint,floatfix,superscriptaddress,urlname]{revtex4-1}
\usepackage[utf8]{inputenc}
\usepackage{amsmath}
\usepackage{color}
\usepackage{graphicx}
\usepackage{comment}
\usepackage{caption}
\usepackage{subcaption}
\usepackage[normalem]{ulem}

\begin{document}
\title{Arithmetic circuit tensor networks, multivariable function representation, and high-dimensional integration}

\author{Ruojing Peng}
\affiliation{Division of Chemistry and Chemical Engineering, California Institute of Technology, Pasadena,
California 91125, USA}
\author{Johnnie Gray}
\affiliation{Division of Chemistry and Chemical Engineering, California Institute of Technology, Pasadena,
California 91125, USA}

\author{Garnet Kin-Lic Chan}
\affiliation{Division of Chemistry and Chemical Engineering, California Institute of Technology, Pasadena,
California 91125, USA}

\begin{abstract}
    Many computational problems can be formulated in terms of high-dimensional functions. Simple representations of such functions and resulting computations with them typically suffer from the ``curse of dimensionality'', an exponential cost dependence on dimension. Tensor networks provide a way to represent certain classes of high-dimensional functions with polynomial memory. This results in computations where the exponential cost is ameliorated or in some cases, removed, if the tensor network representation can be obtained. Here, we introduce a direct mapping from the arithmetic circuit of a function to arithmetic circuit tensor networks, avoiding the need to perform any optimization or functional fit. We demonstrate the power of the circuit construction in examples of multivariable integration on the unit hypercube in up to 50 dimensions, where the complexity of integration can be understood from the circuit structure. We find very favourable cost scaling compared to quasi-Monte-Carlo integration for these cases, and further give an example where efficient quasi-Monte-Carlo cannot be theoretically performed  without knowledge of the underlying tensor network circuit structure.
\end{abstract}

\maketitle

\section{Introduction}

High-dimensional multivariable functions (henceforth, multivariable functions) and their integrals appear in a multitude of areas, ranging from statistical and quantum many-body physics~\cite{popov2001,Yaglom1960,Pierre2000} to applications in machine learning~\cite{LEWIS1994148,liberty16,ElRafey2017,rs11020185,Tyagi2020,Kaltenecker2020}. Simple representations of such functions, for example, on a product grid, require exponential storage, and subsequent manipulation of the functions, e.g. in integration by quadrature, then requires exponential cost in dimension. Many techniques have been introduced to bypass this exponential cost. For example, high-dimensional integration is often carried out by Monte Carlo or quasi-Monte Carlo methods, which sample the function at a set of random or pre-selected points~\cite{Morokoff1994,Tuffin1996,JAMES1997180,Kocis1997}, thereby exchanging the exponential dependence on dimension for weaker guarantees on error.

In the many-body physics community, tensor networks (TN) such as matrix product states, projected entangled pair states, and the multiscale entanglement renormalization ansatz, have long been used to represent multivariable physical quantities, such as quantum states~\cite{ORUS2014117,Verstraete2008,Silvi2019,Orus2019,Felser2021} or Boltzmann densities~\cite{Chen2018,Pastori2019,Li2021}. Similar techniques (although mainly for more restricted classes of tensor networks, such as the CANDECOMP/PARAFAC decomposition~\cite{harshman1970foundations,Carroll1970,Beylkin2002,Beylkin2005,Kolda2009}, Hierarchical Tucker decomposition~\cite{Kolda2009,Grasedyck2010,GrasedyckHackbusch2011}, and tensor trains~\cite{GrasedyckHackbusch2011,Oseledets2011}) have appeared in the applied mathematics community as well, and have been used for high-dimensional function computation and approximation. In both cases, the idea is to represent a (often discretized) high-dimensional function as a connected network of low-dimensional tensors. This can reveal a non-trivial low-rank structure in the function, thereby 
achieving a great reduction in memory. Because TN are further embodied with a notion of approximation (through singular value decomposition of pairs of tensors) the use of such approximations can ameliorate, and in some cases remove, the exponential cost with respect to dimension, when computing with the TN.

For multivariable function computation with TN, we must first obtain a TN representation of the function.
The manner in which such representations is obtained usually involves problem-specific numerical computation. For instance, if the target function can be efficiently evaluated, then determination of the tensor parameters can be formulated as an optimization problem\cite{Schneider2012,Schneider2011,Novikov2016,OSELEDETS201070,DOLGOV2020,Beylkin2002,Beylkin2005,GarciaRipoll2021quantuminspired}. Another common scenario is when the function is implicitly defined from a minimization, in which case the tensor representation can be optimized by the variational principle~\cite{Evenbly2014,Vanderstraeten2016,Haghshenas2022}. When the function satisfies a differential equation with a known initial condition that is easily expressed as a TN, the TN representation can be propagated~\cite{Daley_2004,Iblisdir2007,Phien2015,Lin2021,LUBASCH2018587,Gourianov2022}. In all these cases, therefore, the determination of the tensor network representation of the function involves an associated, potentially large, computational cost. In addition, the tensor networks are generally chosen with a fixed network structure ahead of time (for example, a matrix product state or a tensor train) in order to make the determination of the representation feasible, even though such a structure may not be the most compact~\cite{Glasser2019,Levine2019,Haghshenas2022}.

Here, we introduce an alternative way to construct the TN representation of a function from its  arithmetic circuit. As both classical arithmetic circuits and quantum arithmetic circuits can be viewed as TNs, it is immediately clear that multivariable functions can be represented as TNs through these circuits. However, such circuits carry various disadvantages; for example, classical logic gates lead to extremely sparse tensors where the number of indices is proportional to the number of bits of precision~\cite{Verma2008,Ciesielski2020}, while quantum arithmetic circuits are constrained to unitary tensors, and may be sub-optimal for classical calculations~\cite{Wiebe2014,Haener2018}. Consequently, we introduce a tensor network circuit representation that takes advantage of both the ability to store floating point numbers as entries in the tensors, as well as the lack of unitary constraints. This leads to a concise construction that  avoids the auxiliary computation involved in obtaining the TN representation itself.

A byproduct of the arithmetic circuit TN form is that the structure of the tensor network is dictated by the circuit, rather than specified beforehand. Consequently, the TN structure that arises can have a general connectivity, and can in principle look quite different from those typically encountered in quantum many-body physics or applied mathematics settings. Recent developments in the exact and approximate contraction of tensor networks with more unstructured geometries are discussed, for example, in Refs.~\onlinecite{Gray2021hyperoptimized,Pan2020,Jermyn2020,Chubb2021,Gray2022hyperoptimized}. 

Using the arithmetic circuit TN representation we carry out high dimensional quadrature over the unit hypercube (in up to 50 dimensions) for multivariable polynomials and multivariable Gaussians.
The use of tensor networks in conjunction with high dimensional integration/summation is hardly new; it is one of the main applications of TNs see, for example, Refs.~\onlinecite{Vysotsky2021,OSELEDETS201070,DOLGOV2020,Xavier2022}. However, the availability of the circuit structure of the function leads to new insights into 
this problem. For example, for the multivariable polynomials, we find exact compressibility of the TN circuit in certain limits of the polynomial parameters. This can then be decoded into an exact integration rule (see Ref.~\onlinecite{Oseledets2013} for a related result). Furthermore, away from the exact point, we can relate the difficulty of approximation to the degree of nonlinearity (number of copy operations) in the circuit. In terms of practical performance, because Gaussian quadrature weights can be used in each dimension of the TN quadrature, we find that TN integration converges in accuracy orders of magnitudes more quickly than quasi-Monte Carlo integration for instances of the multivariable polynomial and Gaussian integrals. Finally, we finish with an artificial but instructive case where we  construct a function (based on the multi-scale entanglement renormalization ansatz~\cite{Evenbly2009_mera,Evenbly2014_mera}) that can be integrated efficiently when using knowledge of its internal TN circuit structure, but for which function evaluation is hard, thus making quasi-MC hard if the function is treated as a blackbox. 

\section{Tensor network arithmetic circuits}\label{sec:basic}

\subsection{Function tensor representation and
circuit composition of functions}

\label{sec:circuitcomposition}
We first introduce a tensor representation of individual single- and multi-variable functions and how to use the representation to compose complicated functions from simpler ones. We start with a single-variable function $f(x)$. Assuming for the time being that $x$ is a continuous variable, we introduce the continuous function tensor (function matrix) representation $F$, where the function values are stored in the elements $F_{x\alpha}$, specifically,
\begin{align}
F_{x0} &= 1 \notag\\
F_{x1} &= f(x)  \label{eq:functionrepr}
\end{align}
We refer to $x$ as the variable, and $\alpha$ as the control index/leg (the dimension of the control index is 2). We consider here scalar-valued functions, but vector-valued functions can be similarly defined. In numerical applications $x$ will typically be discretized, e.g. on a grid with $G$ points. The control index is used to perform arithmetic and other gate operations on the functions. 

\begin{figure}%
    \centering
    \includegraphics[width=0.1\linewidth]{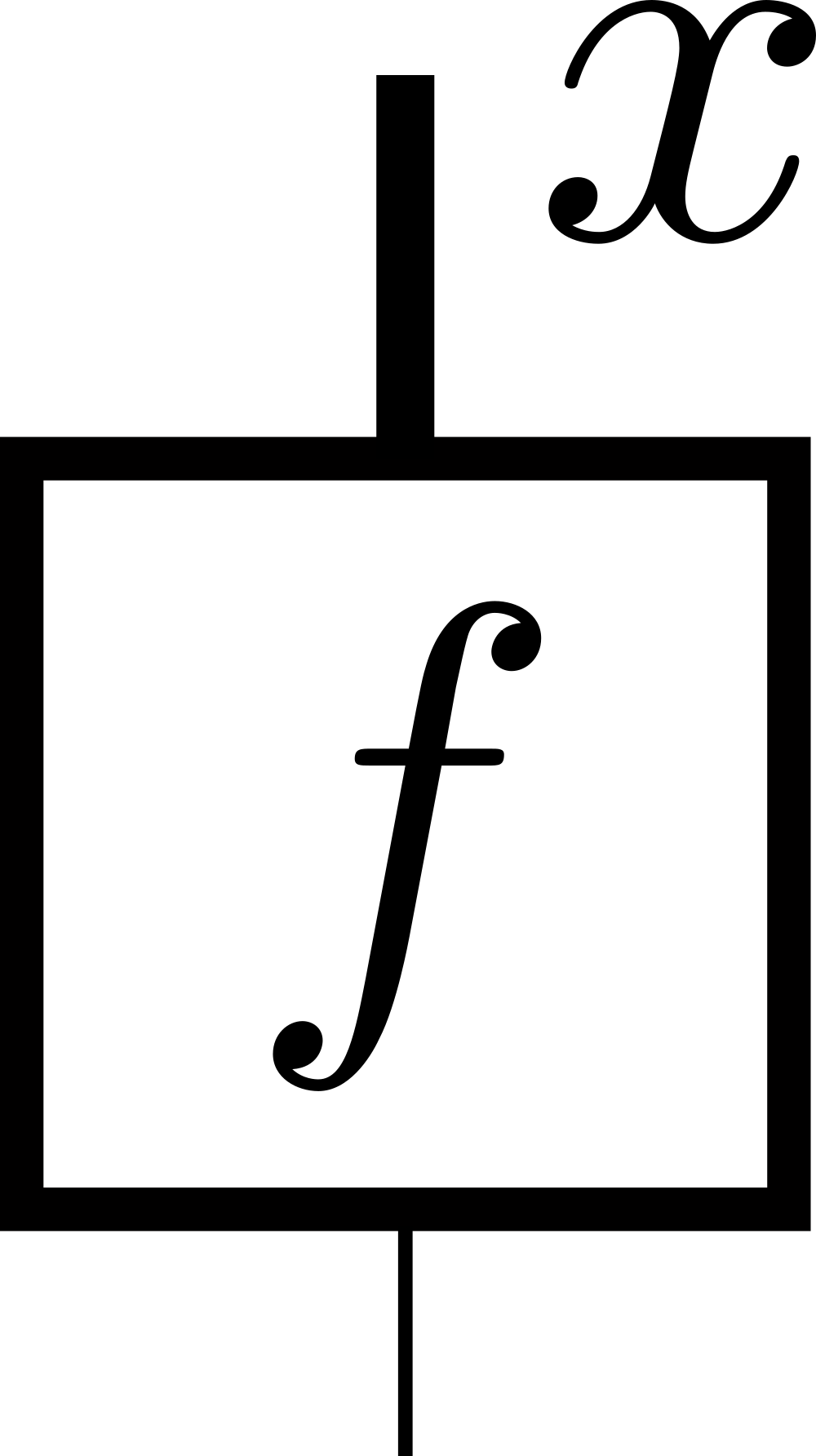}
    \caption{Tensor representation $F$ of the scalar function $f(x)$. The ``control'' index is represented by the thin leg on the bottom.}
    \label{fig:1var_scalar_fxn}
\end{figure}

It will be convenient to use the standard graphical notation of tensor networks. We show a diagram of $F$ in Fig.~\ref{fig:1var_scalar_fxn}. 
We use the convention that labelled lines represent indexed elements and unlabelled lines
between two tensors %(relevant below) 
are summed over for discrete indices/integrated over for continuous indices.

The control index can be thought of as a ``qubit'' index in the tensor product Hilbert space of functions, analogous to the qubit representation used in quantum mechanics; from Eq.~\ref{eq:functionrepr}, we see that $|1\rangle$ is associated with a basis function $f(x)$, and $|0\rangle$ is replaced by the scalar $1$. Then given a set of functions $\{ f_i (x_i) \}_{i=0...{N-1}}$, we represent a monomial of functions as
\begin{align}
    F^{[1]}_{x_0\alpha_0} F^{[2]}_{x_1\alpha_1} F^{[N-1]}_{x_{N-1}\alpha_{N-1}} = f_0(x_0)^{\alpha_0} f_1(x_1)^{\alpha_1} \ldots f_{N-1}(x_{N-1})^{\alpha_{N-1}}
\end{align}
A multivariable function $c(x_0, x_1, \ldots x_{N-1})$ in the product Hilbert space $\prod_{i} \{ 1, f(x_i) \}$ takes the form
\begin{align}
    c(x_0, \ldots, x_{N-1}) = \sum_{\{ \alpha_i\}} \prod_i C_{\alpha_0, \ldots, \alpha_{N-1}} F^{[i]}_{x_i\alpha_i}
\end{align}
Using the qubit analogy, $C_{\alpha_0, \ldots, \alpha_{N-1}}$ may be regarded as a wavefunction amplitude in the computational basis. Note that the above constructs the Hilbert space using single-variable functions, but we could also use more complex building blocks, e.g. a two-variable function tensor, $F_{xy\alpha} \leftrightarrow f(x, y)$. We also note that the use of a product Hilbert space to represent multivariable functions has been considered in other contexts, for instance, for length scale separation in the solution of partial differential equations, or for efficient function parameter storage. For relevant discussions, see e.g. Refs.~\onlinecite{GarciaRipoll2021quantuminspired,LUBASCH2018587,Gourianov2022}.

We use a tensor network to build a circuit representation of the function $c(x_0, \ldots, x_{N-1})$ using contractions of the function tensors either with themselves, or other fixed tensors. The simplest example corresponds to performing classical arithmetic and logic to combine single-variable functions, using control tensors. The addition tensor $(+)$ is a three index control tensor, with elements
\begin{align}\label{eqn:add}
(+)_{\alpha\beta\gamma}=\bigg\{\begin{array}{cc}
    1 & \alpha+\beta=\gamma \\
    0 & \text{otherwise}
\end{array}
\end{align}
Note that unlike usual binary arithmetic, the addition is not modulo 2. 
The multiplication tensor $(\times)$ is also a three index control tensor,
\begin{align}\label{eqn:scalar_mult}
(\times)_{\alpha\beta\gamma}=\bigg\{\begin{array}{cc}%#
    1 & \alpha=\beta=\gamma \\
    0 & \text{otherwise}
\end{array}
\end{align}

\begin{figure}%
    \centering
    \begin{subfigure}[b]{0.49\linewidth}
    \centering
    \includegraphics[width=0.5\textwidth]{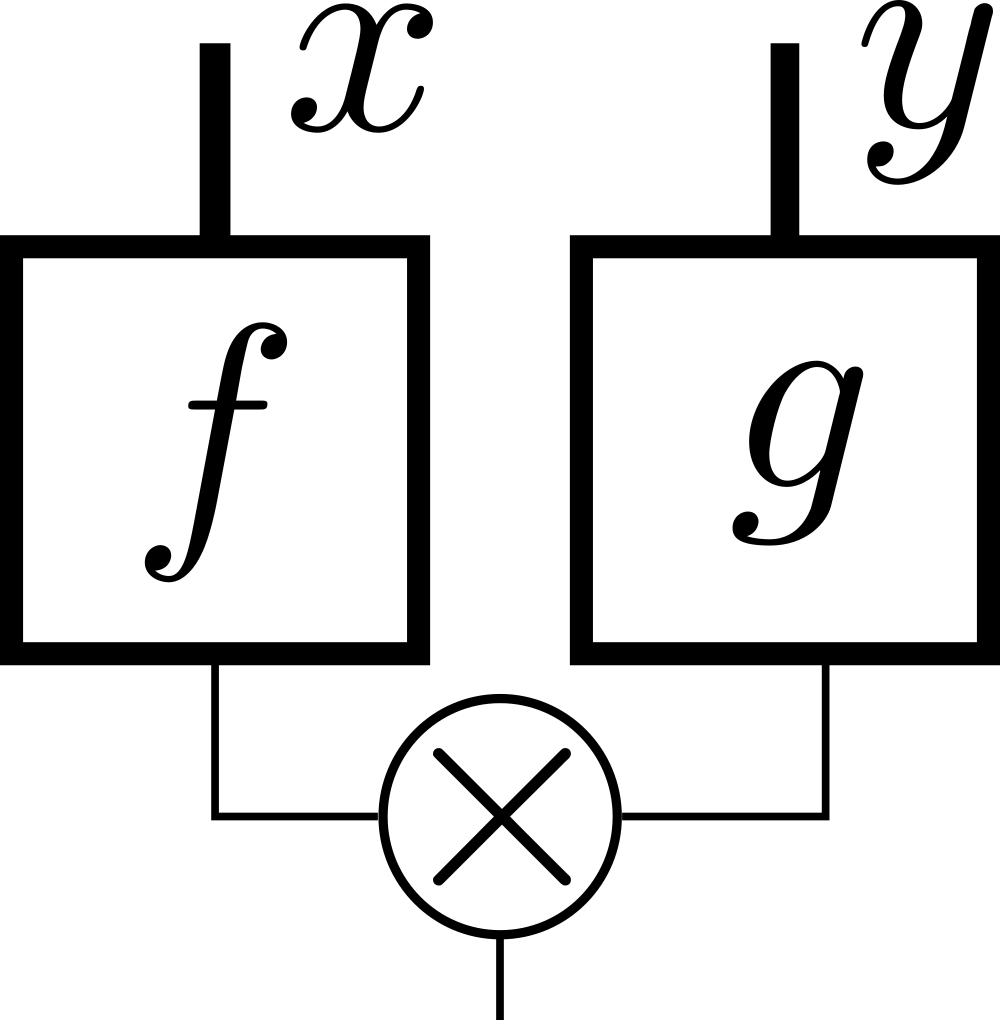}
    \caption{$f(x)g(y)$}
    \end{subfigure}
    \begin{subfigure}[b]{0.49\linewidth}
    \centering
    \includegraphics[width=0.5\textwidth]{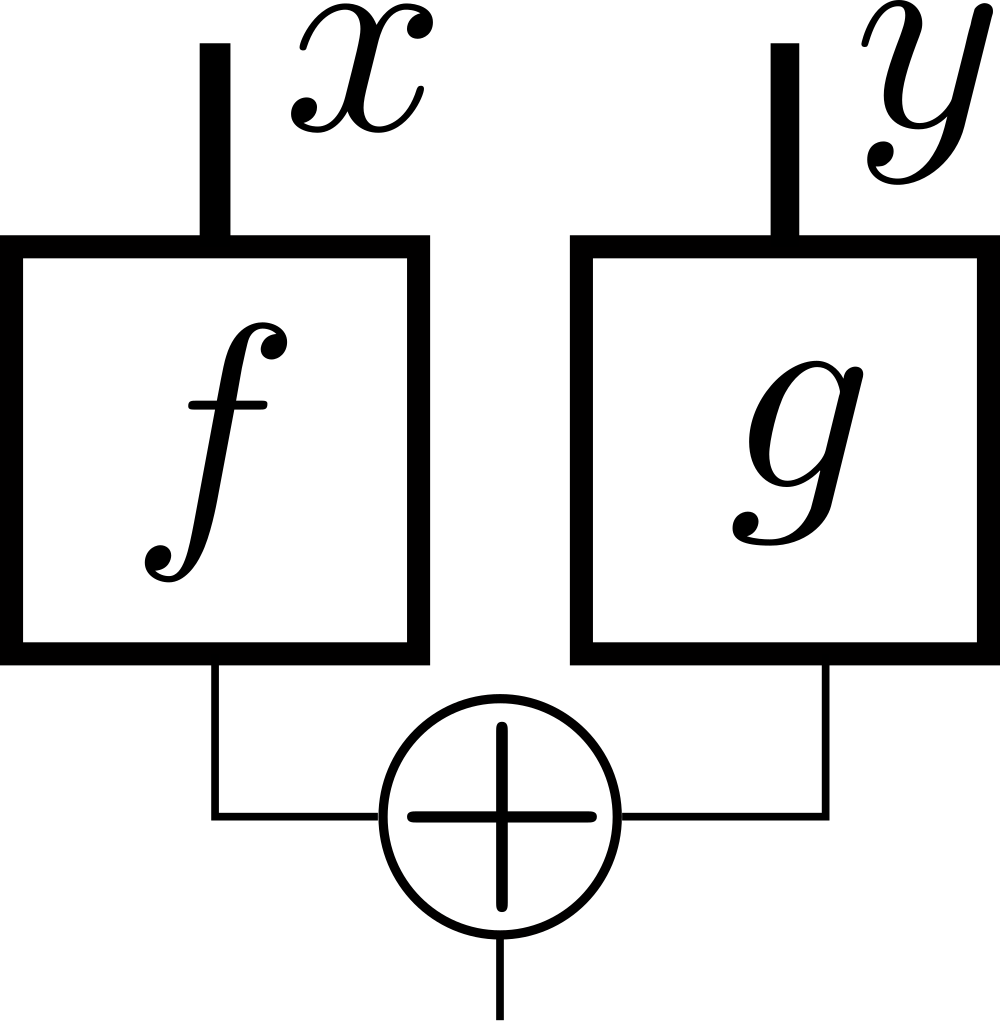}
    \caption{$f(x)+g(y)$}
    \end{subfigure}
    \caption{Arithmetic circuit tensor network representation of $f(x)g(y)$ and $f(x)+g(y)$. %Arrows pointing towards the $(+)$ and $(\times)$ tensors corresponds to $\alpha$ and $\beta$ indices. The arrow pointing away from $(+)$ and $(\times)$ tensor corresponds to $\gamma$ index.
    }%
    \label{fig:scalar_fxn_jux_add}%
\end{figure}

With such definitions we can perform classical arithmetic on functions. Contraction of two function tensors $F$, $G$ with the multiplication or addition tensors yields the higher dimensional objects (shown as a tensor network in Fig.~\ref{fig:scalar_fxn_jux_add}) corresponding to functions of two variables, i.e.
\begin{align*}
&F_{x\alpha}G_{y\beta}(\times)_{\alpha\beta\gamma} \leftrightarrow f(x) g(y)\\
&F_{x\alpha}G_{y\beta}(+)_{\alpha\beta\gamma} \leftrightarrow f(x) + g(y)
\end{align*}
where we have assumed summation over repeated indices. 

The above are simple examples of arithmetic circuit tensor networks.
Within these arithmetic circuits, the $0$ value of the control index is only needed to perform addition; if we know $f(x)$ only enters the circuit via multiplication, we can always choose to fix the control index $\alpha=1$ of the corresponding tensor $F$, and thus omit the control leg entirely. Similarly, scalars have no variable dependence and can thus be specified without their variable leg. 

A second circuit structure involves contraction between the variable legs of the function tensors, which corresponds to integrating a common variable between two functions. For example, given $F_{x\alpha} \leftrightarrow f(x)$ and $G_{k, x} \leftrightarrow g(k, x)$ (we have dropped the control index on $G$ since we are only using it for multiplication), then
\begin{align}
     F_{x\alpha} G_{kx} \leftrightarrow \int dx f(x) g(k, x)
\end{align}
where we have assumed integration over the repeated continuous index on the left. 

\begin{figure}%
    \centering
    \begin{subfigure}[b]{0.49\linewidth}
    \centering
    \includegraphics[width=0.5\textwidth]{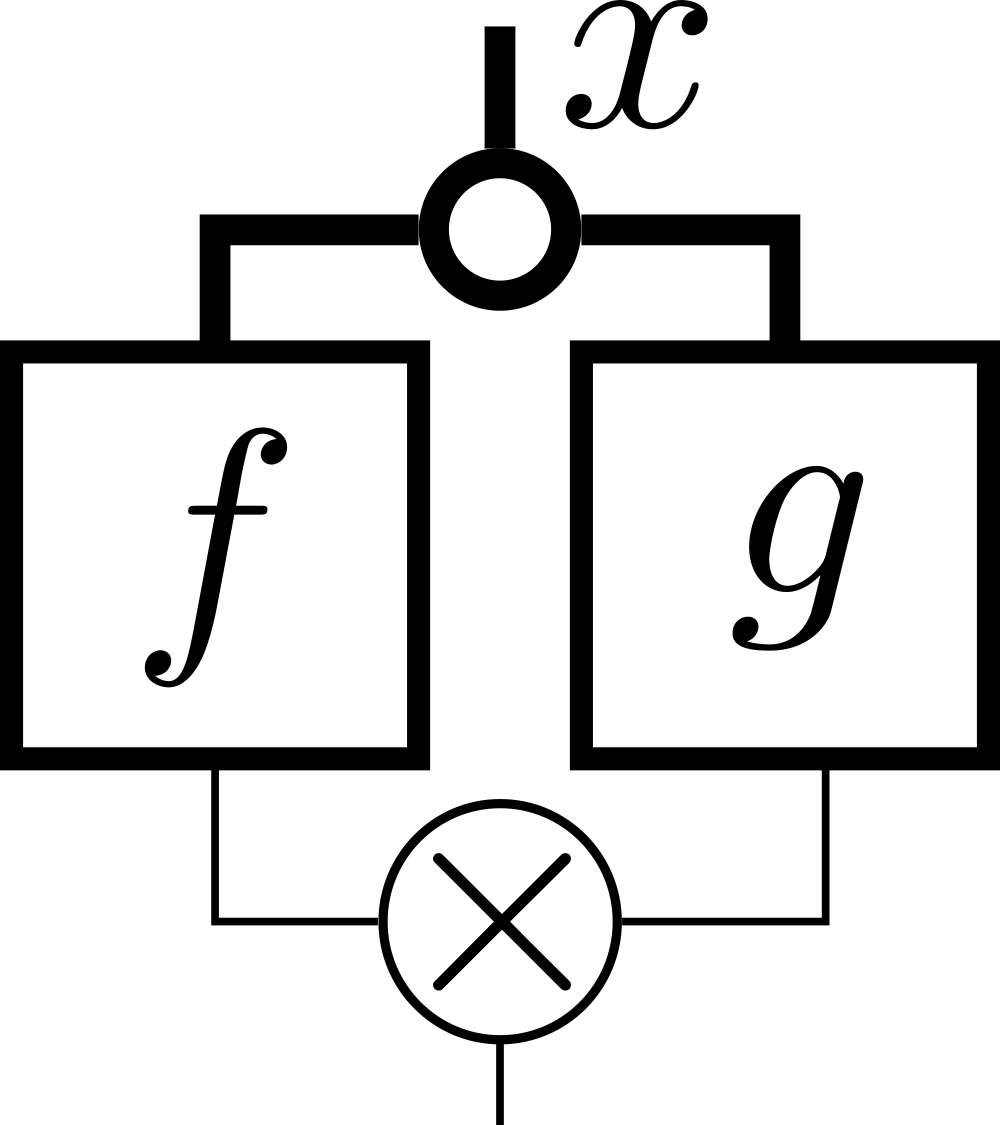}
    \caption{$f(x)g(x)$}
    \end{subfigure}
    \begin{subfigure}[b]{0.49\linewidth}
    \centering
    \includegraphics[width=0.5\textwidth]{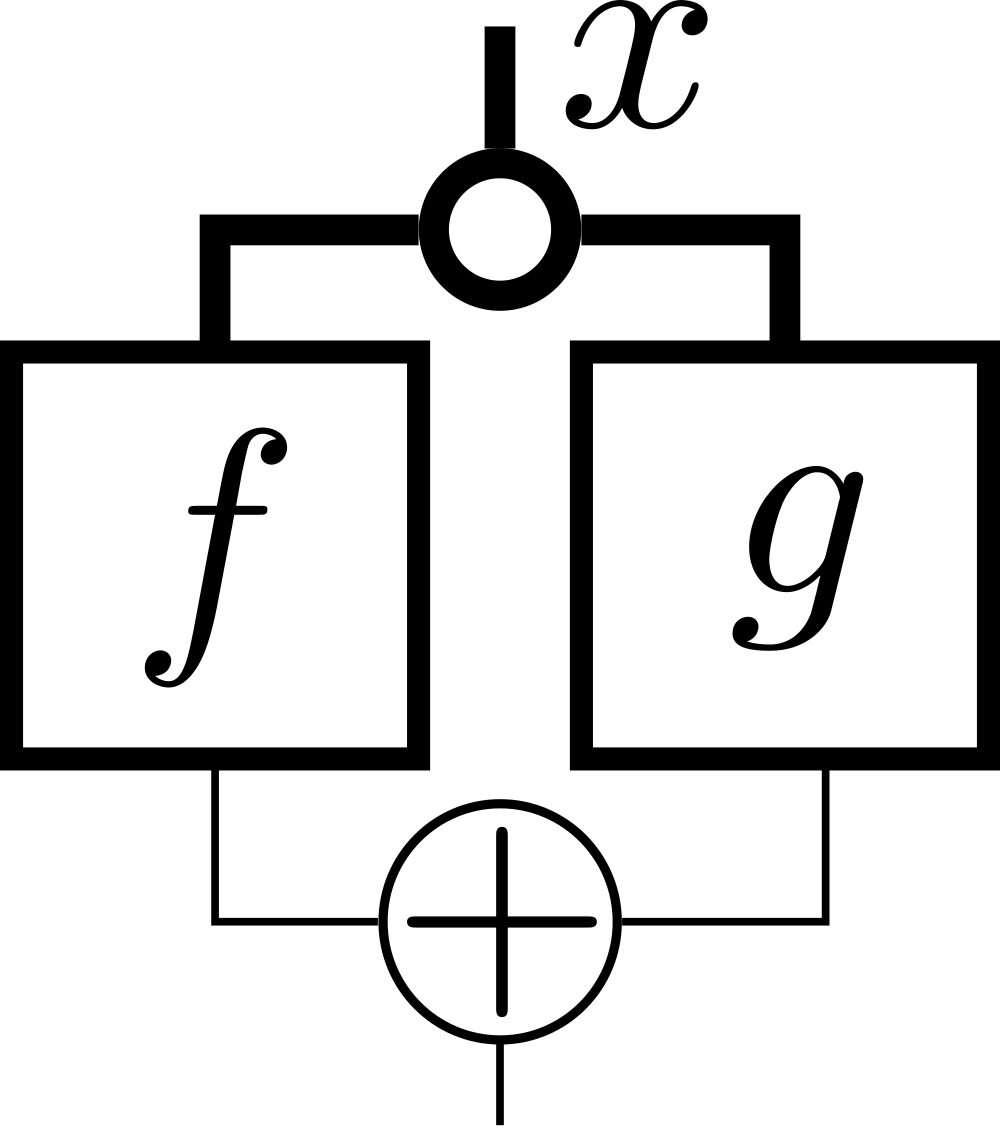}
    \caption{$f(x)+g(x)$}
    \end{subfigure}
    \caption{Arithmetic circuit tensor network representation of $f(x)g(x)$ and $f(x)+g(x)$. The circle represents the $\mathsf{COPY}$ tensor. }%
    \label{fig:scalar_fxn_jux_add_del}%
\end{figure}

Multivariable functions built from the above arithmetic operations will be multilinear in the underlying functions. To build non-linearity, we use multiple functions of the same variable. To do so, we 
define a tensor $\mathsf{COPY}$ on the continuous variable legs, with entries
\begin{align}
(\mathsf{COPY})_{xyz}=\bigg\{\begin{array}{cc}
    1 & x=y=z \\
    0 & \text{otherwise}
\end{array}    
\end{align}
Then, using the copy tensor we can
multiply or add two functions of the same variable, as shown in Fig.~\ref{fig:scalar_fxn_jux_add_del}, corresponding to the contractions
\begin{align*}
&(\mathsf{COPY})_{xyz}F_{y\alpha}G_{z\beta}(\times)_{\alpha\beta\gamma} \leftrightarrow f(x) g(x)\\
&(\mathsf{COPY})_{xyz}F_{y\alpha}G_{z\beta}(+)_{\alpha\beta\gamma} \leftrightarrow f(x) + g(x).
\end{align*}
Note that the contraction of the copy tensor is an example of a third type of circuit operation, namely, the contraction over the variable legs of functions with additional variable leg tensors; further examples of this kind are discussed in Sec.~\ref{sec:variable}.

\begin{figure}%
    \centering
    \includegraphics[width=0.25\linewidth]{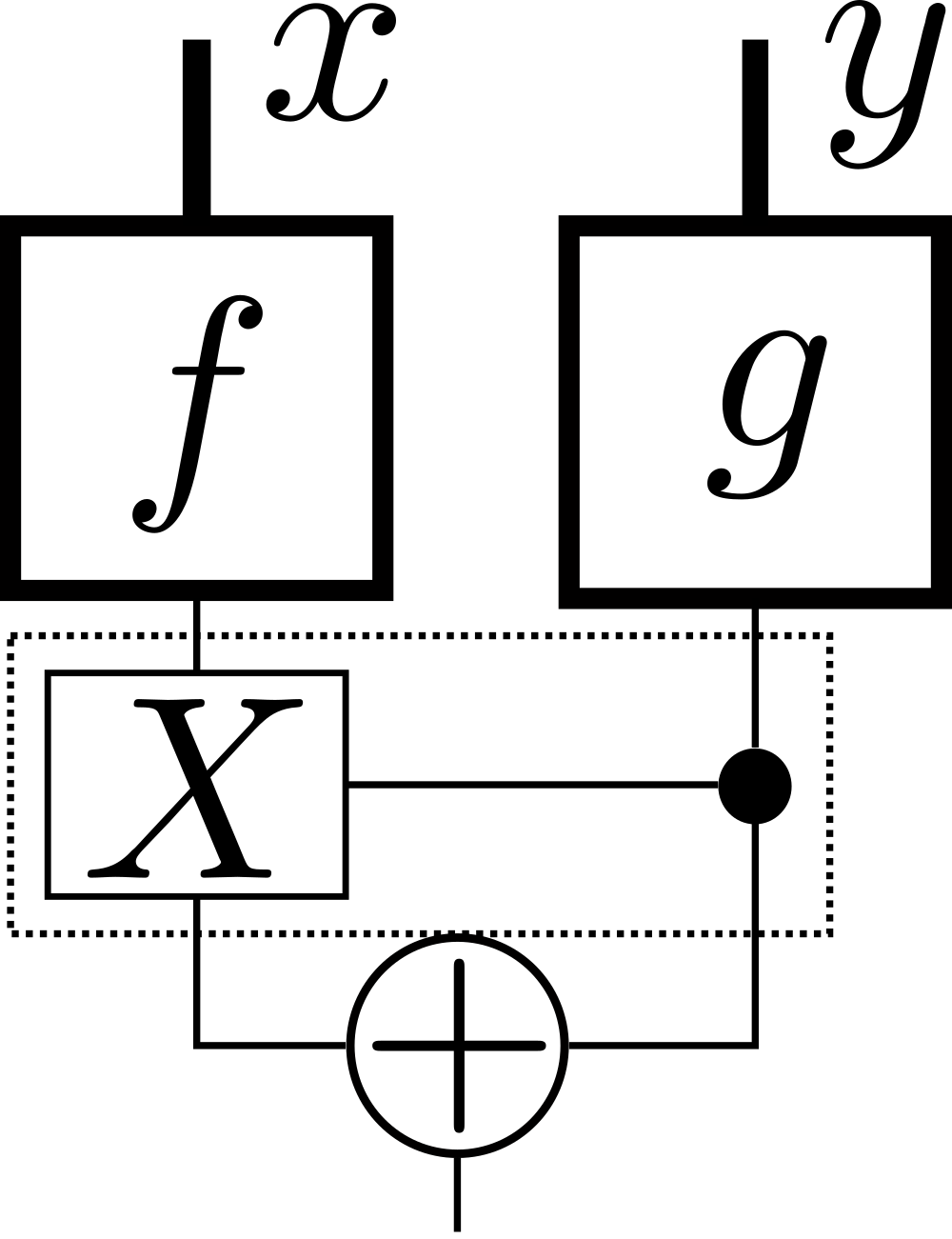}
    \caption{Arithmetic circuit tensor network representation of Eq.~(\ref{eqn:cnot}). The CNOT gate is represented by the dashed box. }
    \label{fig:quantumgate}
\end{figure}

Finally, we briefly mention that the control legs above are contracted with tensors that implement classical arithmetic/logic; contraction  generates a single function output. 
However, one can also apply more general binary operations, such as quantum logic gates where a single set of control inputs may map to multiple control outputs. For example, we can create a circuit using a quantum CNOT gate (see Fig.~\ref{fig:quantumgate}), giving 
\begin{align}\label{eqn:cnot}
    F_{x\alpha} G_{y\beta} (\mathsf{CNOT})_{\alpha\beta,\gamma\delta} (+)_{\gamma\delta\epsilon} \leftrightarrow f(x) + f(x)g(y).
\end{align}

\subsection{Variable circuits, transformations of variables, and other representations}

\label{sec:variable}

The $\mathsf{COPY}$ operation defined above generalizes to other TN circuit representations that use operations only on the variable legs of the function tensors. This can be viewed as using a different representation of functions as tensors, and we briefly examine the resulting circuits, here called variable TN circuits. In this case, the function $y(x)$ is represented by the function tensor with two continuous indices, 
\begin{align}
    \underline{F}_{y'x}= y(x) \ \text{if } y'=y(x)
\end{align}
where we have used the underline to emphasize that this is a variable TN circuit. One can e.g. perform arithmetic in this representation, using the addition and multiplication tensors
\begin{align}
    \underline{(+)}_{xyz} \ \text{if } z = x+y \notag\\
    \underline{(\times)}_{xyz} \ \text{if } z  = xy   \end{align}
In principle, one can represent all operators and all functions from circuits built up this way, however, in a discrete computation, the continuous variables must be discretized in some way. Then, because the precision of the result is stored in the tensor index, rather than the value of a tensor element, the number of bits of precision is limited by the size of the tensor indices. 

\begin{figure}%
    \centering
    \begin{subfigure}[b]{0.49\linewidth}
    \centering
    \includegraphics[width=0.7\textwidth]{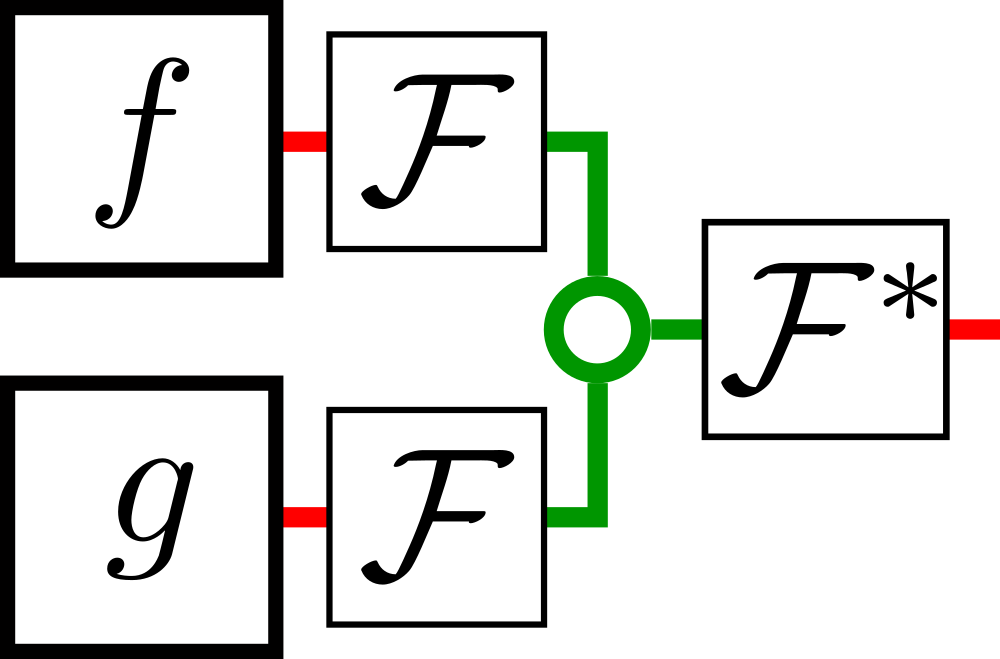}
    \caption{$\int dk \Tilde{f}(k) \Tilde{g}(k) e^{2\pi iky}$}
    \end{subfigure}
    \begin{subfigure}[b]{0.49\linewidth}
    \centering
    \includegraphics[width=0.5\textwidth]{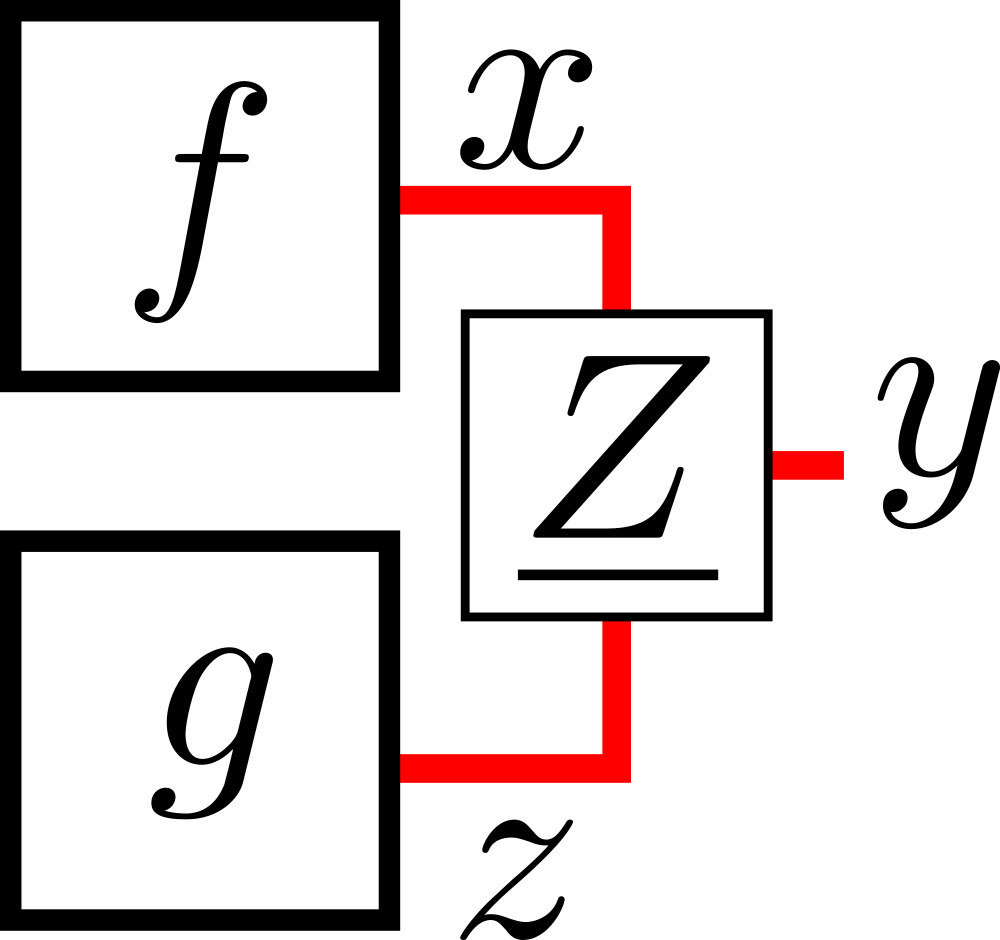}
    \caption{$\int dx f(x) g(y-x)$}
    \end{subfigure}
    \caption{Representation of the convolution theorem. The red legs represent position space variables $x,y,z$, and the green legs represent momentum space variable $k$. $\mathcal{F}$ represents the Fourier transform, and $\mathcal{F}$ and $\underline{Z}$ are both variable tensors.}%
    \label{fig:convolution}%
\end{figure}

In a circuit that contains both function tensors and variable tensors, the variable part of circuit can be seen as performing a transformation or mapping of variables before entering the function representation. For example, consider the convolution $\int dx f(x) g(y-x)$. This can be written as the combination of a function and a variable circuit, as the network
\begin{align}
    F_{x\alpha} G_{z\beta} (\times)_{\alpha\beta\gamma} \underline{Z}_{xyz} 
\end{align}
with 
\begin{align*}
\underline{Z}_{xyz}=\bigg\{\begin{array}{cc}
    1 & z=y-x \\
    0 & \text{otherwise}
\end{array}.
\end{align*}
The convolution theorem $\int dk \Tilde{f}(k) \Tilde{g}(k) e^{2\pi iky} = \int dx f(x) g(y-x)$ can then be expressed as the diagram shown in Fig.~\ref{fig:convolution}.

\subsection{Arithmetic circuit tensor networks and integration}

\begin{figure}%
    \centering
    \begin{subfigure}[b]{0.24\linewidth}
    \centering
    \includegraphics[width=0.9\textwidth]{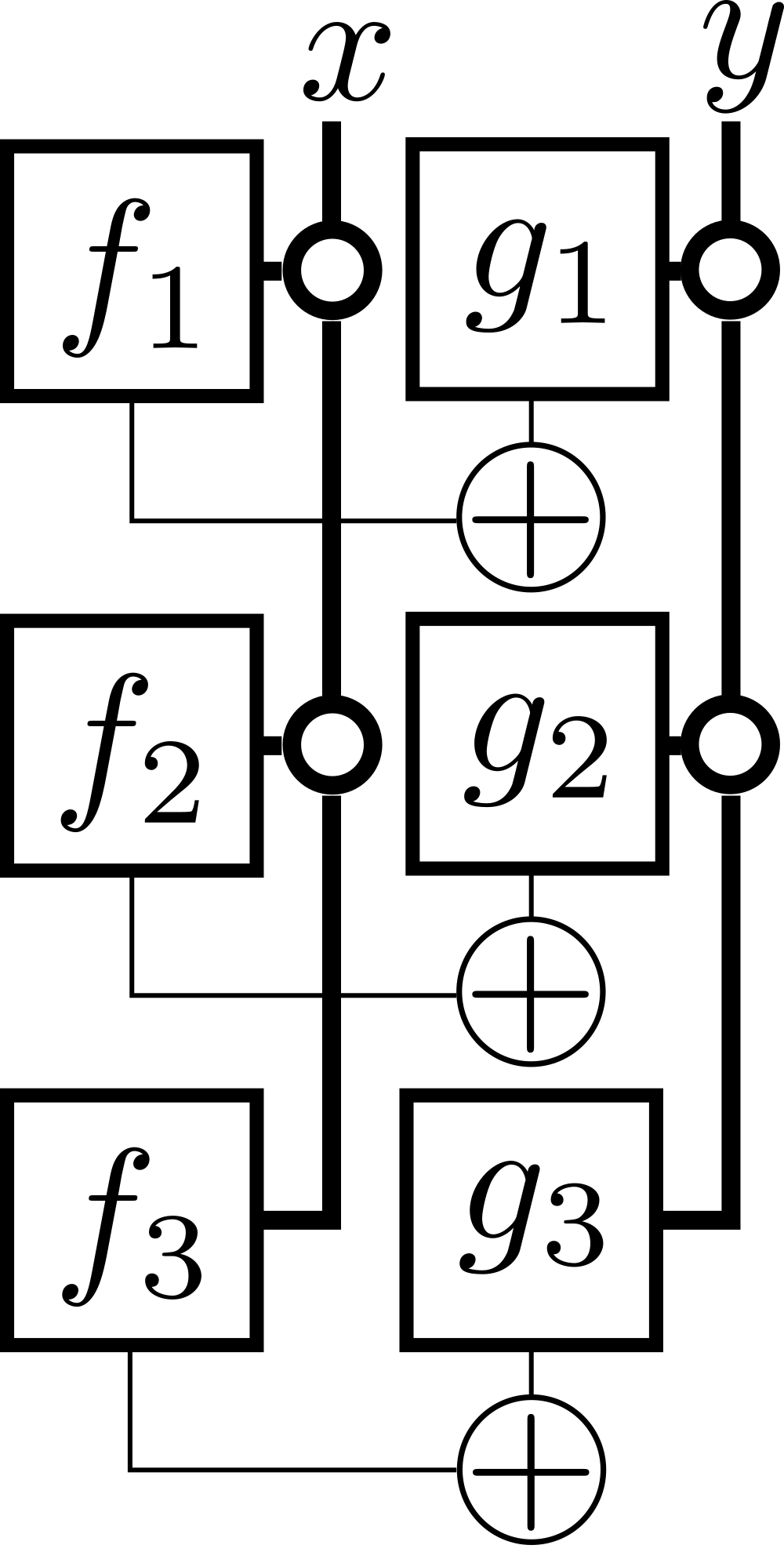}
    \caption{}
    \label{fig:generaltensornetworkcircuit}
    \end{subfigure}
    \begin{subfigure}[b]{0.24\linewidth}
    \centering
    \includegraphics[width=\textwidth]{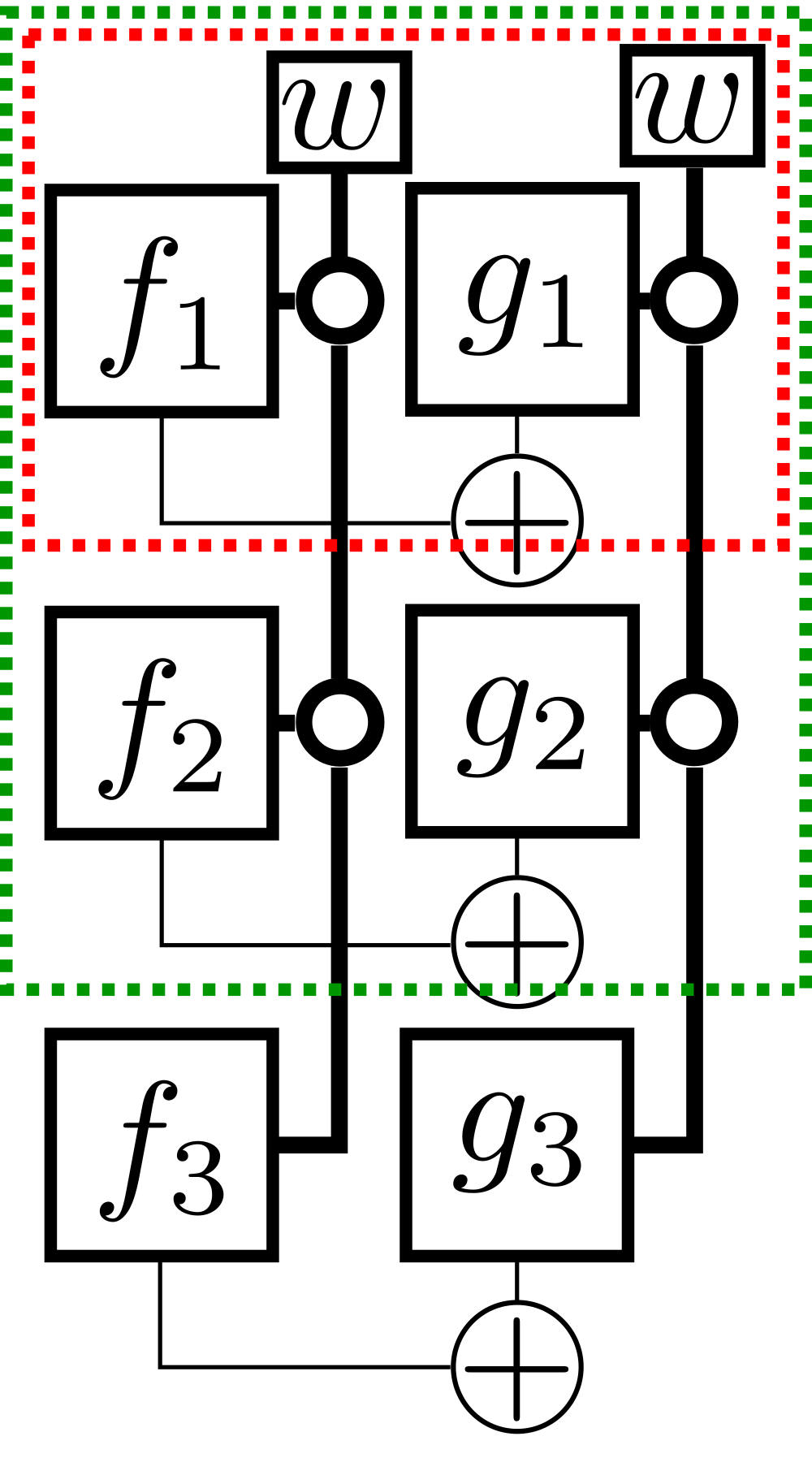}
    \caption{}
    \label{fig:contractionorder1}
    \end{subfigure}
    \begin{subfigure}[b]{0.24\linewidth}
    \centering
    \includegraphics[width=\textwidth]{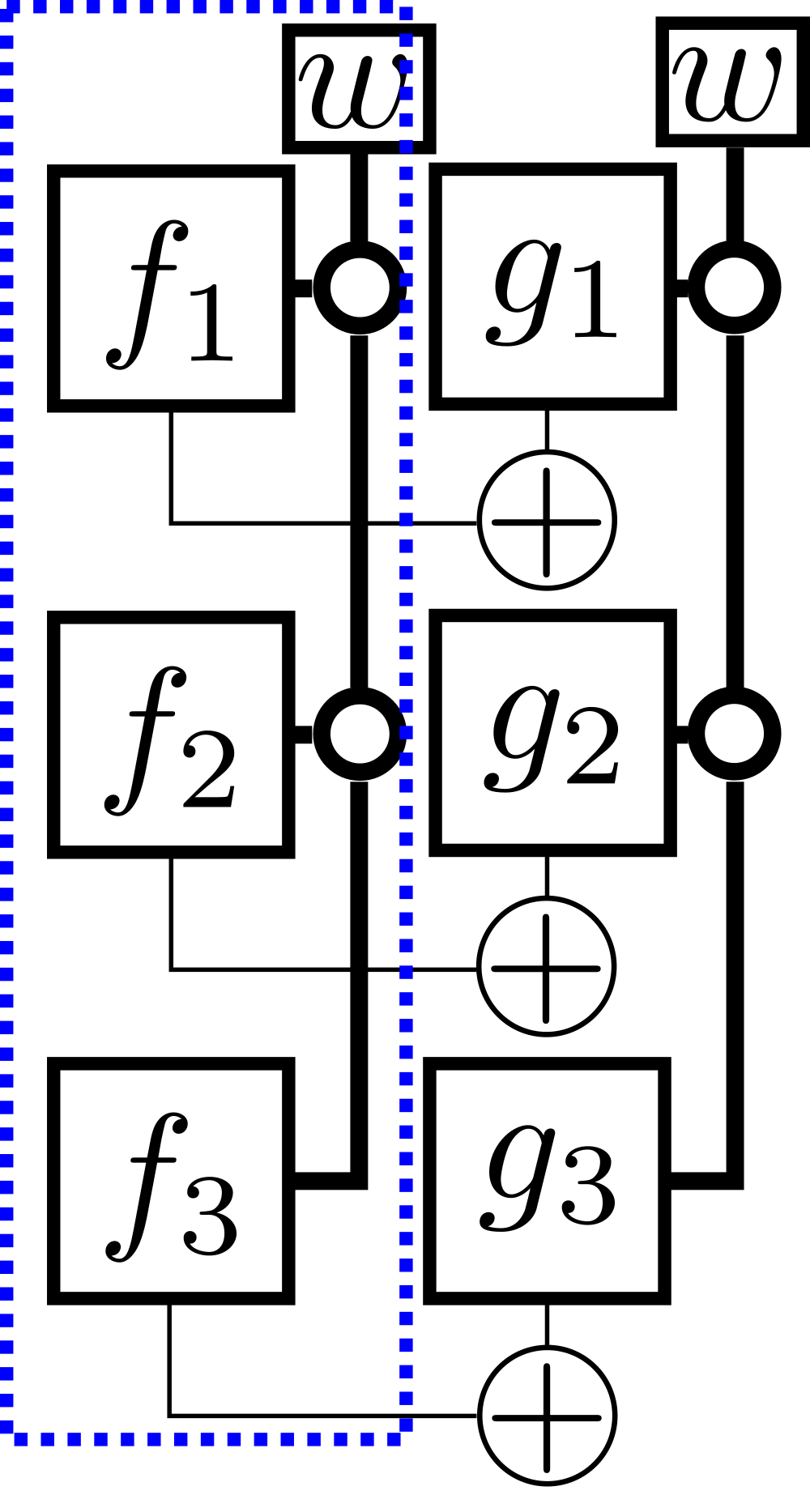}
    \caption{}
    \label{fig:contractionorder2}
    \end{subfigure}
    \begin{subfigure}[b]{0.24\linewidth}
    \centering
    \includegraphics[width=0.9\textwidth]{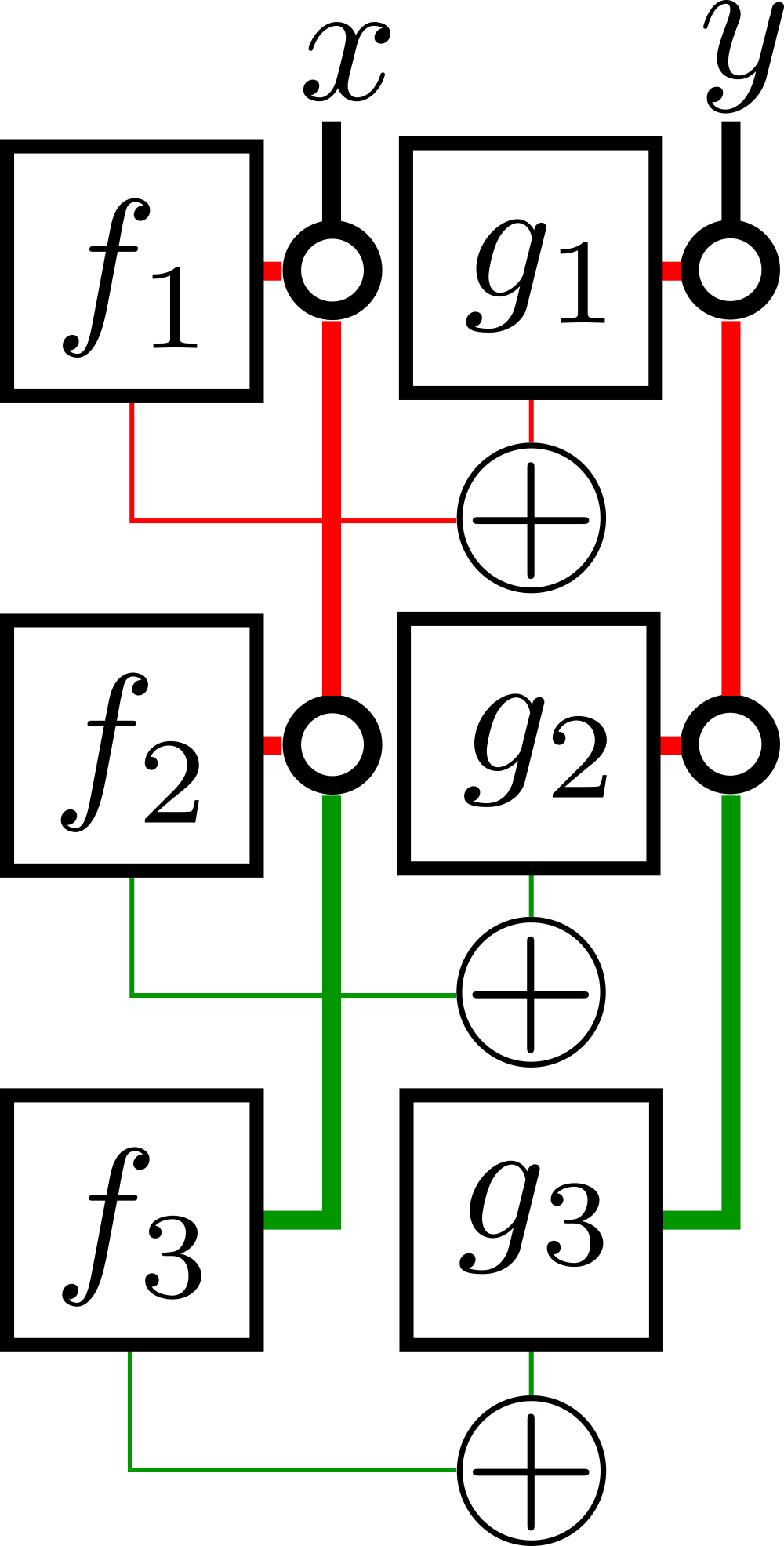}
    \caption{}
    \label{fig:loop}
    \end{subfigure}
    \caption{Representation of the function $\prod_{i=1}^3(f_i(x)+g_i(y))$. (a) TN representation of the function. (b) Integration of the function with one contraction scheme: the tensors in the red box are contracted first, then contracted with the tensors in the green box, then contracted with the remaining tensors. (c) Integration of the function with another contraction scheme: the tensors in the blue box are contracted first, then contracted with the remaining tensors. (d) Representation of the function showing loops (red, green lines) }%
\end{figure}

The composition of the above elements clearly allows us to construct a general multivariable function by operations on function tensors, or on variable tensors. This yields the arithmetic circuit tensor network representation of the multivariable function. A simple example for the function $\prod_{i=1}^3(f_i(x)+g_i(y))$ is shown in Fig.~\ref{fig:generaltensornetworkcircuit}. 

Given the TN representation, it is trivial to define integration over the input variables. Assuming a product quadrature for each variable $x_0, \ldots x_{N-1}$, then each continuous variable is discretized $x_i[p] \leftrightarrow x_i$, where $[p]$ denotes the $p$-th grid point. One can also introduce quadrature weights $w[p]$, giving
\begin{align}
    F^{\int}_{\alpha} = \sum_p w_p F_{p\alpha} \leftrightarrow \int dx f(x) \label{eq:gridweight}
\end{align}
as shown in Figs.~\ref{fig:contractionorder1},~\ref{fig:contractionorder2}. 

\subsection{Contracting the arithmetic circuit}

To obtain the value of the function from the arithmetic circuit one needs to contract the tensors. Assuming a scalar output, then 
\begin{align}
    F_\alpha = \sum_{\{ e\}/\alpha} \prod_i f^{[i]}_{ \{ e\}} 
\end{align}
where $e$ are all the contracted indices, which includes both variable and logical indices in the general case. 

There are many techniques for contracting tensor networks. Here we will view exact and approximate tensor network contraction as mainly black box algorithms, thus we do not discuss the details and only give a brief idea of the fundamentals. Further details of techniques for exact contraction are described in Ref.~\onlinecite{Bravyi2008,Fried2018,Levin2007,Evenbly2015,Evenbly2017,Bal2017,Yang2017,Ran2020,Schindler_2020,Gray2021hyperoptimized}. Additional information on approximate contraction can be found in Refs.~\onlinecite{Pan2020,Jermyn2020,Chubb2021,Denny_2011,Daley_2004,Verstraete2008,ORUS2014117,Orus2019,Silvi2019,Ran2020,Vlaar2021,Felser2021,Lin2021,Haghshenas2022,Gray2022hyperoptimized}. 

Although a typical arithmetic circuit has a flow from input values to output values, the fact that the arithmetic circuit TN encodes all output values for all inputs simultaneously, removes the directionality of the circuit. In particular, the contractions in the circuit can be evaluated in any order (see Figs.~\ref{fig:contractionorder1},~\ref{fig:contractionorder2}), interchanging the order of summations and products. This can lead to drastic changes in the complexity of evaluation. Heuristic techniques exist to search and find good orders for exact contraction which have been applied to tensor network and quantum circuit tensor network contraction problems (see Ref.~\onlinecite{Gogate2012,Pfeifer2014,Strasser2017,Akhremtsev2017,Boixo2017,Kourtis2019,Dudek2019,Gray2021hyperoptimized}).

Tensor networks without loops (i.e. trees) can be contracted exactly with a cost linear in the number of tensors. For arbitrary connectivity, e.g. with loops, the exact tensor network contraction scales exponentially with the number of tensors. However, tensor networks can also be contracted approximately with an approximation error. Such approximate contraction is widely used in tensor network applications in the simulation of quantum systems, and is closely related to low-rank matrix factorization~\cite{Grasedyck2013,Kumar2017,KHOROMSKIJ20121}. During the contraction of the tensor network, tensors will be generated which share an increased number of legs with neighbouring tensors (see Fig.~\ref{fig:contractionorder1}, \ref{fig:contractionorder2}); denote the combined dimension of the legs $D$. In approximate contraction, we use projectors of dimension $D\times \chi$ to project the shared leg of dimension $D$ down to a shared leg of dimension $\chi$, thus controlling the cost of further operations. In practice, these isometric matrices are usually determined from an singular value decomposition (SVD).

A simple example of a compression algorithm is the ``boundary'' compression algorithm for a regular 2D tensor network. Here, rows of tensors (so-called matrix product states and matrix product operators) are contracted together, to form matrix product states with shared bonds; then a series of SVD's are applied to reduce the bonds to dimension $\chi$ (see Refs.~\onlinecite{Daley_2004,Verstraete2008,ORUS2014117,Orus2019,Silvi2019,Ran2020,Vlaar2021,Felser2021,Lin2021,Haghshenas2022,Gray2022hyperoptimized}). We will deploy the boundary contraction algorithm as the approximation contraction algorithm below. However, just as for exact contraction, the choice of order of when to contract and compress can greatly affect cost and accuracy, and similar to exact contraction, there are heuristics to choose an optimized order of contraction and compression. In this work, we do not encounter sufficiently complicated network structures to use these more sophisticated strategies, but the interested reader is referred to Ref.~\onlinecite{Gray2022hyperoptimized}. 

Approximate contraction is critical for applications of arithmetic tensor networks, because it allows subclasses of arithmetic tensor networks to be executed for less than the brute force exponential cost in the number of tensors. This is useful in defining classes of computational problems where the
curse of dimensionality is circumvented, as we now examine in our application to multidimensional integration below.

\section{Applications to integration}

We now apply the arithmetic tensor network formalism to the problem of multivariable integration. We start with some basic intuition about complexity, then proceed to progressively more complicated examples of multivariable polynomial and multidimensional Gaussian integration in the hypercube,
illustrating the power of the method versus another high dimensional technique, namely quasi-Monte Carlo integration. We finish with a specific circuit function with a complexity-theoretic guarantee of hardness with respect to sampling its values (and thus integration by direct application of quasi-Monte Carlo methods) but which can be efficiently integrated if one uses its tensor network structure. 

\subsection{Intuition regarding complexity}

It is well known that multi-variable functions admitting a separation of variables are easy to integrate over a separable range. The simplest examples are
\begin{align}
&\int_{\Omega^2}dxdyf(x)g(y)=\left(\int_{\Omega}dxf(x)\right)\left(\int_{\Omega}dyg(y)\right)\\
&\int_{\Omega^2}dxdy\left(f(x)+g(y)\right)=\Omega\int_{\Omega}dxf(x)+\Omega\int_{\Omega}dyg(y)
\end{align}
where $\Omega$ is the integration range in each variable. In TN language, the separability in the above equations corresponds to the tree structure in the TN diagrams shown in Fig.~\ref{fig:scalar_fxn_jux_add}. As discussed, such loop-free tensor networks are easy to contract. In the case of integration, one integrates over the variable legs (numerically, one sums over the discrete variable index with grid weights as in Eq.~(\ref{eq:gridweight})) and then one repeatedly contracts child tensors into their parents, never creating any shared legs. Note that loop free structures constitute a larger class of functions than separable functions. For example, assuming all operations are adds and multiplies, then the following function
\begin{align}
     ((f_1(x_1) f_2(x_2) + f_3(x_3)) f_4(x_4) + f_5(x_5))f_6(x_6) \ldots
\end{align}
where the nested parentheses reflects a binary tree structure, is easily integrated.

In contrast, the following multivariable integration, illustrated in Fig.~\ref{fig:loop},
\begin{align}
\int_{\Omega^2}dxdy(f_1(x)+g_1(y))(f_2(x)+g_2(y))(f_3(x)+g_3(y))
\end{align}
does not generate a loop-free arithmetic tensor network. In the corresponding TN diagram the \textsf{COPY} tensor for each variable results in the red and green loops, arising from a non-linear dependence of the function on a variable. 
The contraction of tensors that are part of 2 loops with their neighbors can lead to tensors with more legs/larger size. Similarly the use of quantum gates can result in loops. Thus in our arithmetic TN circuits the use of \textsf{COPY} tensors and quantum gates make the resulting tensor networks increasingly hard to contract, and the resulting functions harder to integrate.

\section{Multivariable polynomial of functions}\label{section:pol}

As an instructive example, consider the integral
\begin{align}
    Z = \int_{\Omega^N} dx_1 \ldots dx_N f(x_1,\ldots x_N)
\end{align}
with $f$ a polynomial of the form
\begin{align}\label{eqn:pol}
f(x_1,...,x_N)=\prod_{i=1}^kp_i(x_1,...,x_N)
\end{align}
with
\begin{align}\label{eqn:pol_fac}
p_i(x_1,...,x_N)=\sum_{j=1}^Nq_{ji}(x_j)
\end{align}
and where $q_{ji}(x_j)$ are single-variable functions. This functional form was considered in Ref.~\onlinecite{Fu2012}, which showed that integration over the domain $[0,1]^N$ is NP-hard for arbitrary functions $q_{ji}$. We first construct an arithmetic circuit TN representation and then proceed to investigate the complexity of integration for different choices of the functions $q_{ji}(x_j)$.

\subsection{Arithmetic tensor network representation}

\begin{figure}%
    \centering
    \begin{subfigure}[b]{0.49\linewidth}
    \centering
    \includegraphics[width=0.8\textwidth]{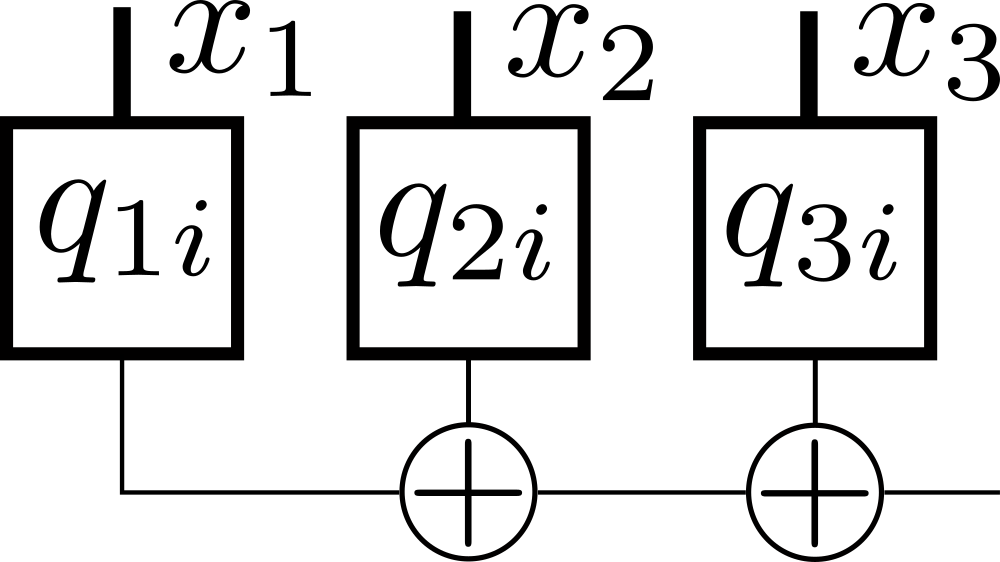}
    \caption{$p_i(x_1,\ldots,x_N)$}
    \label{fig:pol_fac}
    \end{subfigure}
    \begin{subfigure}[b]{0.49\linewidth}
    \centering
    \includegraphics[width=0.8\textwidth]{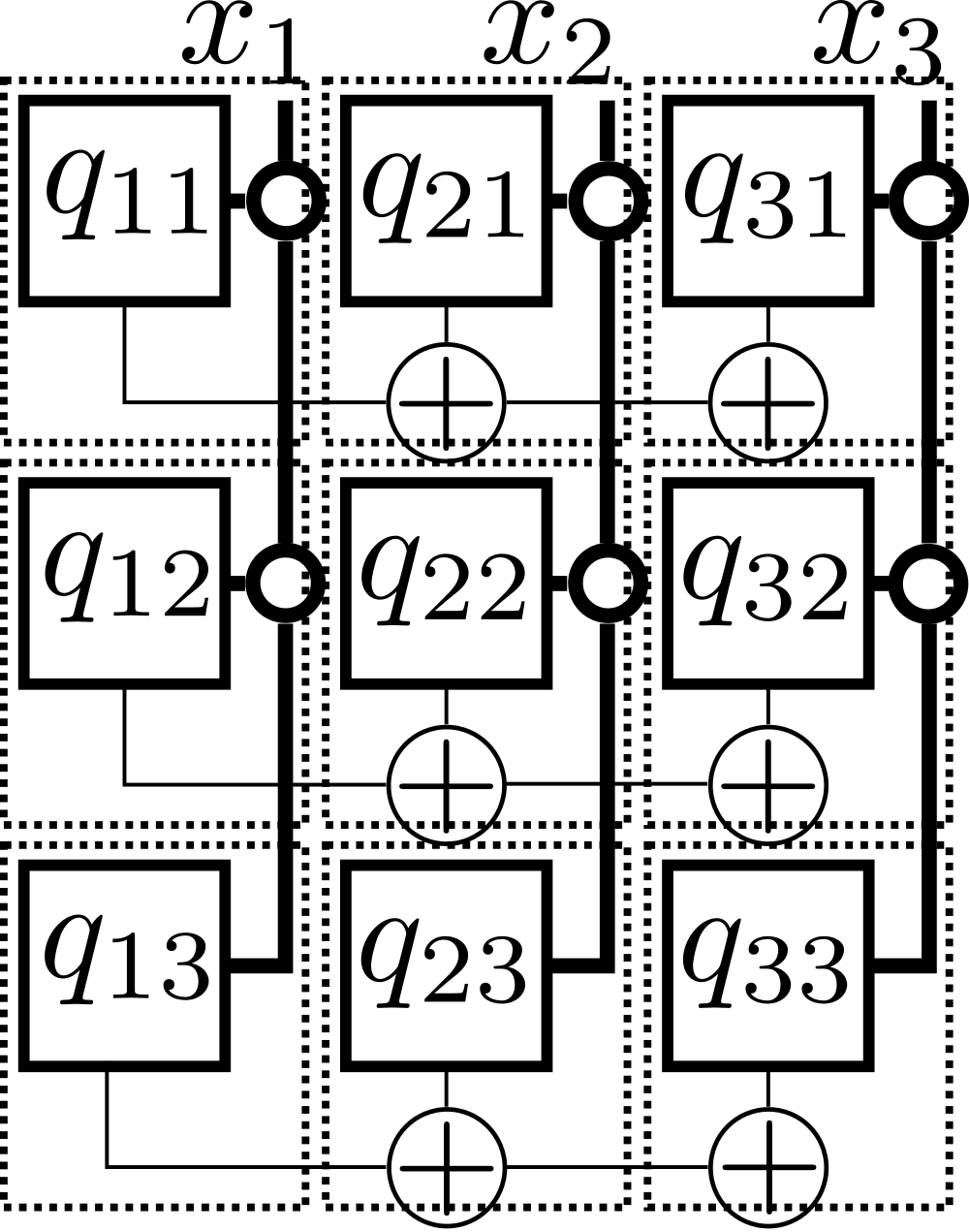}
    \caption{$f(x_1,\ldots,x_N)$}
    \label{fig:pol_full}
    \end{subfigure}
    \caption{Representation of $p_i(x_1,\ldots,x_N)$ and $f(x_1,\ldots,x_N)$ for $N=3$. As discussed in section~\ref{sec:basic}, we have omitted the open control leg for each $p_i(x_1,\ldots,x_N)$ since they enter $f(x_1,\ldots,x_N)$ only through multiplication. }%
\end{figure}

Since the total function is explicitly a product of factors, it is natural to construct a representation for each factor and then multiply them together. 
Consider a single factor $p_i(x_1,\ldots,x_N)$. This is a sum of terms, and we can use a circuit (shown in Fig.~\ref{fig:pol_fac}) to implement Eq.~(\ref{eqn:pol_fac}), where the horizontal bonds of dimension 2 are a manifestation of the low rank structure of $p(x_1,\ldots,x_N)$. The arithmetic circuit for the full tensor network is shown in Fig.~\ref{fig:pol_full}, where the TN for each factor $p_i(x_1,\ldots,x_N)$ becomes a row (with index $i$) and the set of $q_{ji}$ connected by \textsf{COPY} tensors for a given variable $x_j$ becomes a column (indexed by $j$). The \textsf{COPY} tensors make the entire network loopy; as each variable is copied $k$ times, $k$ is a measure of the loopiness of the network. 
If  we first contract the tensors in each dashed box, the resulting circuit has a regular 2D structure, known as a projected entangled pair state (PEPS)~\cite{ORUS2014117} structure. 

\subsection{Exact compressibility and an identity}\label{sec:pol_exact}

\begin{figure}%
    \centering
    \begin{subfigure}[b]{0.32\linewidth}
    \centering
    \includegraphics[width=0.9\textwidth]{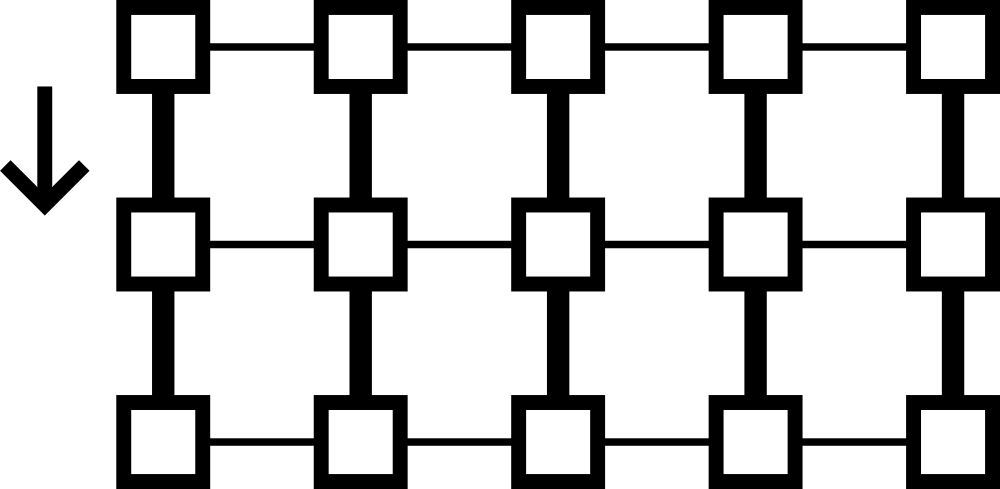}
    \caption{}
    \end{subfigure}
    \begin{subfigure}[b]{0.32\linewidth}
    \centering
    \includegraphics[width=0.9\textwidth]{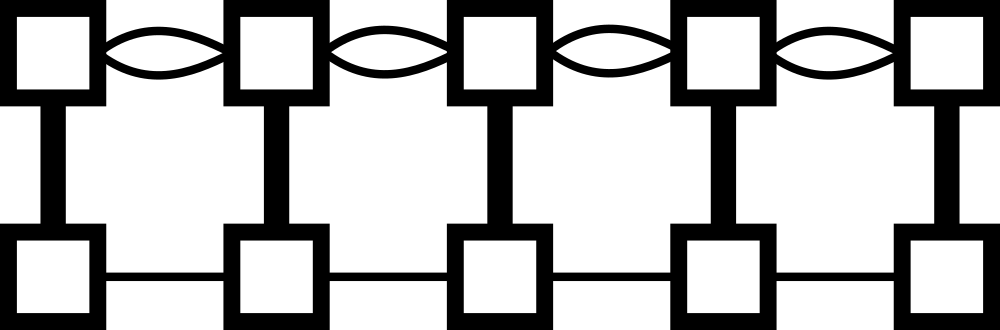}
    \caption{}
    \end{subfigure}
    \begin{subfigure}[b]{0.32\linewidth}
    \centering
    \includegraphics[width=0.9\textwidth]{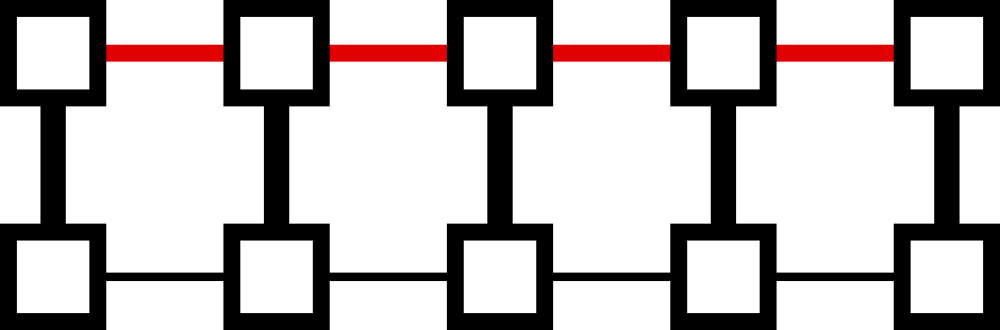}
    \caption{}
    \end{subfigure}
    \caption{Boundary row compression of a PEPS with 3 rows, 5 columns. (a) shows the PEPS, where each tensor $T_{i,j}$ ($i/j$ index row/columns) corresponds to a dashed box in Fig.~\ref{fig:pol_full}. (b) For each column $j$, contract $T_{1,j}$ with $T_{2,j}$. (c) For each pair of tensors $T_{2,j}$, $T_{2,j+1}$ connected by multi-bonds, compress the multi-bonds to a maximum bond dimension $\chi$ (red). }%
    \label{fig:boundary}%
\end{figure}

\begin{figure*}%
    \centering
    \begin{subfigure}[b]{0.3\linewidth}
    \centering
    \includegraphics[width=\textwidth]{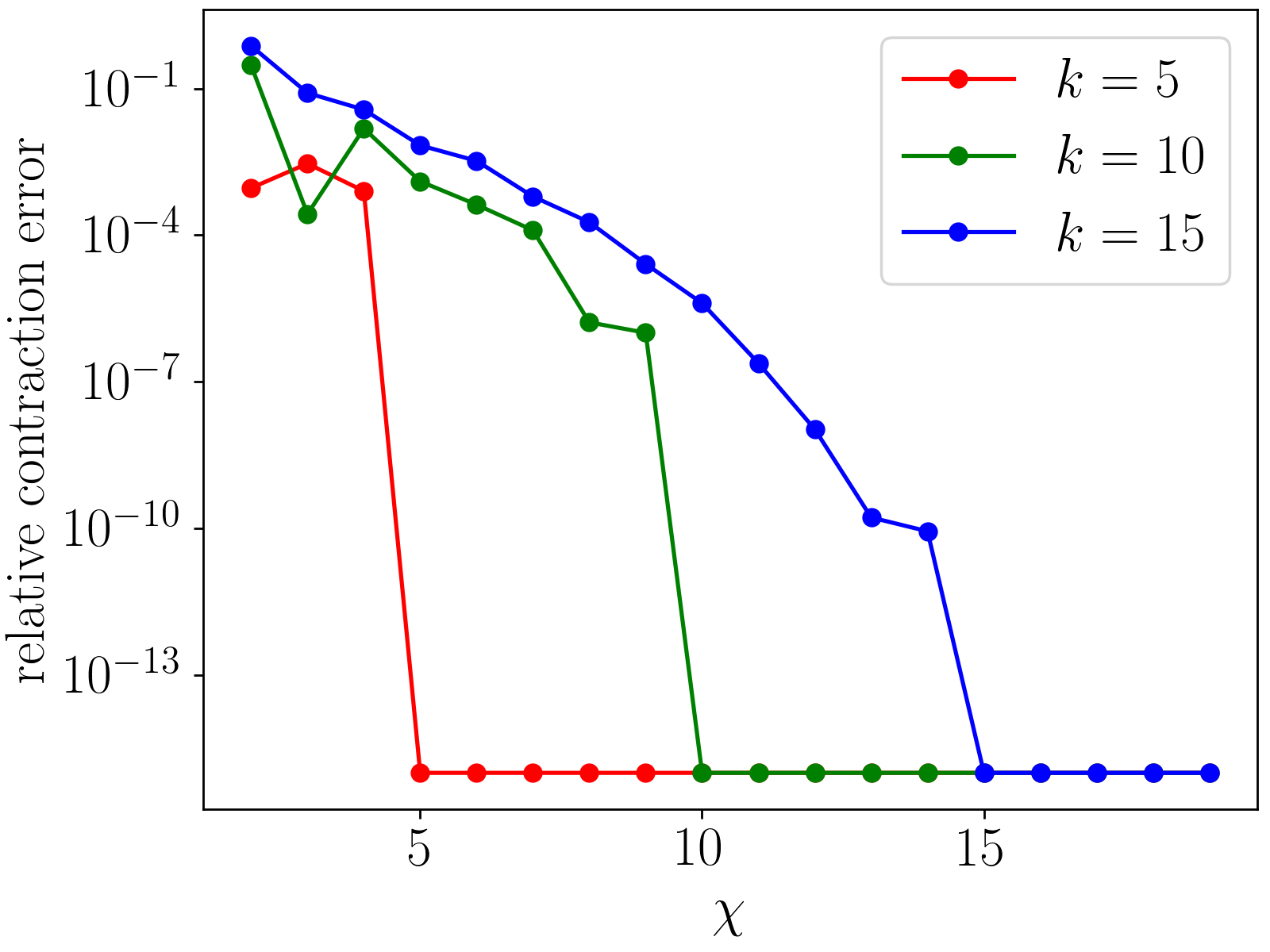}
    \caption{}
    \label{fig:pol_exact}
    \end{subfigure}
    \begin{subfigure}[b]{0.08\linewidth}
    \centering
    \includegraphics[width=0.9\textwidth]{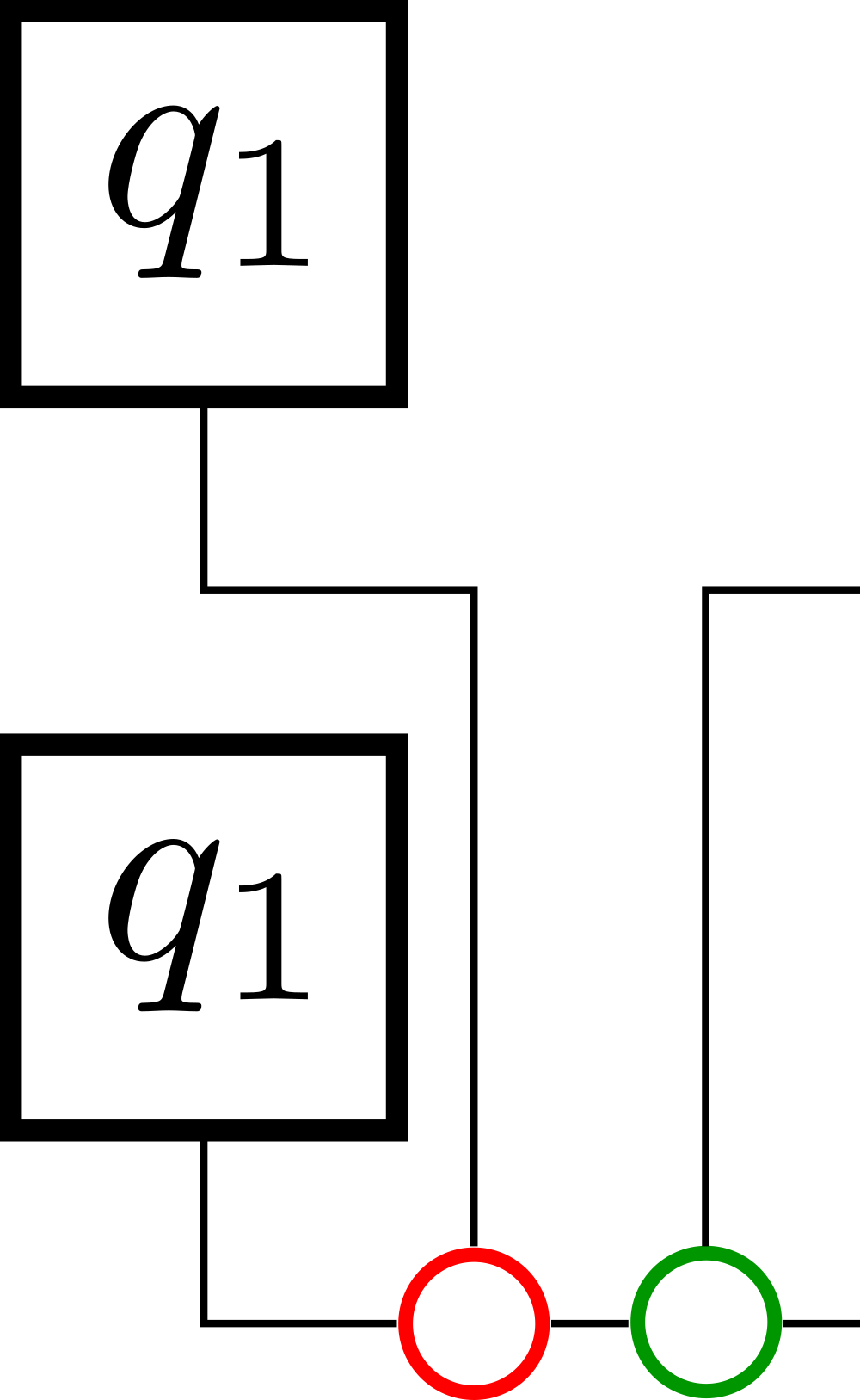}
    \caption{}
    \label{fig:proja}
    \end{subfigure}
    \begin{subfigure}[b]{0.15\linewidth}
    \centering
    \includegraphics[width=0.9\textwidth]{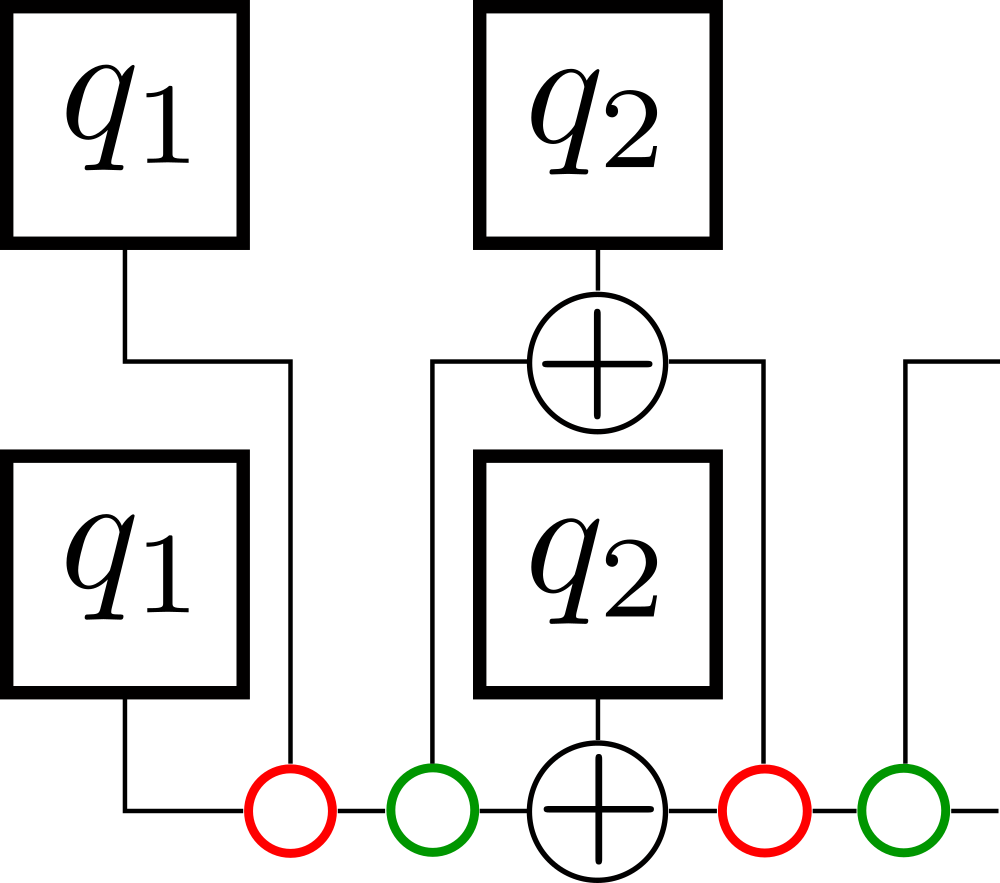}
    \caption{}
    \label{fig:projb}
    \end{subfigure}
    \begin{subfigure}[b]{0.15\linewidth}
    \centering
    \includegraphics[width=0.9\textwidth]{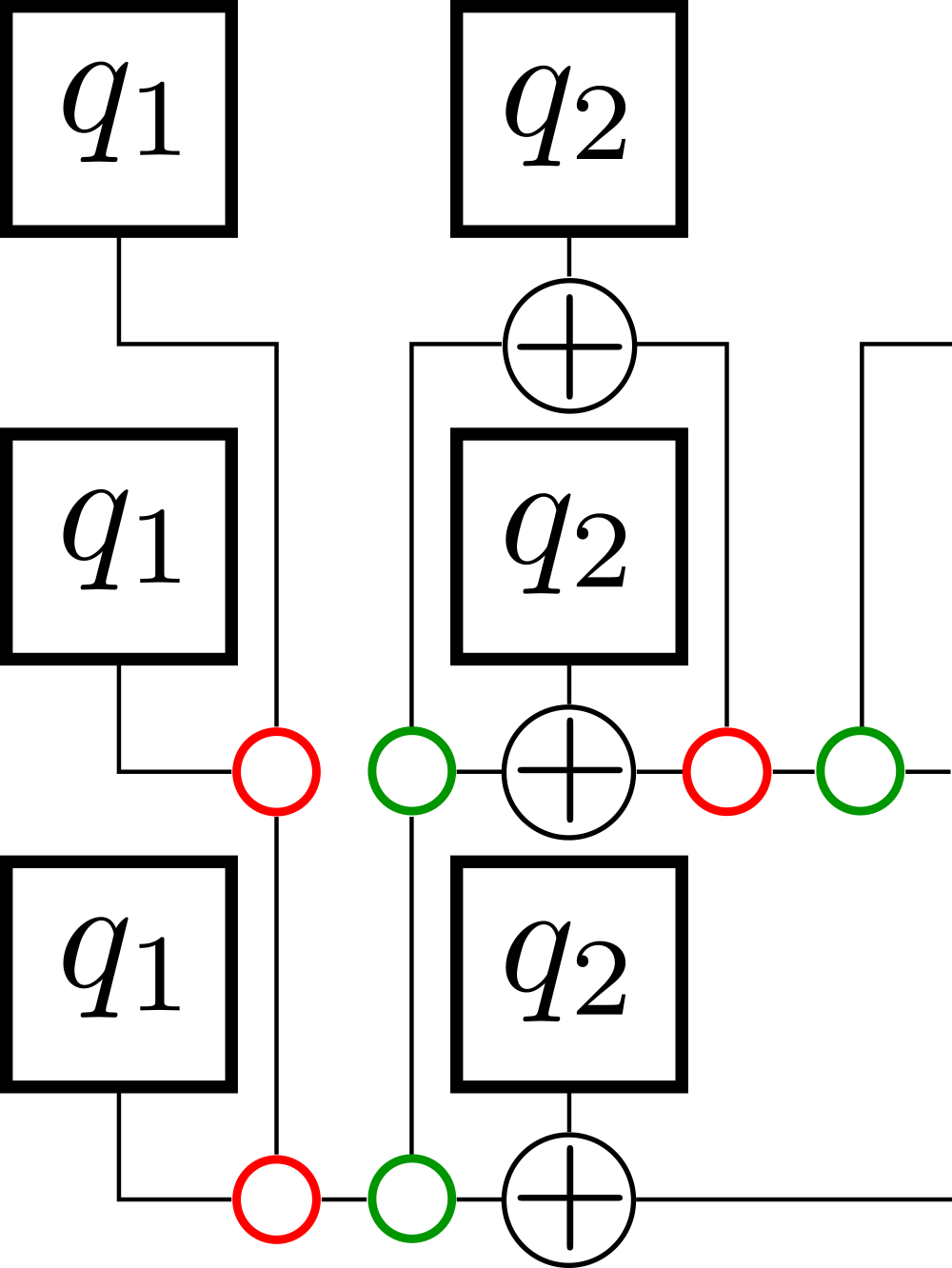}
    \caption{}
    \label{fig:projc}
    \end{subfigure}
    \begin{subfigure}[b]{0.25\linewidth}
    \centering
    \includegraphics[width=0.9\textwidth]{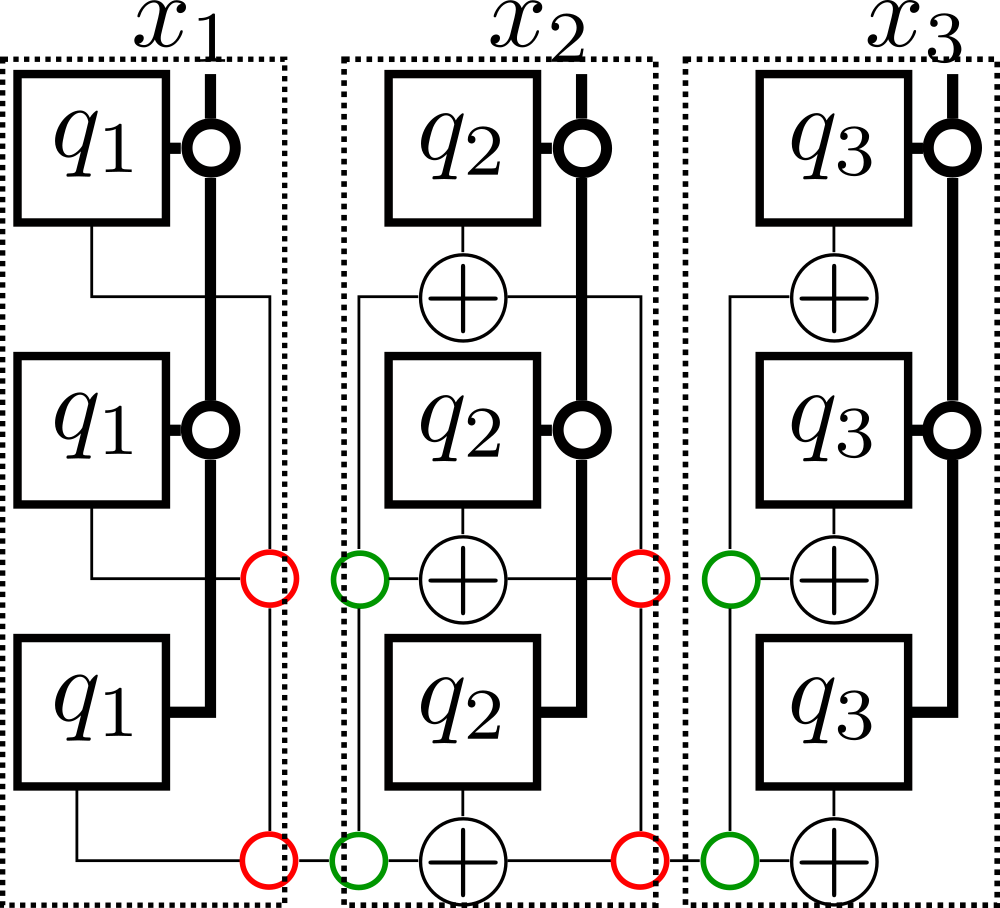}
    \caption{}
    \label{fig:projd}
    \end{subfigure}
    \caption{(a) Accuracy of TN integration (relative contraction error) for the integral of the multivariable function polynomial in Eq.~(\ref{eqn:pol_exact}). This demonstrates exact compressibility via SVD for the 2D arithmetic circuit TN, with number of variables $N=20$, number of points per variable $G=10$. The exact compressibility is equivalent to inserting projectors $P_L[i,j]$ (red), $P_R[i,j]$ (green) as in diagrams (b)$\sim$(e). (b) inserts projector $P_L[2,1]$ and $P_R[2,1]$. (c) inserts $P_L[2,2]$ and $P_R[2,2]$. (d) inserts projectors $P_L[3,1]$, $P_R[3,1]$. (e) Full TN for $N=3$, $k=3$ with projectors inserted. The result from each grouped column can be related to a recursive computation of the integral. See Sec.~\ref{sec:pol_exact} for details. }%
    \label{fig:proj}%
\end{figure*}

We first consider the case where each factor is identical, i.e. $p_i(x_1, \ldots, x_N) = p(x_1, \ldots, x_N)$. Then the total function takes the form
\begin{align}\label{eqn:pol_exact}
f(x_1,...,x_N)=(q_1(x_1)+...+q_N(x_N))^k,
\end{align}
To perform the integration, we first discretize each variable $x_i$ on a grid of $G$ points, and allow each single variable function on the grid to take random values between $-1, 1$ (thus each $q_i$ represents a discretized version of a function that is oscillating between $-1$ and $1$). We use an equally weighted quadrature.

The contraction of the 2D tensor network corresponds to the contraction of the PEPS, thus we use an approximate tensor network contraction strategy commonly used in PEPS, compressing to a finite bond dimension $\chi$. Within the PEPS structure, the horizontal bonds have dimension $2$, while the vertical bonds associated with the grid points have dimension $G$. For $G\gg 2$, the cost of contracting along the vertical dimension is much cheaper than contracting along the horizontal dimension. Thus we perform the contraction by a boundary row contraction method, contracting rows into rows (see Fig.~\ref{fig:boundary}).

After a row is contracted into a row, the dimension of the tensors along the boundary increases; the maximum bond dimension of the boundary tensors after $k$ rows is $2^k$. However, when performing compression, we immediately find that the boundary tensors are \emph{exactly compressible} to $\chi < 2^k$ after each row contraction. In fact, the entire tensor network can be exactly contracted with bond dimension $\chi=k$, as illustrated in Fig.~\ref{fig:pol_exact}, which shows the error in the computed integral (due to compressing to finite bond dimension $\chi$); we see the error of contraction drops to $0$ for $\chi = k$. 
Overall this means that the integral can then be computed exactly with cost linear in the number of variables $N$, and polynomial in the function nonlinearity $k$. 

The ability to exactly contract the network with small $\chi$ implies the existence of an exact algebraic identity. 
Each compression corresponds to the insertion of isometries/projectors into the 2D TN. 
For each row $i=2, \ldots, k$ and column $j=1,\ldots,N$, we insert left and right projectors $P_L[i,j]$ (of dimension $i \times 2 \times (i+1)$) and $P_R[i,j]$ (of dimension $(i+1) \times 2 \times i$) between row $i-1$, $i$ and column $j-1$, $j$ where
\begin{align}
&(P_L[i,j])_{ab,c}=\bigg\{\begin{array}{cc}
    1 & a=c,b=0 \\
    1 & a=j-1,c=j,b=1 \\
    0 & \text{otherwise}
\end{array}\\
&(P_R[i,j])_{c,ab}=\bigg\{\begin{array}{cc}
    1 & a+b=c \\
    0 & \text{otherwise}
\end{array}.
\end{align}
This is shown step by step in Figs.~\ref{fig:proja}-\ref{fig:projd}.
For example, in Fig.~\ref{fig:proja}, we insert a pair of projectors $P_L[2,1]$, $P_R[2,1]$ where $P_L[2,1]_{ab,c}$ takes two input vectors (the $a$, $b$ indices) $[1, q_1]$, $[1, q_1]$ and maps them to an output vector (the $c$ index) $[1, q_1, q_1^2]$;  $P_R[2,1]_{c,ab}$ takes the input vector ($c$ index) $[1,q_1,q_1^2]$ and maps it to two output vectors ($a$, $b$ indices) $[1,q_1]$, $[1,q_1]$. Similarly, in Fig.~\ref{fig:projb}, we insert the pair of projectors $P_L[2,2]$, $P_R[2,2]$, where $P_L[2,2]$ takes two input vectors $[1,q_1+q_2]$, $[1,q_1+q_2]$ and maps them to an output vector $[1,q_1+q_2,(q_1+q_2)^2]$; $P_R[2,2]$ takes the input vector $[1,q_1+q_2,(q_1+q_2)^2]$ and maps it to two output vectors $[1,q_1+q_2]$, $[1,q_1+q_2]$. In general, 
\begin{align}
&P_L[i,j]:\left[\begin{array}{c}
    1 \\
    s_j \\
    ... \\
    s_j^{i-1}
\end{array}\right]\otimes\left[\begin{array}{c}
    1 \\
    s_j
\end{array}\right]\to\left[\begin{array}{c}
    1 \\
    s_j \\
    ... \\
    s_j^i
\end{array}\right]\\
&P_R[i,j]:\left[\begin{array}{c}
    1 \\
    s_j \\
    ... \\
    s_j^i
\end{array}\right]\to\left[\begin{array}{c}
    1 \\
    s_j \\
    ... \\
    s_j^{i-1}
\end{array}\right]\otimes\left[\begin{array}{c}
    1 \\
    s_j
\end{array}\right]. 
\end{align}
where $s_j=\sum_{j'=1}^jq_{j'}$ is the sum of single variable functions up to the $j$-th term. Thus, we see the role of the left projector is to retain only the non-redundant polynomials of the single variable functions $q_j$, and the right projector redistributes non-redundant powers into a direct product. 

The tensor network exact compression then corresponds to computing the integral through a recursive formula (see Ref.~\onlinecite{Oseledets2013} for a related constructive approach). For example, we can recursively define a set of integrals over a subset of variables 
\begin{align}
I_1^{k_1}=\int dx_1q_1(x_1)^{k_1}
\end{align}
\begin{align}
I_2^{k_2}
&=\int dx_1dx_2(q_1(x_1)+q_2(x_2))^{k_2}\notag\\
&=\int dx_1dx_2\sum_{k_1=0}^{k_2}\binom{k_2}{k_1}q_1(x_1)^{k_1}q_2(x_2)^{k_2-k_1}\notag\\
&=\sum_{k_1=0}^{k_2}\binom{k_2}{k_1}I_1^{k_1}\int dx_2q_2(x_2)^{k_2-k_1}
\end{align}
\begin{align}
I_3^{k_3}
&=\int dx_1dx_2dx_3(q_1(x_1)+q_2(x_2)+q_3(x_3))^{k_3}\notag\\
&=\int dx_1dx_2dx_3\bigg[\sum_{k_2=0}^{k_3}\binom{k_3}{k_2}\notag\\
&(q_1(x_1)+q_2(x_2))^{k_2}q_3(x_3)^{k_3-k_2}\bigg]\notag\\
&=\sum_{k_2=0}^{k_3}I_2^{k_2}\int dx_3q_3(x_3)^{k_3-k_2}
\end{align}
and so on up to 
\begin{align}
I_N^{k}
&=\int dx_1...dx_N(q_1(x_1)+...+q_N(x_N))^{k}\notag\\
&=\int dx_1...dx_N\bigg[\sum_{k_{N-1}=0}^k\binom{k}{k_{N-1}}\notag\\
&(q_1(x_1)+...+q_{N-1}(x_{N-1}))^{k_{N-1}}q_N(x_N)^{k-k_{N-1}}\bigg]\notag\\
&=\sum_{k_{N-1}=0}^k\binom{k}{k_{N-1}}I_{N-1}^{k_{N-1}}\int dx_Nq_N(x_N)^{k-k_{N-1}}
\end{align}
where each iterative step involves integration of single variable integrals. The quantities $I_{j}^{k_j}$ can be related to the TN with the projectors inserted in each layer.  Fig.~\ref{fig:projd} shows a simple example for $N=3$, $k=3$.
From left to right, the open legs of the first two dashed boxes correspond to $I_1^{k_1}$, $I_2^{k_2}$ respectively for $k_1,k_2=0,\dots,3$, and the full contraction of the TN corresponds to $I_3^k$ for $k=3$. 

\subsection{Perturbations away from exact compressibility}\label{section:pol_gen}

\begin{figure}%
    \centering
    \begin{subfigure}[b]{0.8\linewidth}
    \centering
    \includegraphics[width=\textwidth]{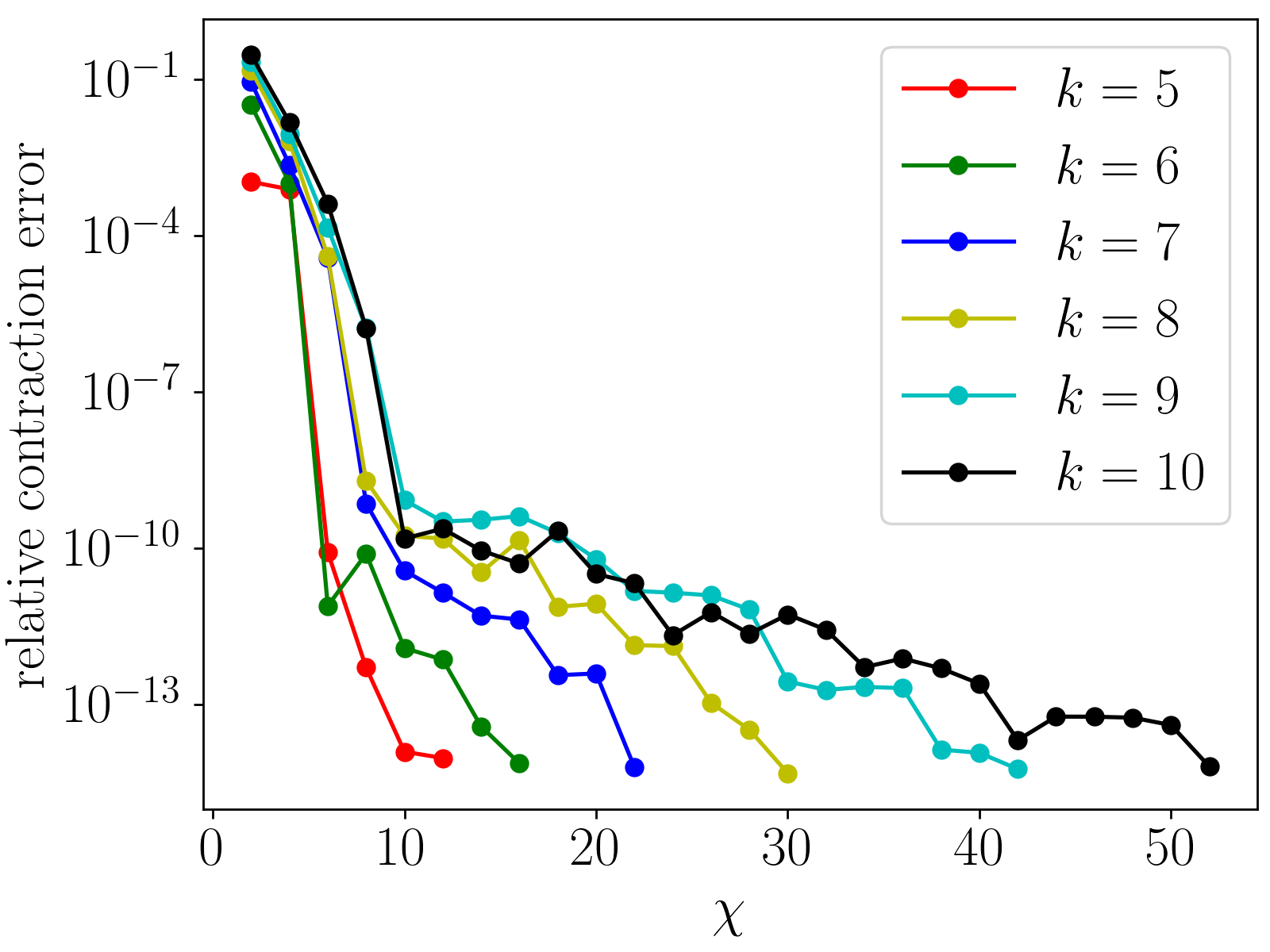}
    \caption{$\delta=0.01$}
    \end{subfigure}
    \begin{subfigure}[b]{0.8\linewidth}
    \centering
    \includegraphics[width=\textwidth]{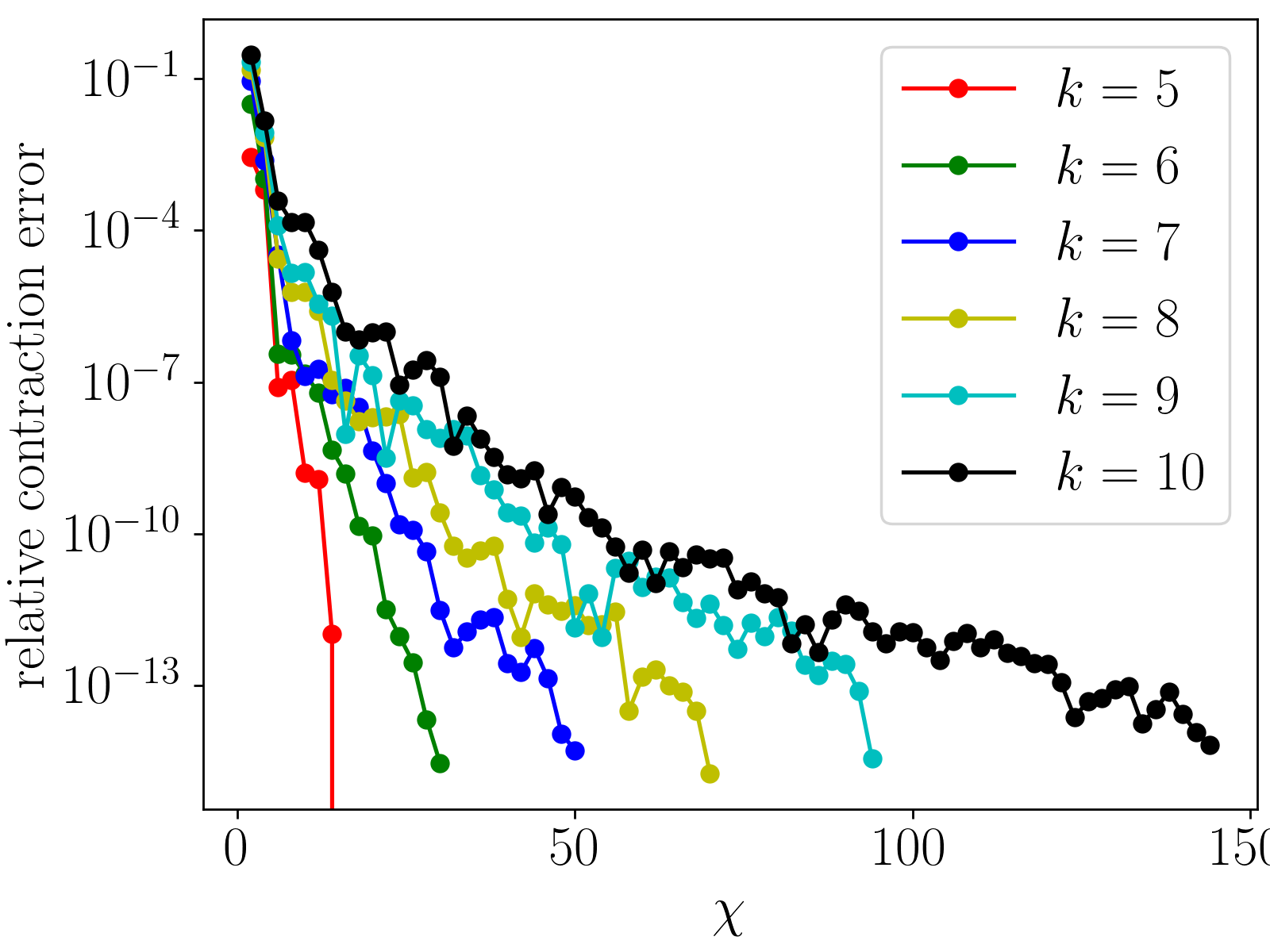}
    \caption{$\delta=0.1$}
    \end{subfigure}
    \caption{Accuracy of TN integration for the integral of a multivariable function polynomial perturbed away from the exactly compressible point as a function of contraction bond dimension $\chi$. Magnitude of perturbation $\delta$, number of variables,  $N=20$, number of points per variable, $G=10$. $k$ is a measure of the nonlinearity expressed in the TN. Note that there is fast convergence with $\chi$ up to a critical precision, which depends on $\delta$ and $k$. See Sec.~\ref{section:pol_gen} for details.}%
    \label{fig:pol_perturb}%
\end{figure}

In the above, the compression of the circuit tensor network for the multivariable polynomial function allows us to identify an efficiently integrable case. However, when inserting more general polynomials in Eq.~(\ref{eqn:pol}), we cannot expect the tensor network to be exactly compressible since we know the general case is NP hard. However, we might expect that functions that are close to the efficiently integrable case remain efficiently integrable up to some accuracy.

We thus now consider polynomial functions in Eq.~(\ref{eqn:pol}) that are obtained by perturbing away from the exactly compressible case. We do this by defining 
\begin{align}\label{eqn:pol_perturb}
q_{ji}(x_j)=q_{j1}(x_j)+\delta\cdot r_{ji}(x_j)
\end{align}
where the single variable functions $q_{j1}(x_j)$, $r_{ji}(x_j)$ are random length-$G$ vectors with values in [-1,1] (recall, $G$ is the number of grid points). For $i=1$, we set $\delta=0$ for $j=1,\ldots,N$. Then for $i>1$, the difference between $q_{ji}(x_j)$ and $q_{j1}(x_j)$ is controlled by $\delta$. In Fig.~\ref{fig:pol_perturb} we use perturbations of magnitude $\delta=0.01$ and $\delta = 0.1$ for $N=20$, $G=10$, and we monitor the error in the multivariable integral as a function of $\chi$. In both cases we see that the integration error shows two regions of convergence; first there is a rapid convergence to some finite error (around $10^{-10}$ for $\delta=0.01$, and $10^{-4}$ for $\delta=0.1$), followed by a much slower convergence to smaller error. In the rapidly converging regime, the required bond dimension for a given relative error grows as $\chi \sim O(k)$; in the slow convergence regime, the required bond dimension $\chi \sim \exp(k)$. In all cases, however, the cost to  integrate is linear in the number of variables $N$, and this is independent of the desired precision (not shown for this case, but see next section). Exponential dependence of the complexity appears instead in the cost to improve the precision past a critical threshold. 

\subsection{General case}\label{sec:gen_gen}

\begin{figure*}%
    \centering
    \begin{subfigure}[b]{0.32\linewidth}
    \centering
    \includegraphics[width=\textwidth]{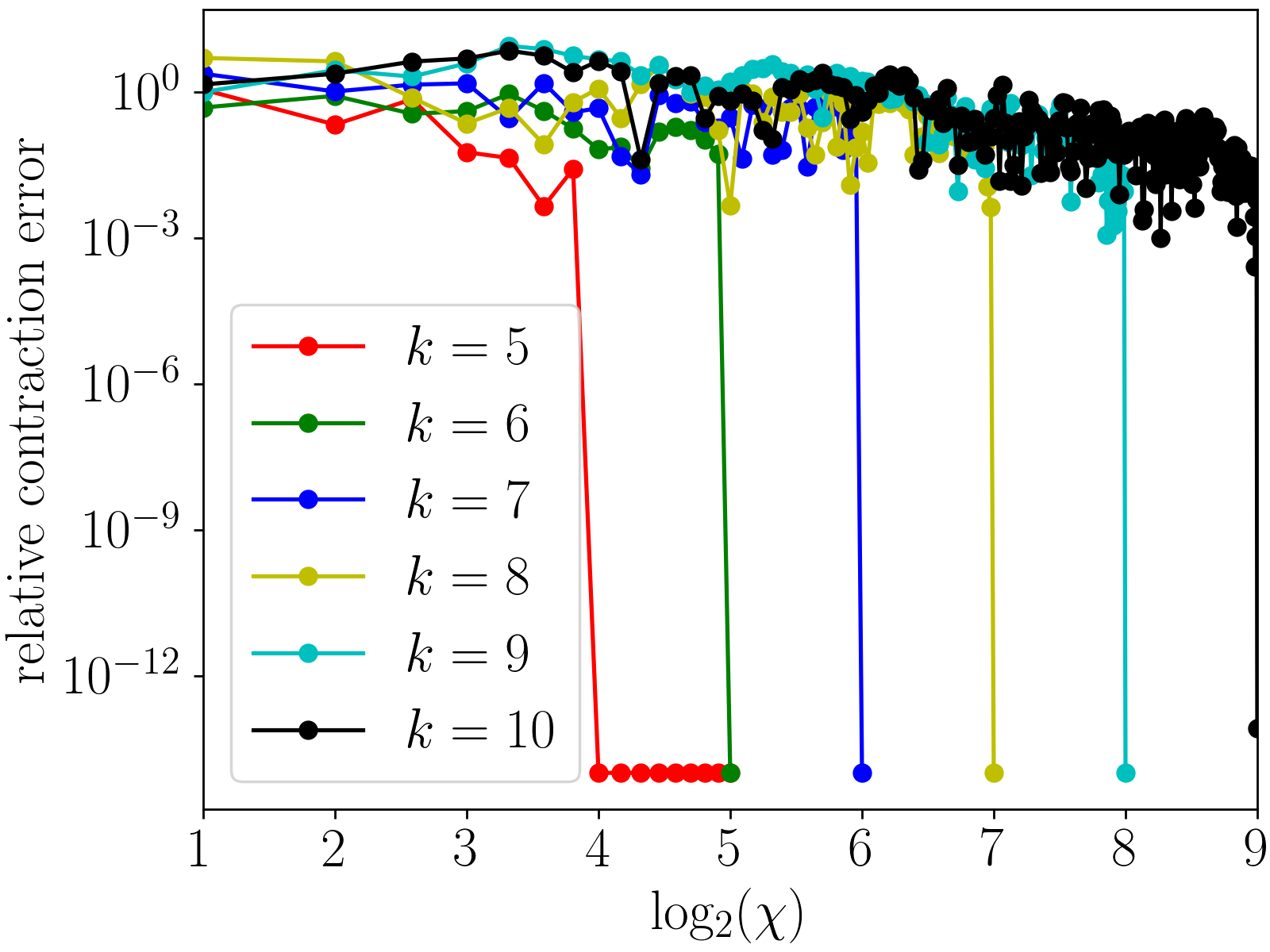}
    \caption{$\lambda=-1$}
    \label{fig:pol_randa}
    \end{subfigure}
    \begin{subfigure}[b]{0.32\linewidth}
    \centering
    \includegraphics[width=\textwidth]{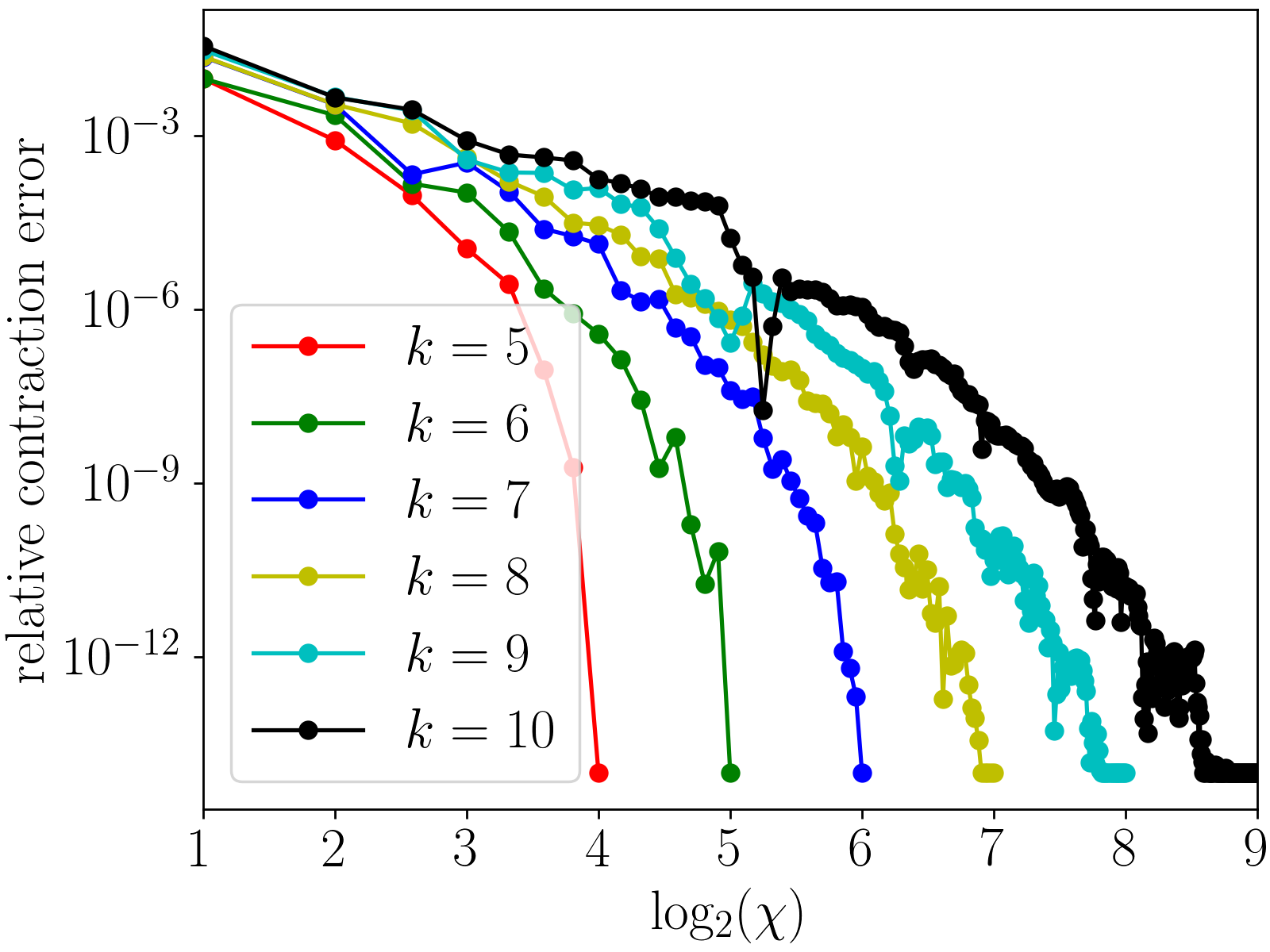}
    \caption{$\lambda=-0.5$}
    \label{fig:pol_randb}
    \end{subfigure}
    \begin{subfigure}[b]{0.32\linewidth}
    \centering
    \includegraphics[width=\textwidth]{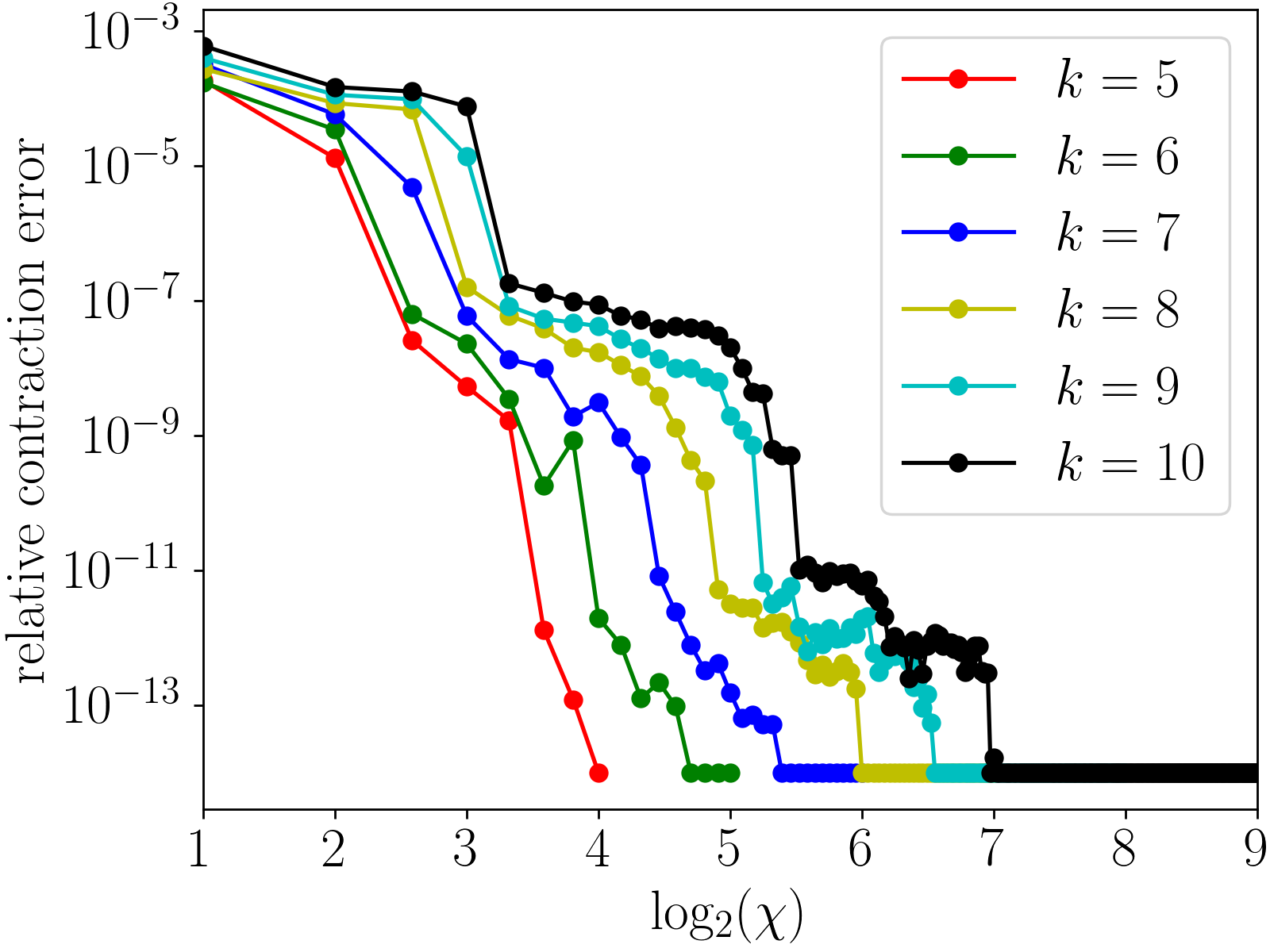}
    \caption{$\lambda=0$}
    \label{fig:pol_randc}
    \end{subfigure}
    
    %\centering
    \begin{subfigure}[b]{0.32\linewidth}
    \centering
    \includegraphics[width=\textwidth]{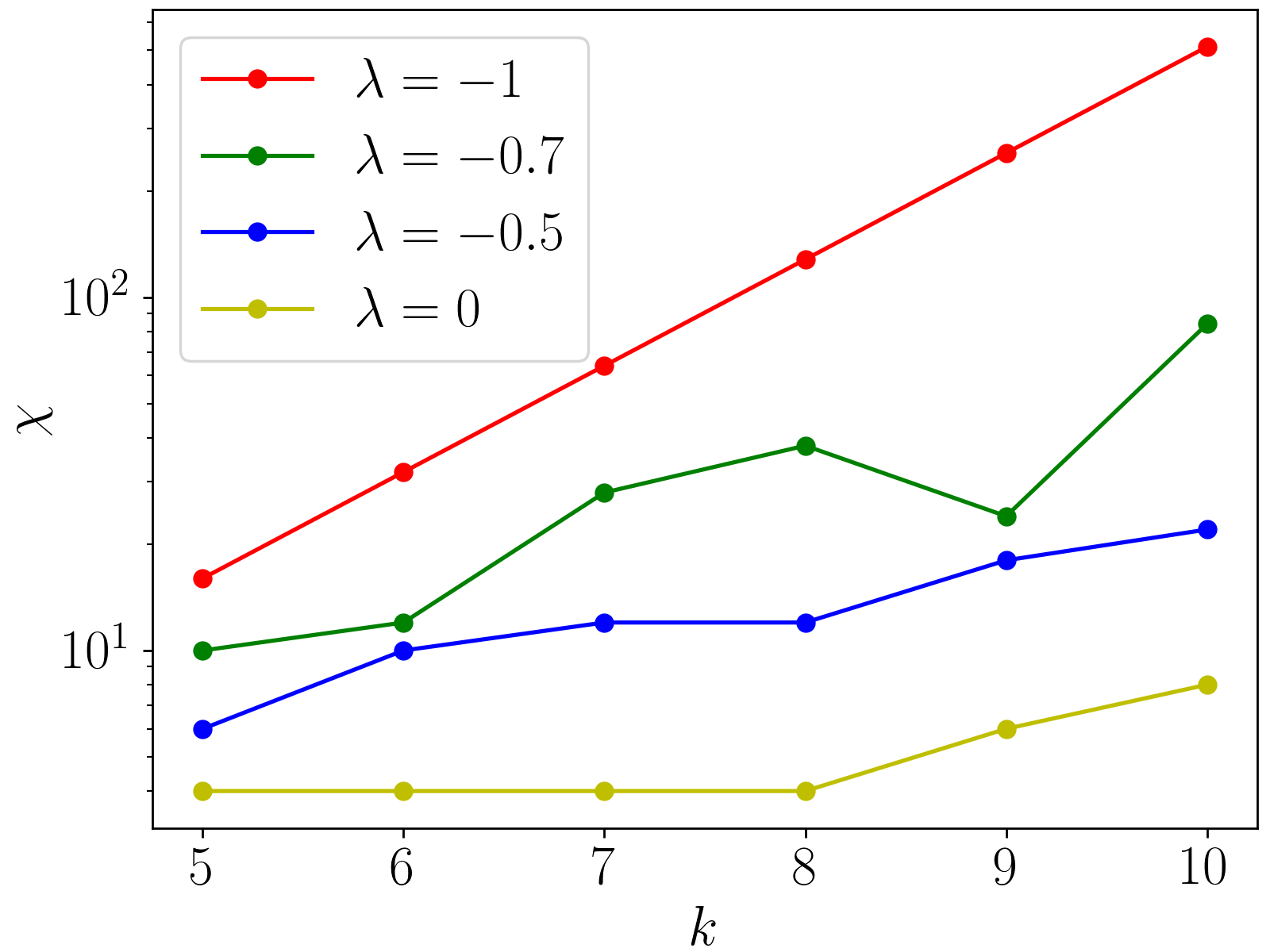}
    \caption{}
    \label{fig:pol_gena}
    \end{subfigure}
    \begin{subfigure}[b]{0.32\linewidth}
    \centering
    \includegraphics[width=\textwidth]{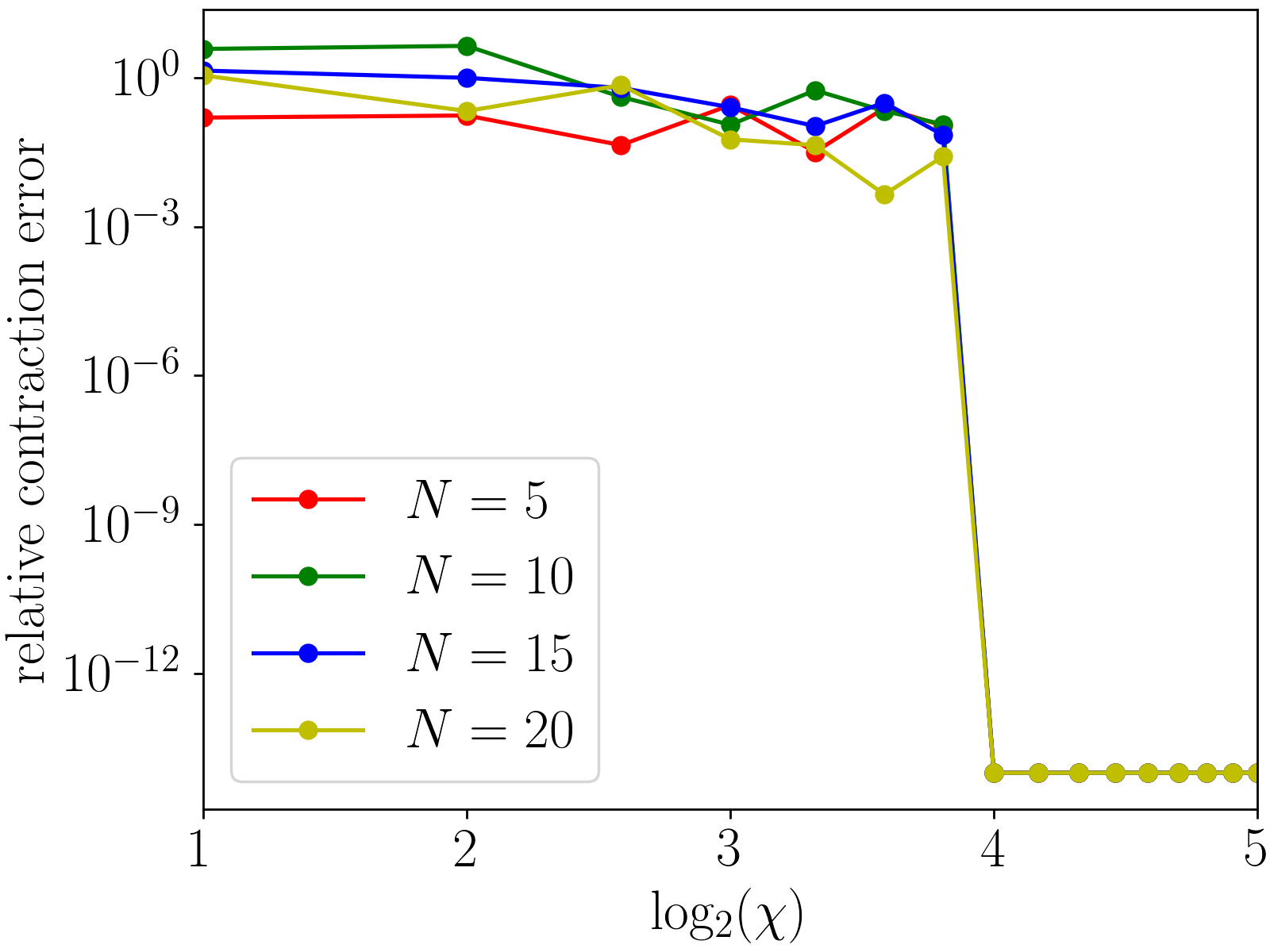}
    \caption{}
    \label{fig:pol_genb}
    \end{subfigure}
    \begin{subfigure}[b]{0.32\linewidth}
    \centering
    \includegraphics[width=\textwidth]{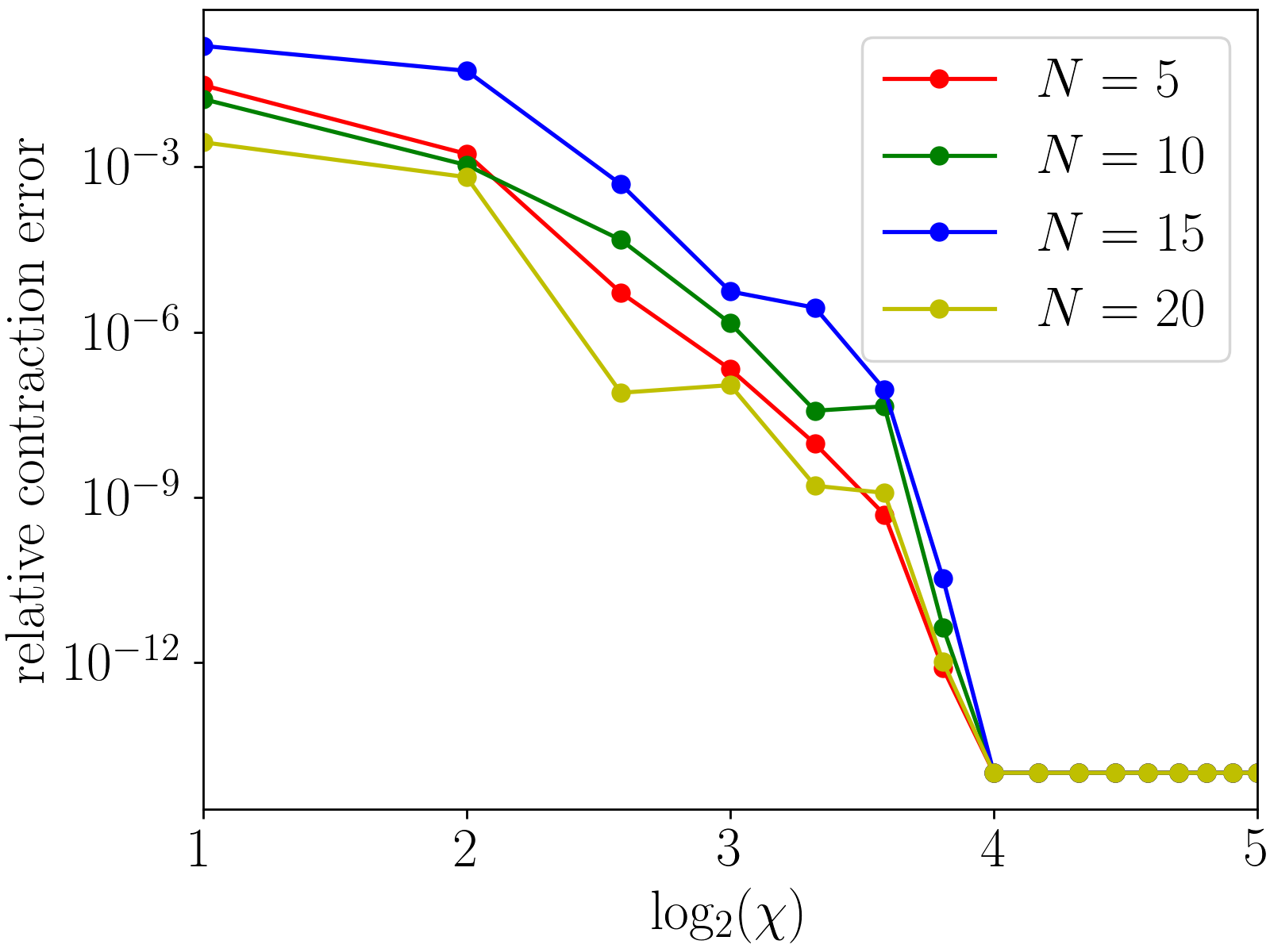}
    \caption{}
    \label{fig:pol_genc}
    \end{subfigure}
    \caption{Accuracy of TN integration of multivariable function polynomials.  For (a), (b), (c), we plot accuracy v.s. $\chi$ for number of variables $N=20$, number of points per variable $G=10$, with random single variable functions $q_{ji}$ with values $\in[\lambda,1]$, for three different $\lambda$; the difficulty of integration changes with the positivity of the integrand.
    (d) Bond dimension $\chi$ v.s. nonlinearity $k$ to achieve a relative accuracy $\leq10^{-4}$ for various $\lambda$ for $N=20$, $G=10$. For (e)/(f), we plot relative accuracy v.s. $\chi$ for the general/perturbative ($\delta=0.1$) function polynomial case with $k=5$, $G=10$ and single variable functions $q_{ji},r_{ji}\in[-1,1]$, demonstrating the independence of accuracy with respect to the number of variables $N$. See Sec.~\ref{sec:gen_gen} for details.}%
    \label{fig:pol_rand}
\end{figure*}

For the general polynomial in Eq.~(\ref{eqn:pol}), we can expect the bond dimension for a given relative precision to scale $\sim \exp(k)$ as discussed; the tensor network is, in some sense, incompressible. However, in some cases, we may still be able to approximate the integral to good accuracy. 
This depends strongly on the integrand range, as illustrated in Fig.~\ref{fig:pol_rand}.  In Fig.~\ref{fig:pol_randa}, we take $N=20$, $G=10$ and single variable function values chosen randomly in the range $[-1,1]$. Then, the corresponding 2D TN representation of the integral appears almost completely incompressible: only until $\chi$ reaches the bond dimension of the exact contraction $\sim \exp(k)$ do we suddenly see a significant improvement in the integral error (although for rough estimates, e.g. to $10^{-1}$ precision, it is possible to take $\chi$ orders of magnitude below that required for exact contraction). We can make the single variable functions less oscillatory by increasing the lower bound of the range of the single variable functions, i.e. $[\lambda, 1]$, where $\lambda$ is increased from $-1$ to $0$. As we do so,  the integration problems become easier, as can be seen from Figs.~\ref{fig:pol_randa},~\ref{fig:pol_randb},~\ref{fig:pol_randc}. The decrease of relative error with increasing bond dimension is much faster as we raise the lower bound $\lambda$ of the the single variable functions. 

Monte Carlo integration can face difficulties with oscillatory functions with small or vanishing integrals. In the circuit tensor network representation, there is an exponential dependence for such oscillatory functions, but unlike in quadrature, the exponential cost is in the non-linearity (parameterized by $k$), not in the number of variables $N$. Fig.~\ref{fig:pol_gena} plots the required bond dimension for a relative error $10^{-4}$ with $N=20$, $G=10$ as a function of the non-linearity $k$ for various single variable function ranges $q_{ji}\in[\lambda, 1]$, where $\lambda$ is increased from $-1 \ldots 0$. We see that the required bond dimension increases exponentially in the non-linearity. However,
the exponent is much larger for oscillatory integrands ($\lambda=-1$) than for more positive integrands $\lambda > -1$. 

We finally confirm that the error of compressed contraction is essentially independent of the number of variables using boundary contraction as can be seen from Figs.~\ref{fig:pol_genb},~\ref{fig:pol_genc}. Either for the general functional form Eq.~(\ref{eqn:pol}) with highly oscillatory values from random single variable functions, or for the perturbation from the exact case Eq.~(\ref{eqn:pol_perturb}) with $\delta=0.1$, the rate of decrease of relative error with increasing bond dimension is comparable as we change the number of variables $N$. 

\subsection{Comparison with quasi-Monte Carlo}\label{sec:compare_qmc}

\begin{figure*}[]%
    \centering
    \begin{subfigure}[b]{0.4\linewidth}
    \centering
    \includegraphics[width=\textwidth]{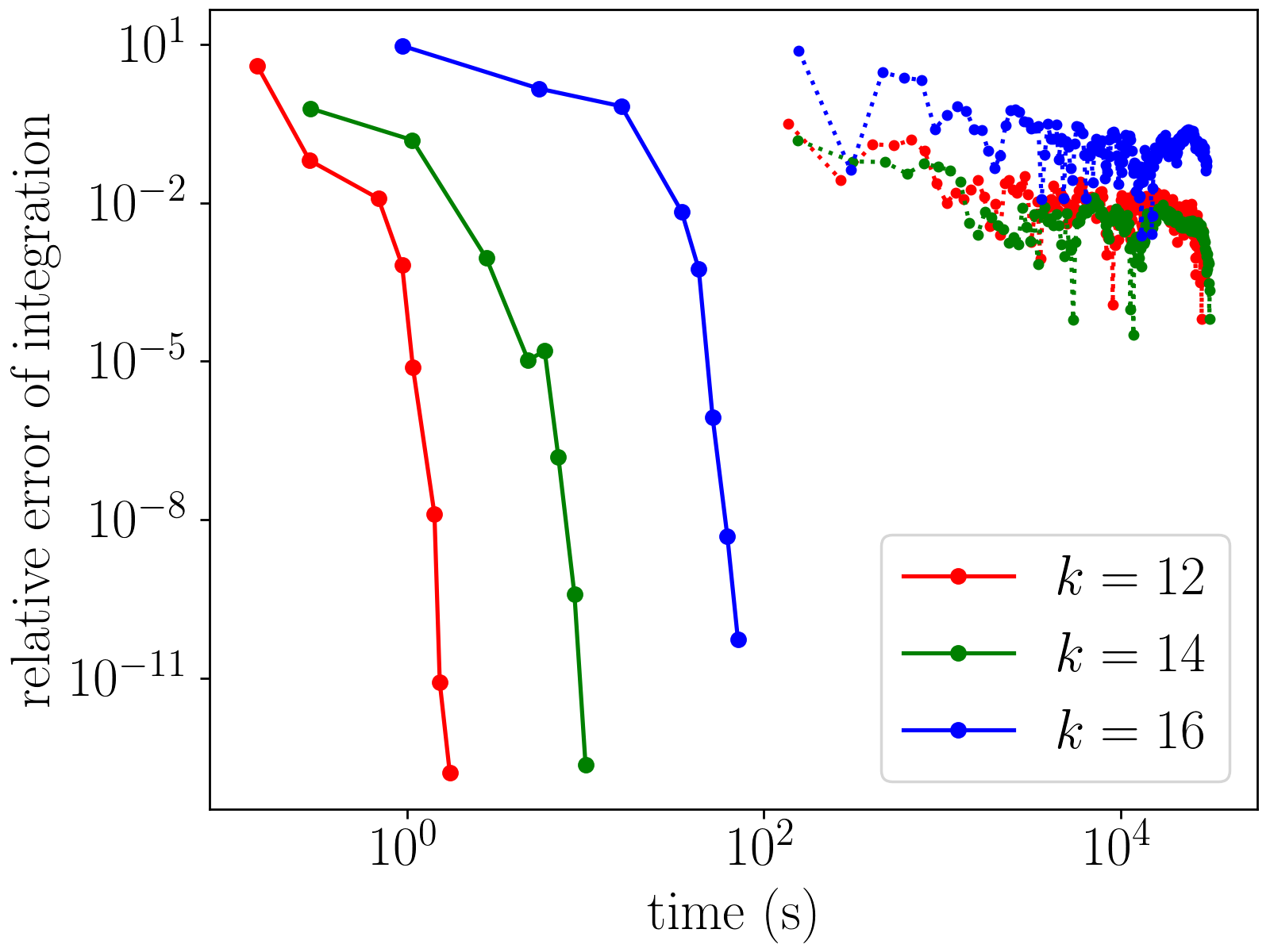}
    \caption{}
    \label{fig:pol_qmca}
    \end{subfigure}
    \begin{subfigure}[b]{0.4\linewidth}
    \centering
    \includegraphics[width=\textwidth]{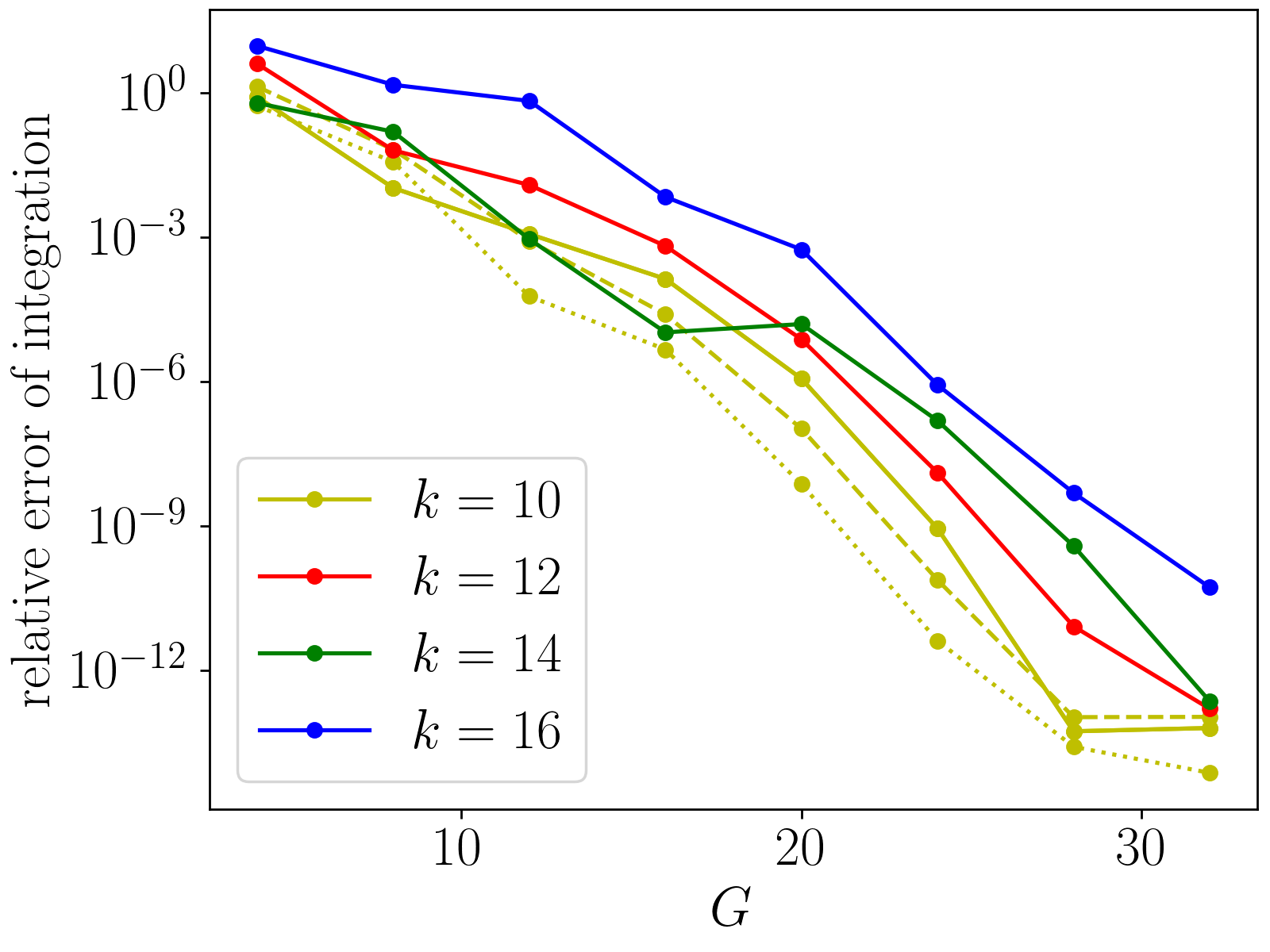}
    \caption{}
    \label{fig:pol_qmcb}
    \end{subfigure}
    %\centering
    
    \begin{subfigure}[b]{0.4\linewidth}
    \centering
    \includegraphics[width=\textwidth]{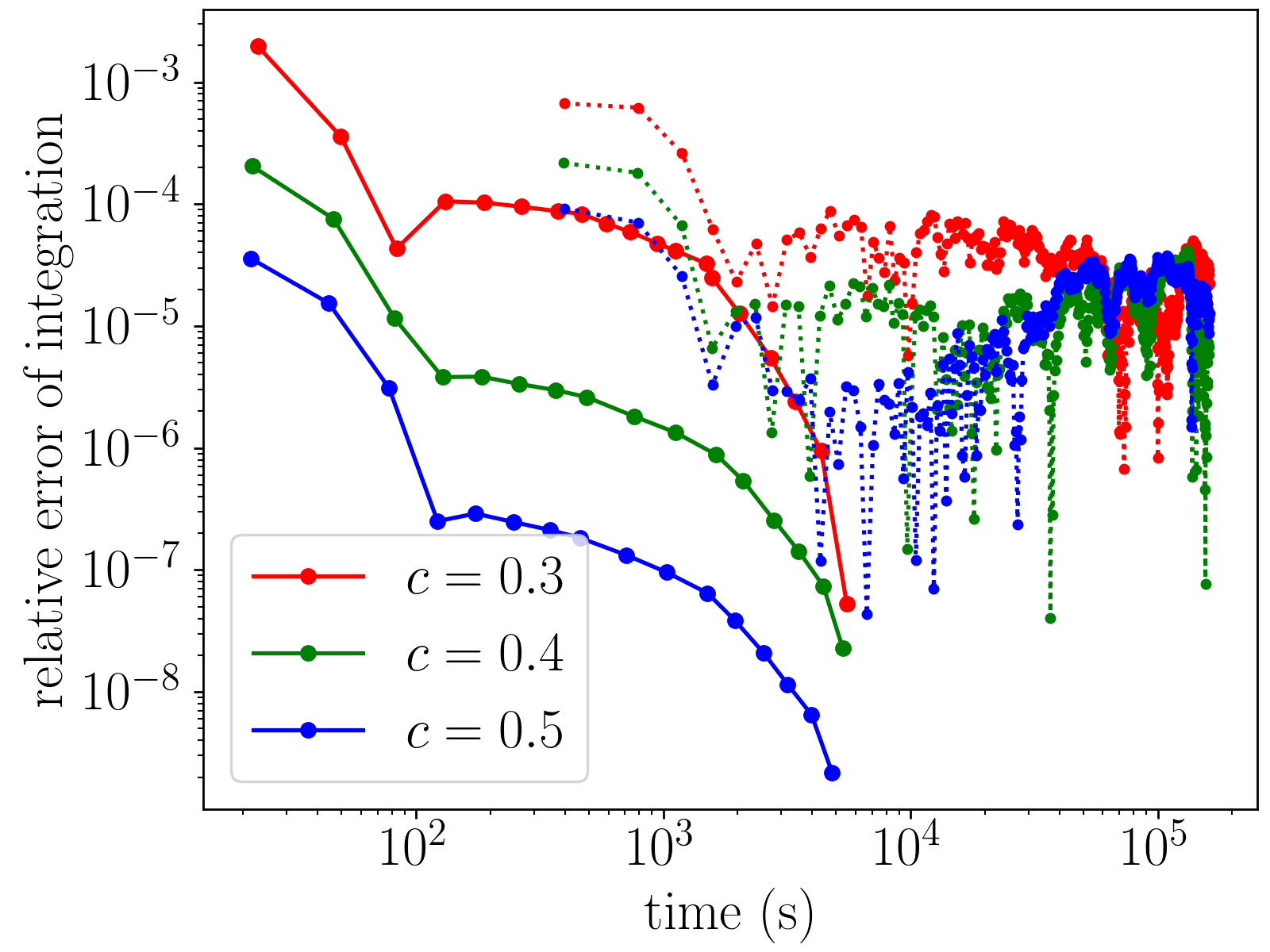}
    \caption{}
    \label{fig:pol_qmcc}
    \end{subfigure}
    \begin{subfigure}[b]{0.4\linewidth}
    \centering
    \includegraphics[width=\textwidth]{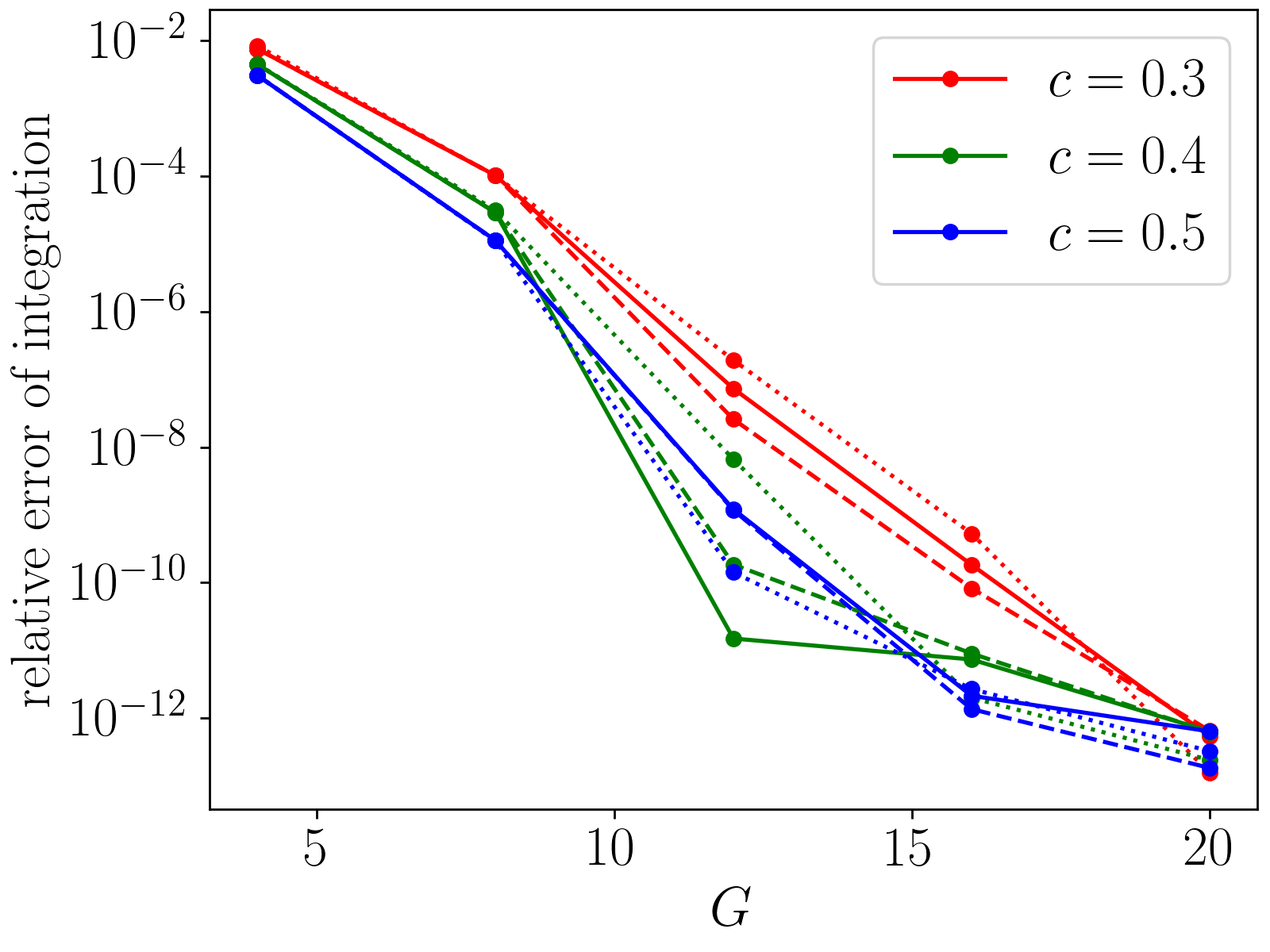}
    \caption{}
    \label{fig:pol_qmcd}
    \end{subfigure}
    \caption{Comparison of TN and quasi-MC integration for the multivariable function polynomial in Eqs.~\ref{eqn:pol},~\ref{eqn:pol_fac},~\ref{eq:contpolynomial}. $N$ is the number of variables, $k$ is the number of factor powers in the function, and $c$ is a function parameter that controls the positivity of the integrand. $G$ is the number of quadrature points per variable in the TN integration. (a) TN and quasi-MC integral convergence with runtime for various $k$ at $N=10$, $c=0$. Exact result taken as TN integral at $G=36$, bond dimension $\chi=2^k$. (b) TN integral convergence w.r.t. $G$ at $\chi=2^k$ for various $k$ at $c=0$. Solid/dashed/dotted lines are for $N=10/20/40$. (c) TN and quasi-MC integral convergence with time for various $c$ at $N=50$, $k=30$. Exact result taken as TN integral at $G=12$, $\chi=560$. (d) TN integral convergence w.r.t. $G$ for various $c$ at $N=50$, $k=30$. Solid/dashed/dotted lines are for $\chi=180/120/60$. See Sec.~\ref{sec:compare_qmc} for further details. }%
    \label{fig:pol_qmc}%
\end{figure*}

To understand the concrete performance of integration using the arithmetic circuit TN, we now compare costs with that of quasi-MC. For this, we construct the polynomial in Eq.~(\ref{eqn:pol}) with the following multivariable function
\begin{align}
q_{ji}(x_j)=\sin{(2\pi(x_j+a_{ji}))}+c \label{eq:contpolynomial}
\end{align}
for random $a_{ji}\in[0,1]$ and some constant $c$, and integrate over $x_j\in[0,1]$. We compare the results from  TN contraction to quasi-MC integration in Fig.~\ref{fig:pol_qmc} (we report representative results for random $a_{ji}$) changing both the number of variables $N$ as well as the polynomial power (nonlinearity) $k$. All calculation are done on Intel(R) Xeon(R) CPU E5-2697 v4 processors. For the quasi-MC calculations, the integrand in Eq.~(\ref{eqn:pol}) was coded in Python then compiled by JAX\cite{jax2018github}, and the sample points were generated using the QMCPY\cite{QMCPy} package with Sobol generating matrices using a batch size of $10^7$, with the function values for each batch computed on a single core. The TN timings are reported as the runtime for contraction using the \textsc{Quimb}~\cite{gray2018quimb} package on a single core or two cores, with the time normalized to a single core time. The reference exact result for the integration was taken as the converged TN result w.r.t $\chi$ and $G$ (specified in the caption of Fig.~\ref{fig:pol_qmc}). 

Fig.~\ref{fig:pol_qmca} plots the convergence of the integral versus time for TN integration and quasi-MC for $N=10$, $c=0$ and various $k$. The TN data points represent increasing grid order $G$ using Gauss-Legendre quadrature and bond dimension $\chi=2^k$. For moderate $k$, quasi-MC already faces convergence difficulties as a result of the highly oscillatory function values around 0. On the other hand, the TN result converges quickly with $G$ as shown in Fig.~\ref{fig:pol_qmcb}. 
Although details of implementation make it difficult to interpret small differences in absolute timing, the rapid  convergence of the TN integration means that it is orders of magnitude more efficient than quasi-MC for high accuracies. However, the TN contraction cost becomes prohibitive at large $k$, due to using a $\chi=2^k$ bond dimension. 

For $c>0$, the integrand is more positive and we expect the TN to be more compressible, making it possible to integrate the function for larger $k$. In Fig.~\ref{fig:pol_qmcc}, we plot the convergence of the integral versus time for TN integration and quasi-MC for $N=50$, $k=30$ and various $c$. 
The TN data points correspond to $G=12$ Gauss-Legendre quadrature for increasing bond dimension $\chi\in[60,560]$. (To justify the choice of fixed $G=12$, Fig.~\ref{fig:pol_qmcd} shows the convergence of the TN contraction result as a function of $G$ for various $c$ at fixed $\chi$ (solid/dashed/dotted lines correspond to $\chi=180/120/60$); at $G=12$, the quadrature errors are below $10^{-6}$ for $c=0.3$ and below $10^{-8}$ for $c=0.4$ and $c=0.5$, and these errors remain essentially constant as $\chi$ increases from 60 to 180). As expected, for a fixed computational time, the integration error decreases for both methods with increasing $c$. Both methods display relatively quick convergence to a loose threshold, with slower convergence to tighter thresholds. The TN integration converges relatively smoothly to high accuracy, while it appears  difficult to converge quasi-MC systematically to high accuracy, which again opens up a significant timing advantage for TN integration over quasi-MC.

\section{Multivariable Gaussian integrals in a hypercube}

\begin{figure}%
    \centering
    \begin{subfigure}[b]{0.15\linewidth}
    \centering
    \includegraphics[width=\textwidth]{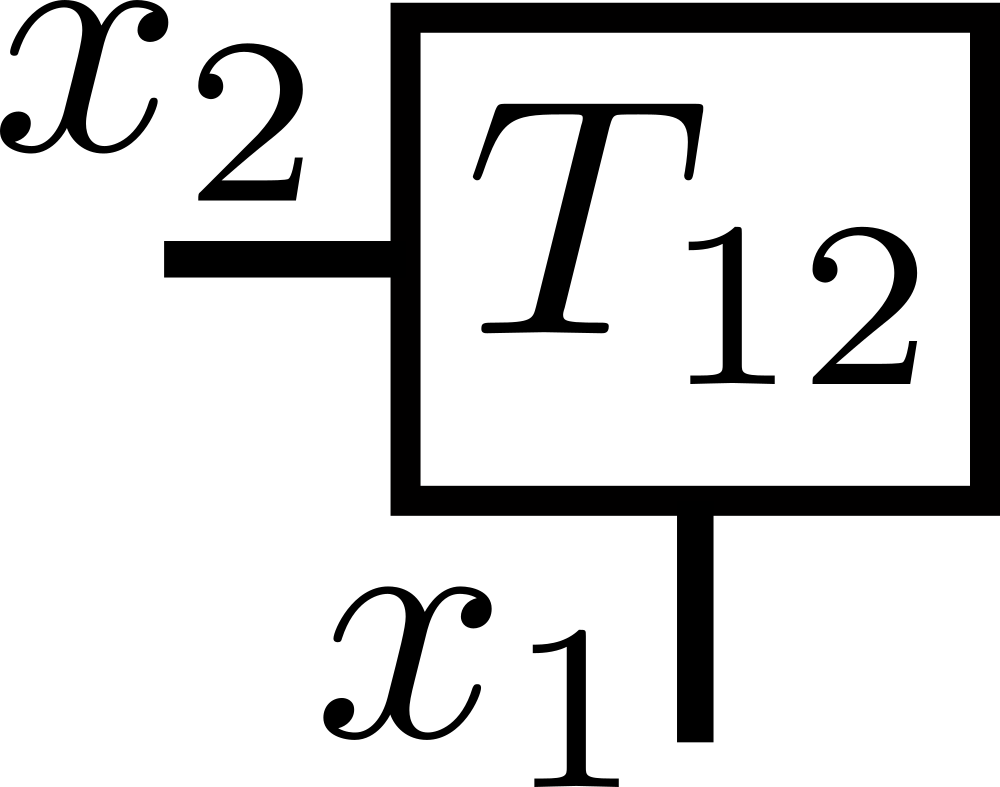}
    \caption{}
    \label{fig:gaussTNa}
    \end{subfigure}
    \begin{subfigure}[b]{0.3\linewidth}
    \centering
    \includegraphics[width=\textwidth]{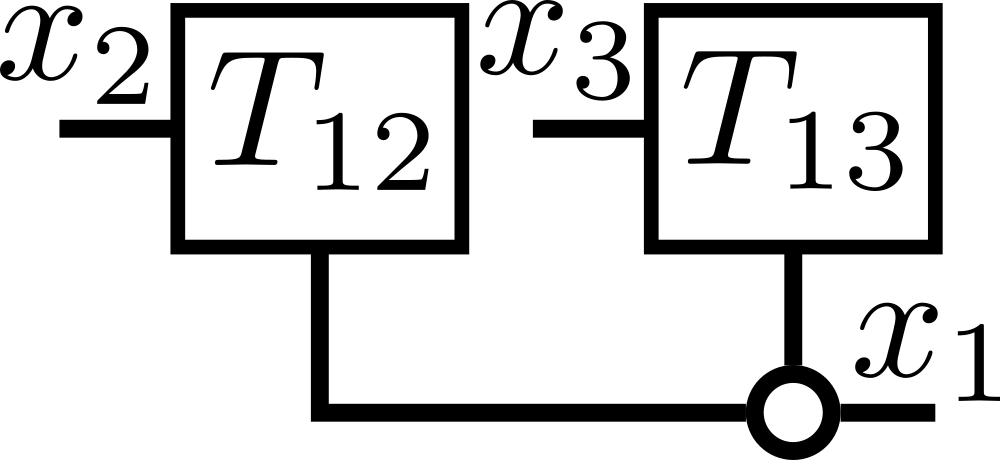}
    \caption{}
    \label{fig:gaussTNb}
    \end{subfigure}
    \begin{subfigure}[b]{0.45\linewidth}
    \centering
    \includegraphics[width=\textwidth]{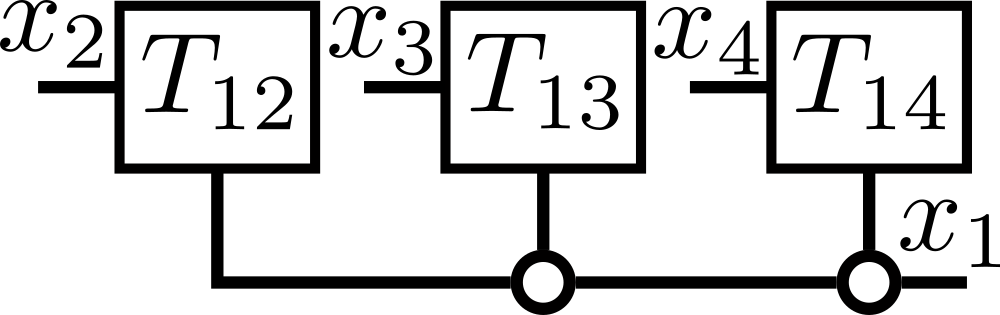}
    \caption{}
    \label{fig:gaussTNc}
    \end{subfigure}
    
    \begin{subfigure}[b]{0.32\linewidth}
    \centering
    \includegraphics[width=\textwidth]{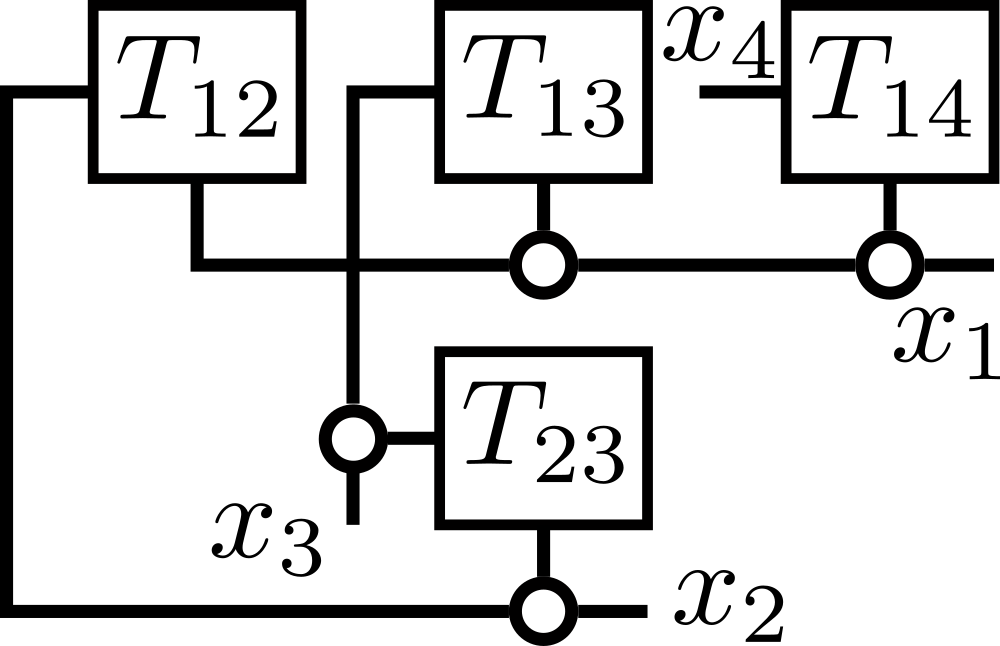}
    \caption{}
    \label{fig:gaussTNd}
    \end{subfigure}
    \begin{subfigure}[b]{0.32\linewidth}
    \centering
    \includegraphics[width=\textwidth]{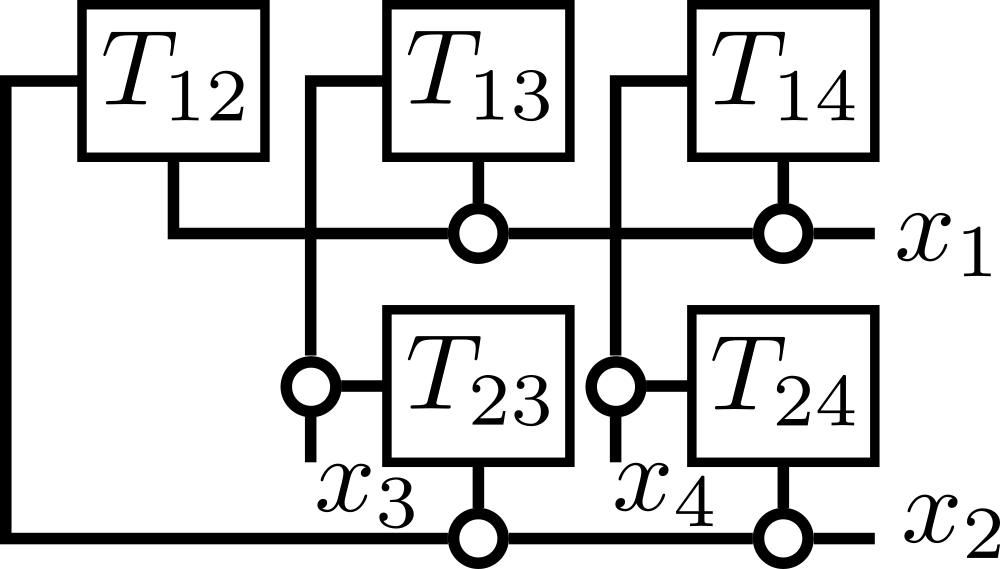}
    \caption{}
    \label{fig:gaussTNe}
    \end{subfigure}
    \begin{subfigure}[b]{0.32\linewidth}
    \centering
    \includegraphics[width=\textwidth]{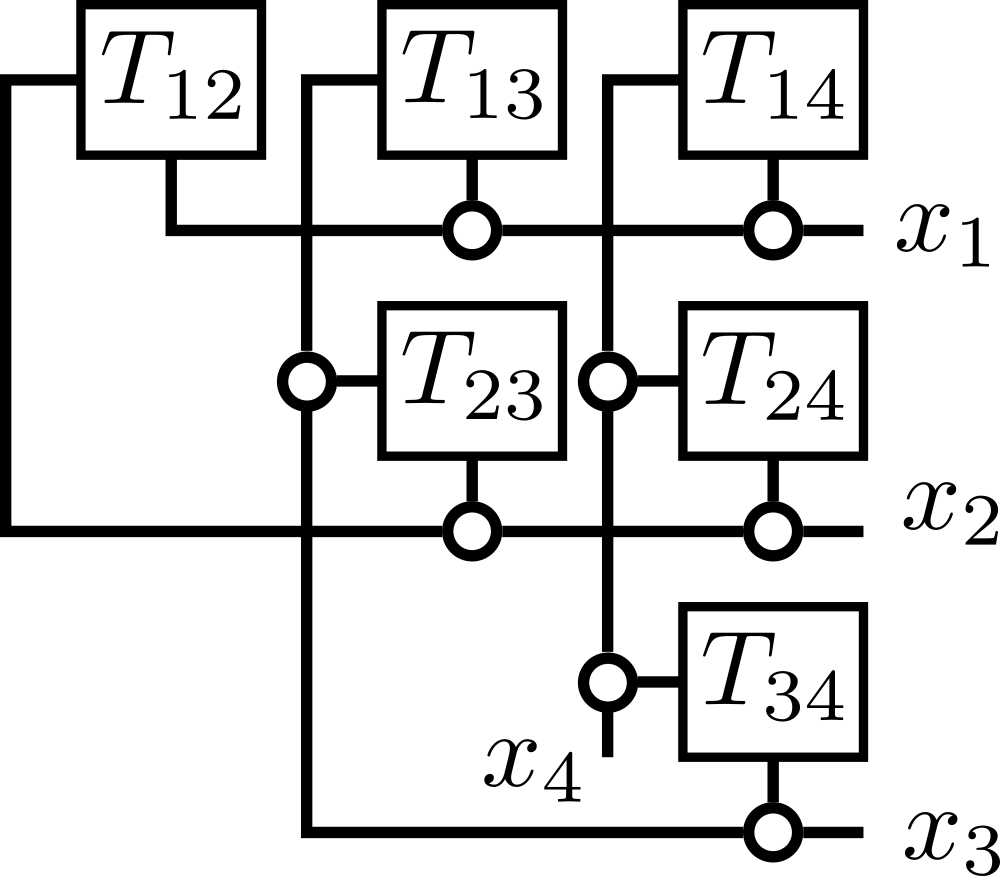}
    \caption{}
    \label{fig:gaussTNf}
    \end{subfigure}
    
    \begin{subfigure}[b]{0.55\linewidth}
    \centering
    \includegraphics[width=\textwidth]{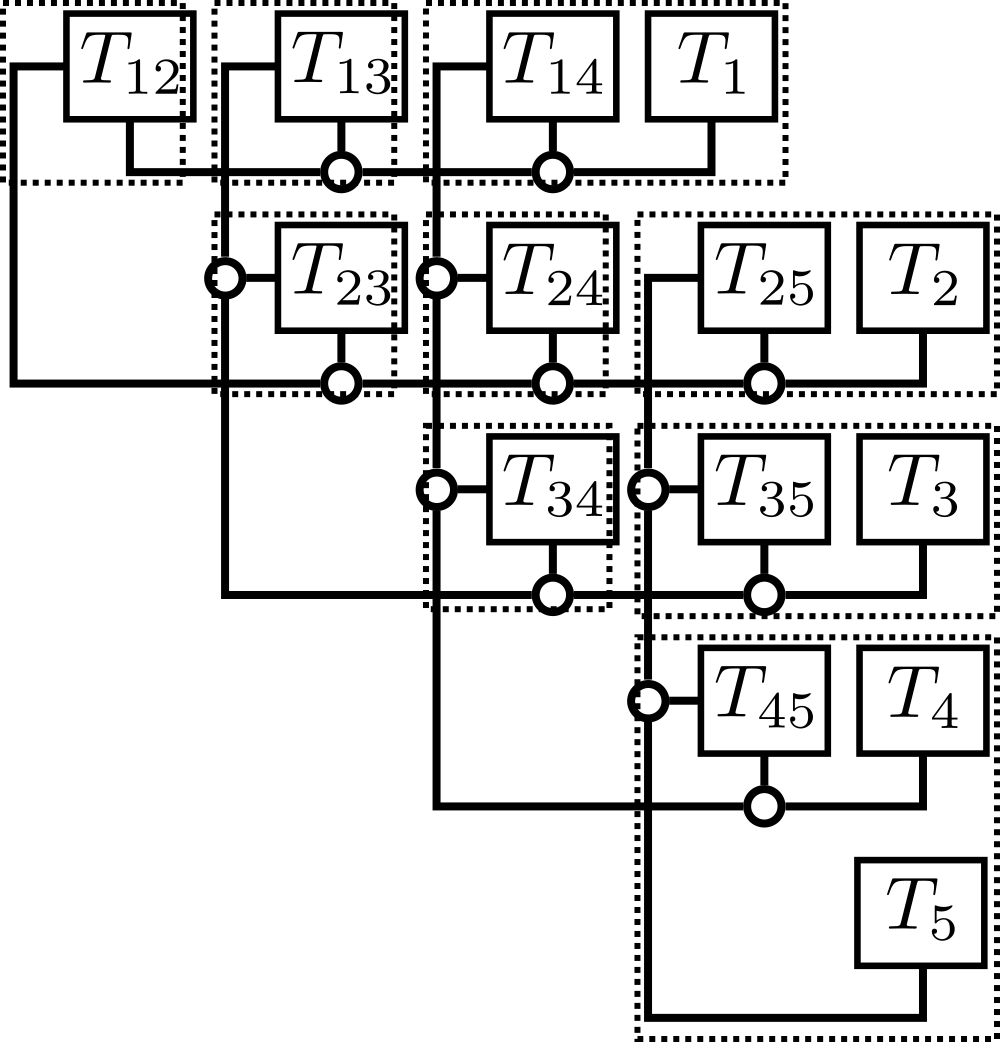}
    \caption{}
    \label{fig:gaussTNg}
    \end{subfigure}
    \begin{subfigure}[b]{0.4\linewidth}
    \centering
    \includegraphics[width=\textwidth]{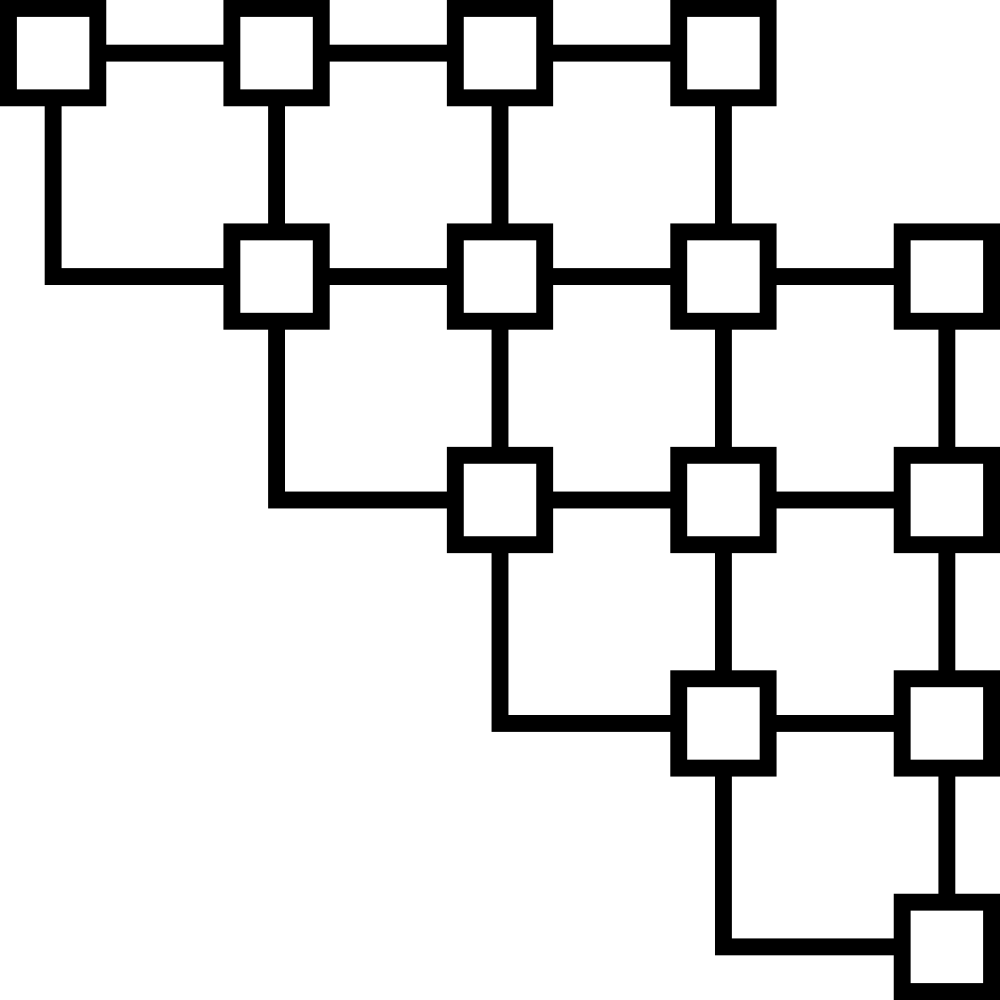}
    \caption{}
    \label{fig:gaussTNh}
    \end{subfigure}
    \caption{Steps for building the arithmetic circuit TN for a multivariable Gaussian with $N=5$ variables, and width $W=3$. Circles represent $\mathsf{COPY}$ tensors. }%
    \label{fig:gaussTN}%
\end{figure}

As another example, we consider the multivariable Gaussian integral over a finite hypercube
\begin{align}\label{eqn:gaussian}
Z=&\int_{\Omega}dx_1...dx_N\exp{(-\sum_{ij} A_{ij}x_ix_j)}
\end{align}
where $A$ is a $N\times N$ matrix and $\Omega=[-1,1]^N$. The expression directly corresponds to a 
tensor network contraction of tensors $(T_i)$ and $(T_{ij})$ for all $i<j$ where
\begin{align}
&(T_i)_{x_i}=\exp{(-A_{ii}x_i^2)}\\
&(T_{ij})_{x_i,x_j}=\exp{(-(A_{ij}+A_{ji})x_ix_j)}.
\end{align}
Note that all tensors enter into the final function via multiplication only, thus there is no need for control legs as discussed in section~\ref{sec:basic}, and they are omitted in the figures. 

The structure of $A$ plays an important role in the the cost of approximability of the integral. In the following we consider band-diagonal $A$ with width $W$, i.e. $A_{ij}\ne0$ only if $|i-j|<=W$ (dense $A$ corresponds to $W=N-1$). In Fig.~\ref{fig:gaussTN}, we show a systematic construction of the TN representation of the integral for $N=5$, $W=3$, which can easily be generalized to arbitrary $N$, $W$. We start with a TN consisting of only $T_{12}$ (representing $\exp{(-A_{12} x_1 x_2)}$) as in Fig.~\ref{fig:gaussTNa}. For each $j=2,\ldots,1+W$, we add to the TN the tensor $T_{1j}$ and a $\mathsf{COPY}$ tensor to account for the additional occurrence of $x_1$ as in Figs.~\ref{fig:gaussTNb},~\ref{fig:gaussTNc}. For example, in Fig.~\ref{fig:gaussTNc}, the
4 open legs represent the variables $x_1,x_2, x_3, x_4$, and the overall tensor network represents $\exp(-(A_{12}x_1x_2 + A_{13}x_1x_3 + A_{14}x_1x_4))$. 

Next, for $i=2$, and for each $j=i+1,\ldots,i+W$, we add to the TN the tensor $T_{ij}$ and two $\mathsf{COPY}$ tensors, e.g. in Figs.~\ref{fig:gaussTNd},~\ref{fig:gaussTNe}, the $\mathsf{COPY}$ tensors connected by the horizontal bonds correspond to copying the variable $x_i$, and the $\mathsf{COPY}$ tensors connected by the vertical bonds correspond to copying $x_j$ for each $j=3,4,...$. Iterating the previous step for each $i=3,\ldots,N-1$ as in Fig.~\ref{fig:gaussTNf}, one adds all $T_{ij}$ to the TN. Finally we add $T_i$ to the TN for each $i=1,\ldots,N$ as in Fig.~\ref{fig:gaussTNg}. Contracting the tensors in each dashed box, one obtains a 2D-TN shown in Fig.~\ref{fig:gaussTNh}. The TN obtained from the procedure above clearly has a direct correspondence with the structure of $A$. If $A$ is dense, then the graph is a triangular network, and if $A$ is banded, then the network is also banded.

\subsection{Fixed width compressibility}
\label{sec:gaussian_fixed}

\begin{figure}%
    \centering
    \begin{subfigure}[b]{0.24\linewidth}
    \centering
    \includegraphics[width=\textwidth]{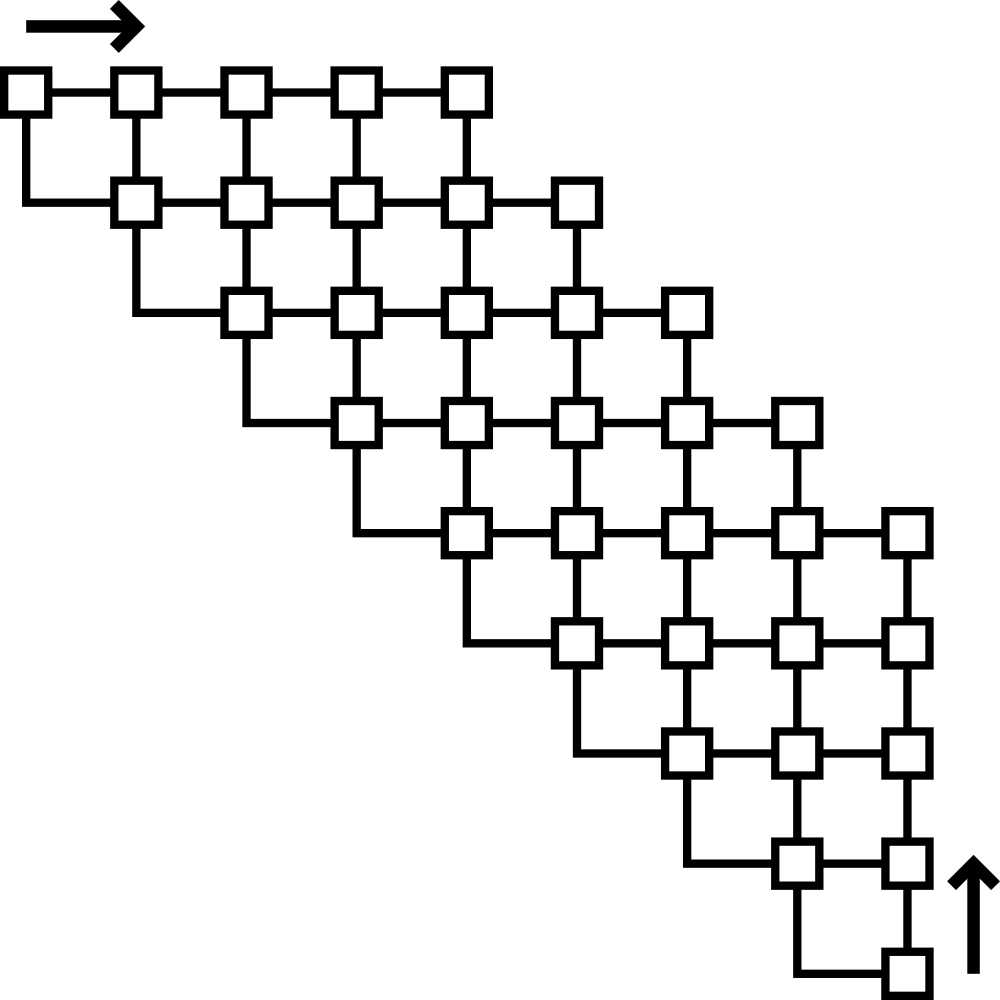}
    \caption{}
    \end{subfigure}
    \begin{subfigure}[b]{0.24\linewidth}
    \centering
    \includegraphics[width=\textwidth]{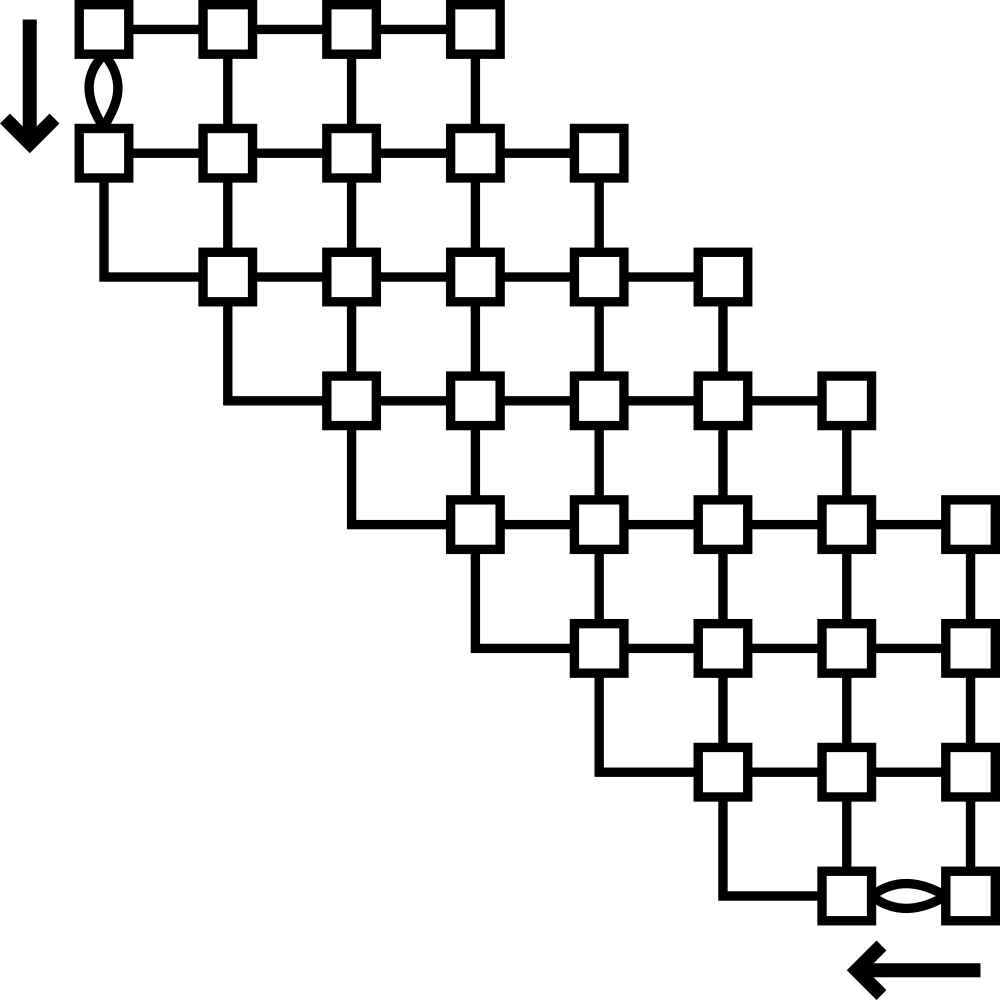}
    \caption{}
    \end{subfigure}
    \begin{subfigure}[b]{0.24\linewidth}
    \centering
    \includegraphics[width=\textwidth]{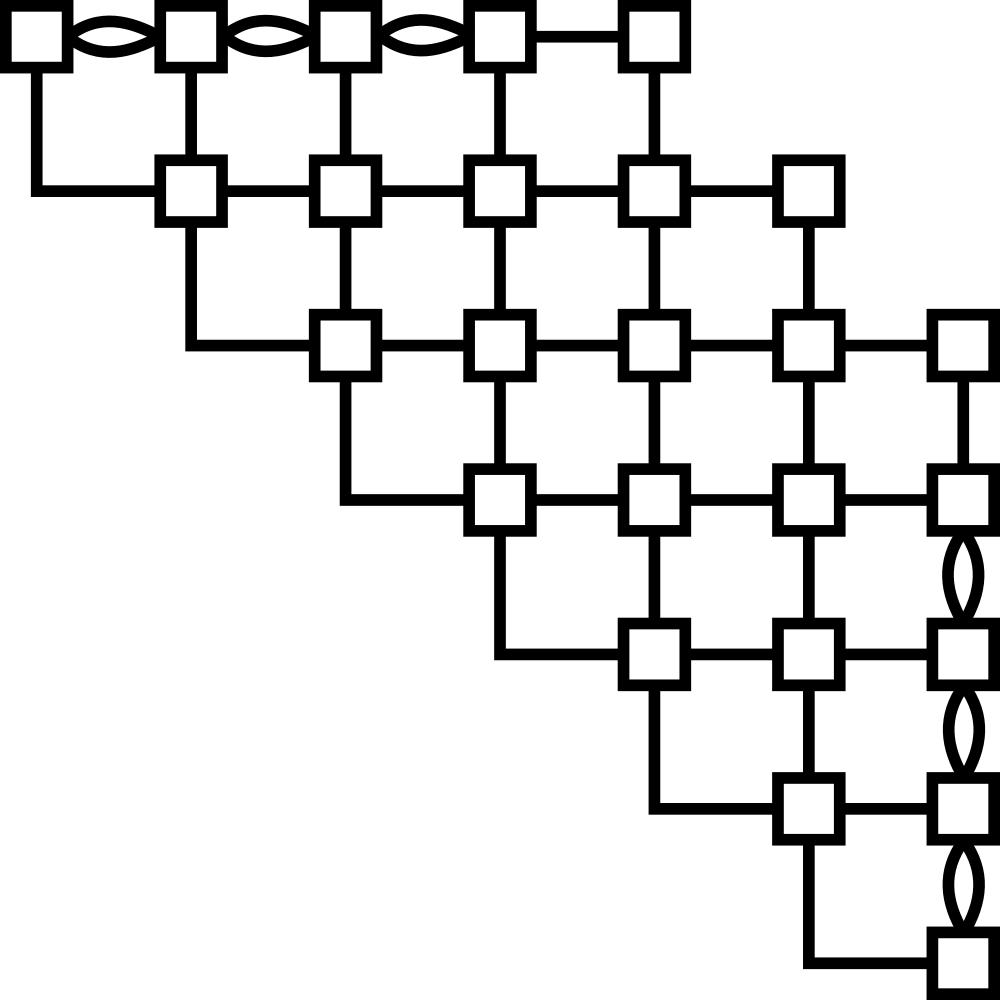}
    \caption{}
    \end{subfigure}
    \begin{subfigure}[b]{0.24\linewidth}
    \centering
    \includegraphics[width=\textwidth]{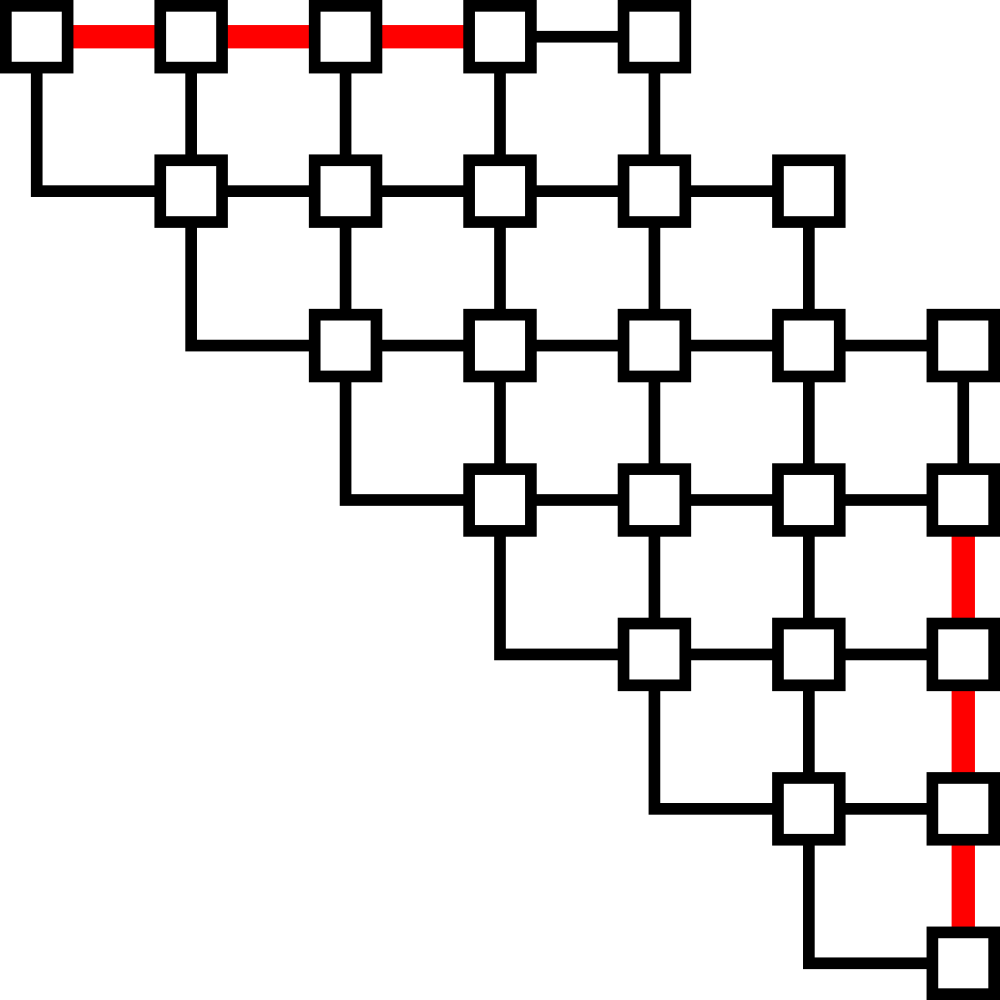}
    \caption{}
    \end{subfigure}
    \caption{Variant of boundary contraction used to contract the arithmetic circuit TN for the multivariable Gaussian integral. }%
    \label{fig:gaussian_compress}%
\end{figure}

\begin{figure}[]
    \centering
    \begin{subfigure}[b]{0.8\linewidth}
    \centering
    \includegraphics[width=\textwidth]{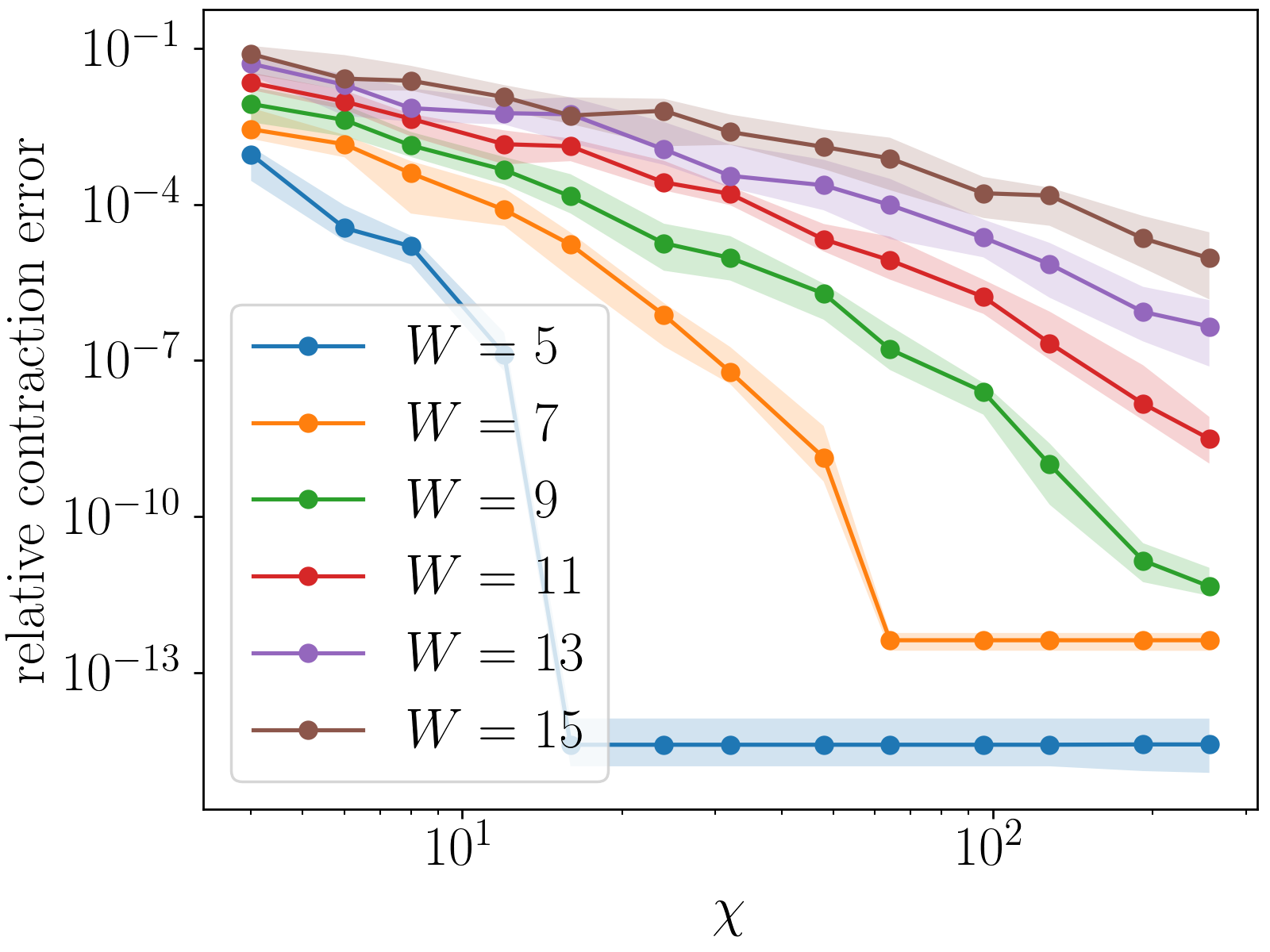}
    \caption{}
    \label{fig:gaussian_W1}
    \end{subfigure}
    \begin{subfigure}[b]{0.8\linewidth}
    \centering
    \includegraphics[width=\textwidth]{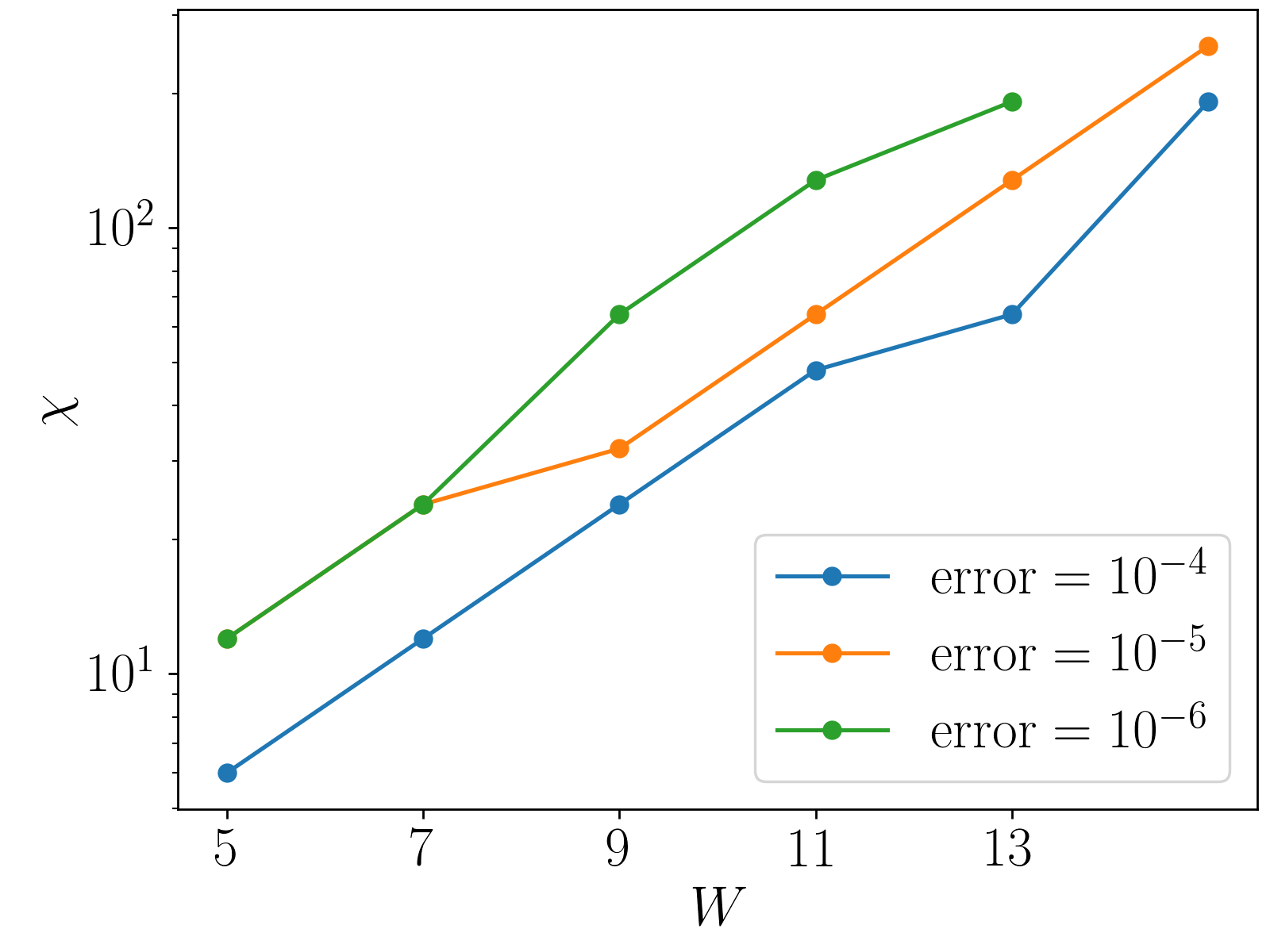}
    \caption{}
    \label{fig:gaussian_W2}
    \end{subfigure}
    \caption{TN integration for a multivariable Gaussian with exactly banded matrix $A$ in the hypercube.  Accuracy v.s. bond dimension $\chi$ for number of variables $N=30$, number of points per variable $G=4$. (a) Median and interquartile range of error over 20 random instances of the Gaussian matrix $A$. (b) Bond dimension to achieve various median contraction errors. For details see Sec.~\ref{sec:gaussian_fixed}.}
    \label{fig:gaussian_W}
\end{figure}

For banded $A$ (fixed $W$ with $N$), the TN has a quasi-1D structure along the diagonal direction. Thus it is always possible to contract the network {\textit{exactly}} with cost linear in $N$ and exponential in $W$ along the diagonal direction. If $W$ is large, it may still be too expensive to use exact contraction, but the regular structure lends itself to a variant of boundary contraction as shown in Fig.~\ref{fig:gaussian_compress}. 

In this case we limit the maximum bond dimension $\chi$ during the approximate contraction, with $\chi$ expected to scale like $\sim e^W$ to achieve a fixed relative error. We show the relative contraction error w.r.t. $\chi$ for various $W$ in Fig.~\ref{fig:gaussian_W1}, and the bond dimension to achieve various median contraction error as a function of $W$ in Fig.~\ref{fig:gaussian_W2}, for a $30\times30$ band-diagonal matrix $A$ with width $W$ and random nonzero elements in $[-1,1]$. 
We see that the difficulty of compression indeed scales exponentially with $W$ from the linear trend in the log-linear plot. 

\subsection{Approximately banded case}
\label{sec:gaussian_approx}

\begin{figure}[]
    \centering
    \begin{subfigure}[b]{0.8\linewidth}
    \centering
    \includegraphics[width=\textwidth]{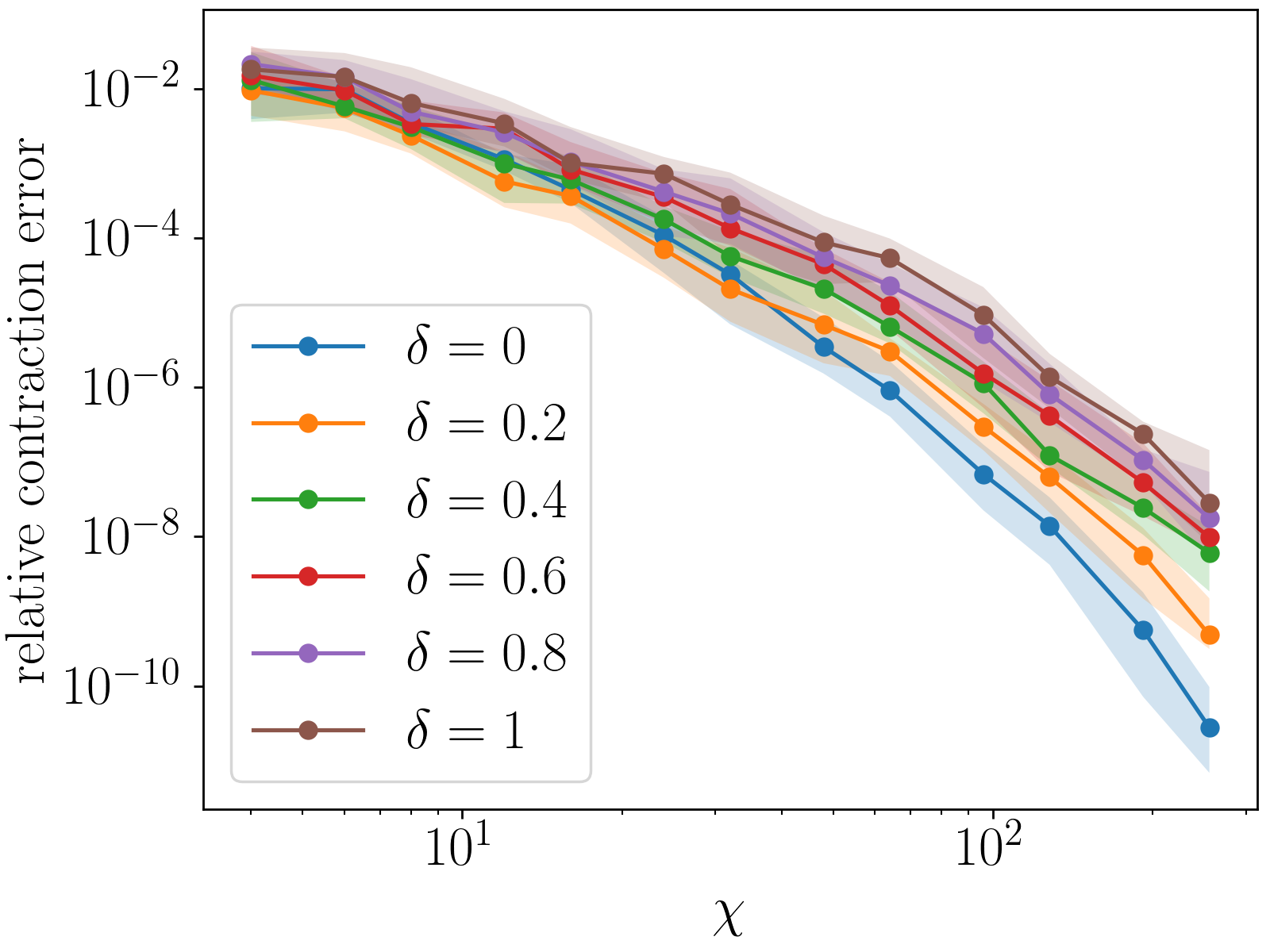}
    \caption{}
    \label{fig:gaussian_approxW1}
    \end{subfigure}
    \begin{subfigure}[b]{0.8\linewidth}
    \centering
    \includegraphics[width=\textwidth]{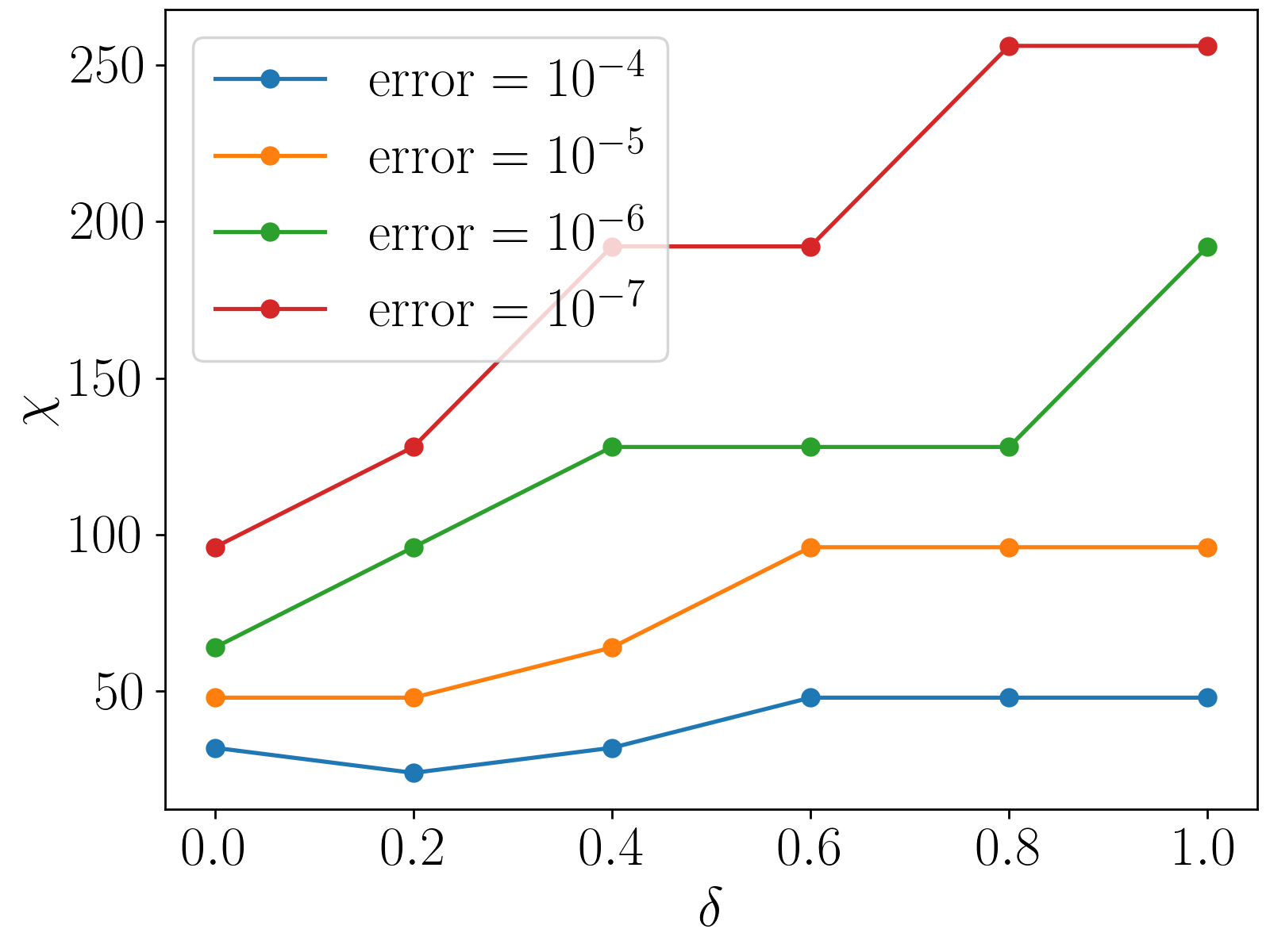}
    \caption{}
    \label{fig:gaussian_approxW2}
    \end{subfigure}
    \caption{TN integration for a multivariable Gaussian with approximately banded matrix $A$ in the hypercube.  Accuracy v.s. bond dimension $\chi$ for number of variables $N=30$, number of points per variable $G=4$
     (a) Median and interquartile range over 20 random instances. (b) Bond dimension to achieve various median contraction error. For details see Sec.~\ref{sec:gaussian_approx}.}
    \label{fig:gaussian_approxW}
\end{figure}

Above, we demonstrate that for moderate width $W$, the quasi-1D TN can be exactly contracted. As in the perturbed case for the polynomial example, we can also ask if the TN for an approximately band-diagonal $A$ remains compressible. We thus consider $A$ of the form
\begin{align}\label{eqn:gaussian_perturb}
A=A_W+\delta\cdot \Tilde{A}_W
\end{align}
where $A_W$ is an $N\times N$ band-diagonal matrix with width $W$ such that $(A_W)_{ij}\ne0$ only if $|i-j|\leq W$, $\Tilde{A}_W$ is an $N\times N$ matrix such that $(\Tilde{A}_W)_{ij}\ne0$ only if $|i-j|>W$, and non-zero elements of $A_W$, $\Tilde{A}_W$ take random values in $[-1,1]$. Then the magnitude of the perturbation of $A$ away from $A_W$ is controlled by $\delta$. Note that the existence of the perturbation makes the TN a full 2D triangle, rather than quasi-1D. 

We show the relative contraction error w.r.t $\chi$ for various $\delta$ in Fig.~\ref{fig:gaussian_approxW1} and the bond dimension required to achieve various median errors as a function of $\delta$ in Fig.~\ref{fig:gaussian_approxW2} for $A$ of dimension $20\times20$. We require increasingly large $\chi$ to reach the same error as $\delta$ increases. However, Fig.~\ref{fig:gaussian_approxW2} is a linear-linear plot, showing the required $\chi$ appears to grow only linearly with the size of the perturbation $\delta$. 

\subsection{Comparison with quasi-Monte Carlo}

\begin{figure}%
    \centering
    \begin{subfigure}[b]{0.8\linewidth}
    \centering
    \includegraphics[width=\textwidth]{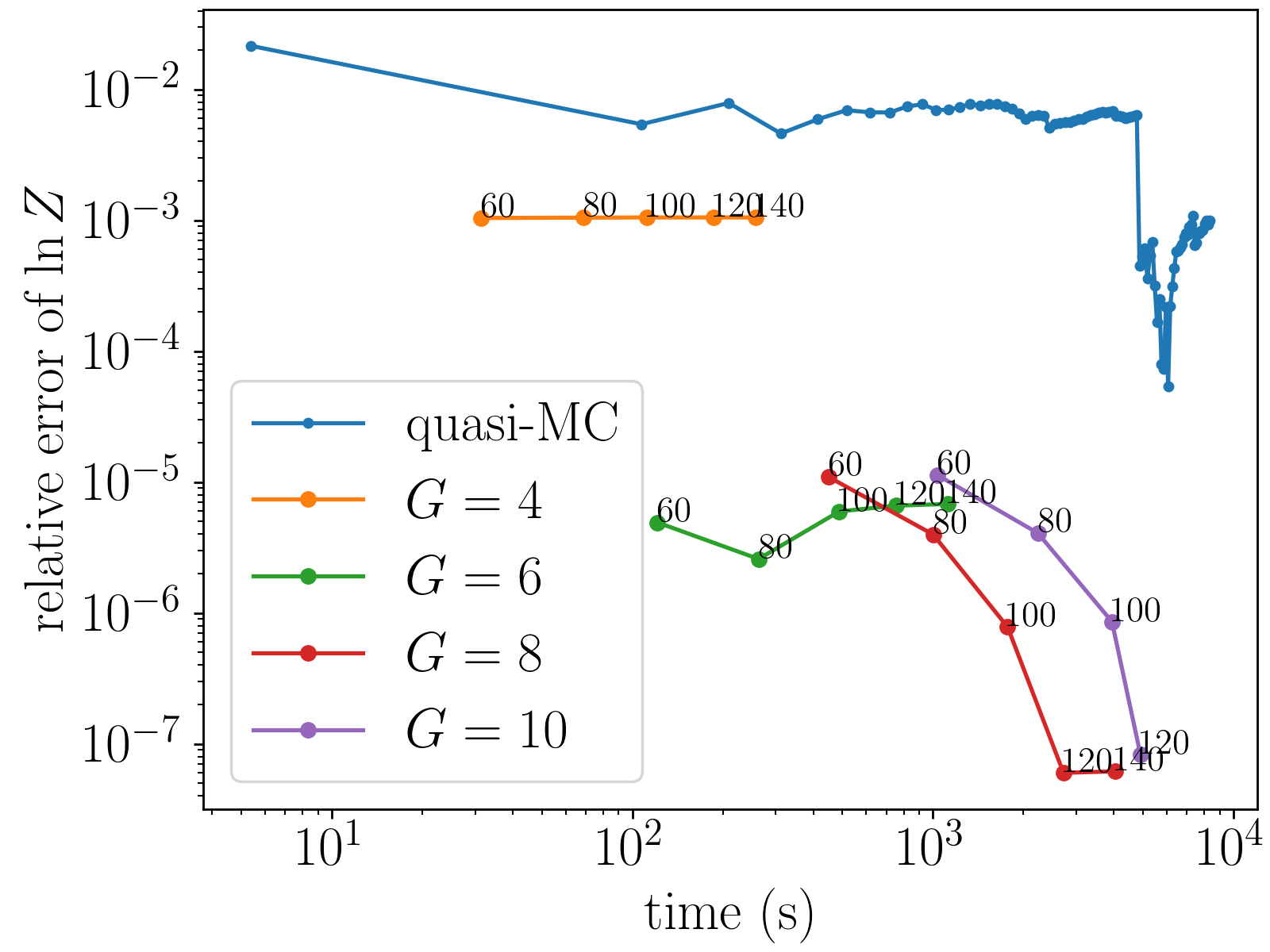}
    \caption{}
    \label{fig:gaussian_qmc1}
    \end{subfigure}
    \begin{subfigure}[b]{0.8\linewidth}
    \centering
    \includegraphics[width=\textwidth]{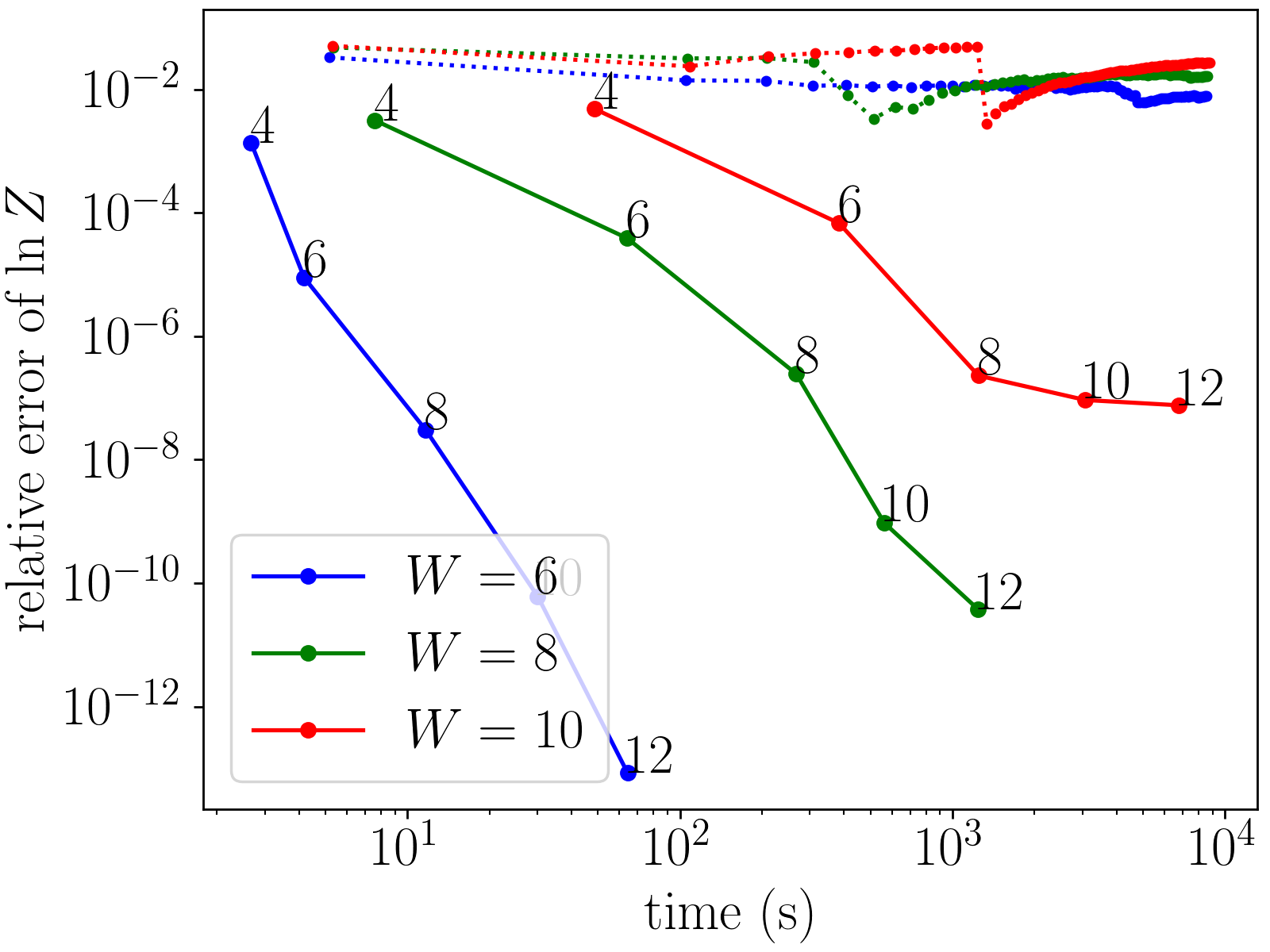}
    \caption{}
    \label{fig:gaussian_qmc2}
    \end{subfigure}
    \caption{Comparison of TN and quasi-MC integration errors for the multivariable Gaussian integral $Z$ in the hypercube, plotted against runtime. (Error plotted is for $\ln Z$ to better distinguish the curves).  $N$ is the number of variables, $G$ is number of quadrature points per variable in the TN. (a) Exact result taken as TN partition function at $G=10$, bond dimension $\chi=140$ (bond dimension error below $5\times10^{-7}$). The TN errors are plotted for $\chi=60,80,100,120,140$ for each $G$, which are the labels next to the TN data points. (b) Exact result is taken at $G=14$ for converged $\chi=80/140/200$ for $W=6/8/10$ with contraction error $\leq10^{-15}/10^{-9}/10^{-6}$. TN/quasi-MC results are represented by solid/dotted lines. The label next to each data point corresponds to its $G$ value. }%
    \label{fig:gaussian_qmc}%
\end{figure}

We now compare the efficiency of TN integration for Gaussian integrals with banded (or approximately banded) $A$ with that of quasi-Monte Carlo. We consider two cases for $A$ of the form Eq.~(\ref{eqn:gaussian_perturb}) (a) $N=50$, $W=5$, $\delta=0.1$ with random nonzero elements of $A_W$, $\Tilde{A}_W$ in $[-1,1]$ and case (b) $N=50$, $W=6/8/10$, $\delta=0$ with random nonzero elements of $A_W$ in $[-1,1]$. 

In Fig.~\ref{fig:gaussian_qmc}, we show the relative error of $\ln{Z}$ from TN and quasi-MC integration for cases (a) and (b). For the quasi-MC results, we directly computed the exponent $-\sum_{ij}A_{ij}x_ix_j$ for each sample point $x$ with batched linear algebra code using Python and JAX\cite{jax2018github}, and accumulated the log of the integral using the identity
\begin{align*}
\ln{(e^a+e^b)}=a+\ln{(1+e^{b-a})}.
\end{align*}
The sample points were generated using QMCPY\cite{QMCPy} with Niederreiter generating matrices in batches of size $10^7$, and plotted every 20 batches. The TN integration was performed using \textsc{Quimb}~\cite{gray2018quimb} with Gauss–Legendre quadrature and the exact value of the integral was estimated from
TN integration using converged $G$ and $\chi$ (indicated in the caption of Fig.~\ref{fig:gaussian_qmc}). All calculation are done on  on a single Intel(R) Xeon(R) CPU E5-2697 v4 processor. 

For the TN integration data in Fig.~\ref{fig:gaussian_qmc}a (corresponding to case (a)) each line labeled by $G$ is a sequence of TN estimates of $\ln Z$ at the corresponding $G$ with increasing $\chi=60,80,\ldots,140$. For $G=4$ and $G=6$, the value of $\ln Z$ plateaus at large $\chi$, corresponding to the intrinsic quadrature error for that $G$. For $G=8$ and $G=10$, the TN results are well converged with quadrature, and the small perturbation strength $\delta=0.1$ allows for quick convergence with $\chi$. In comparison, quasi-MC struggles with this integrand, and in fact shows little systematic convergence behavior over the $\sim 1600\times10^{7}$ samples. 

For case (b), shown in Fig.~\ref{fig:gaussian_qmc}b, we plot the TN data points for increasing $G=4,6,\ldots,12$ at converged $\chi$ for each $W$ (solid lines). This allows us to examine the speed of convergence of the quadrature, which is very fast. Similar to case (a) quasi-MC struggles to reduce the relative error below $10^{-2}$ even with more than $1600\times10^{7}$ samples. Thus, in both cases, we see a substantial advantage of TN integration over quasi-MC for this class of integrands.

\section{A theoretical example of speedup versus quasi-Monte Carlo}

\begin{figure}%
    \centering
    \begin{subfigure}[b]{0.49\linewidth}
    \centering
    \includegraphics[width=0.9\textwidth]{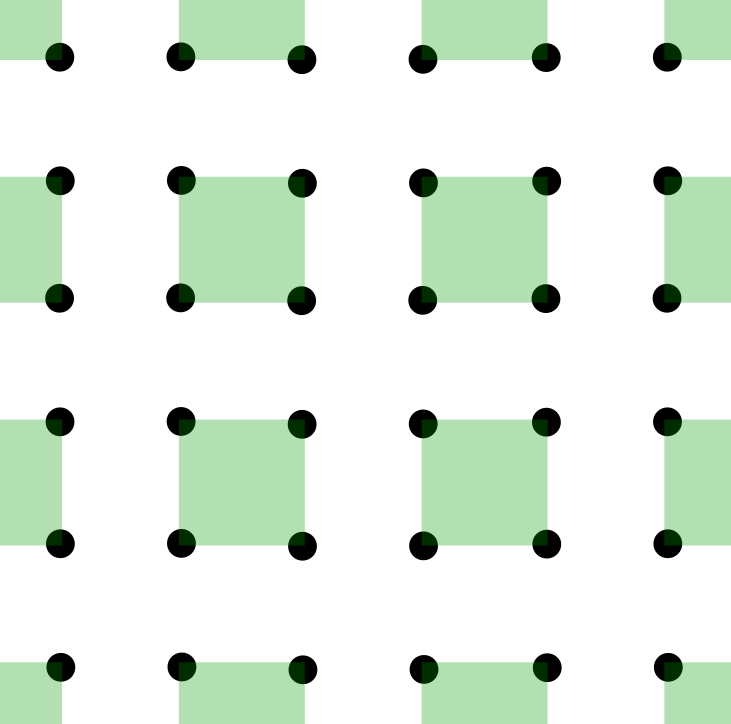}
    \caption{}
    \end{subfigure}
    \begin{subfigure}[b]{0.49\linewidth}
    \centering
    \includegraphics[width=0.9\textwidth]{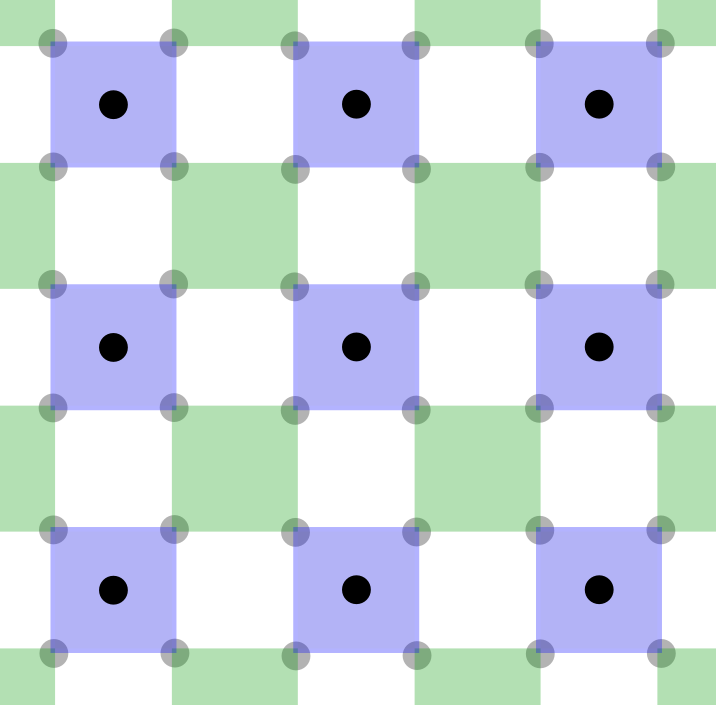}
    \caption{}
    \end{subfigure}
    
    \begin{subfigure}[b]{0.49\linewidth}
    \centering
    \includegraphics[width=\textwidth]{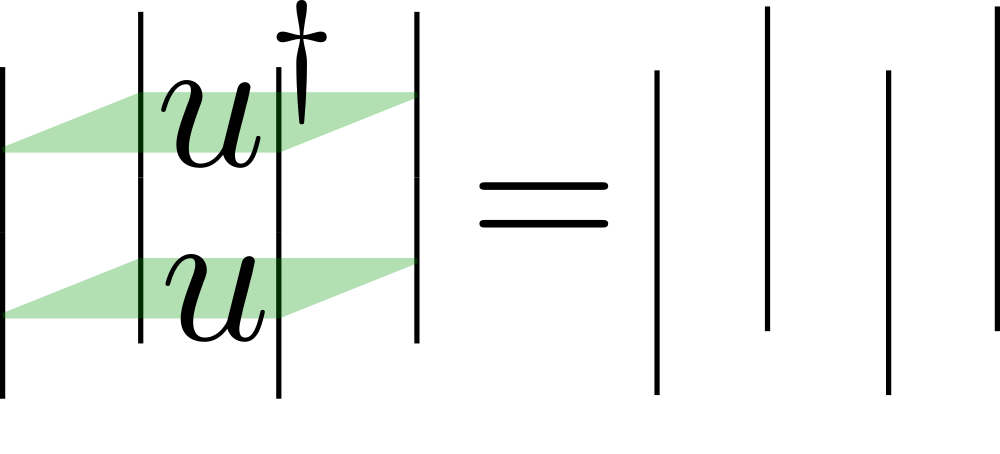}
    \caption{}
    \label{fig:unitary}
    \end{subfigure}
    \begin{subfigure}[b]{0.49\linewidth}
    \centering
    \includegraphics[width=0.8\textwidth]{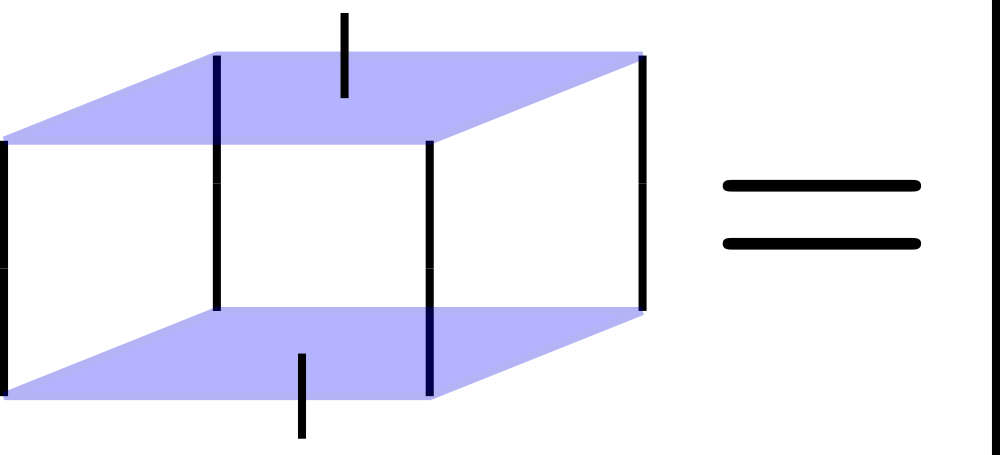}
    \caption{}
    \label{fig:isometry}
    \end{subfigure}
    \caption{$\mathrm{TN}_{i_1,\ldots,i_N}$ as a 2D-MERA, where the legs $i_1,\ldots,i_N$ lie on a 2D-lattice. (a) Part of the 2D-lattice, where each solid dot represent a leg in $\{i_1,\ldots,i_N\}$. Each green square represent a unitary tensor. (b) After the layer of unitaries, a layer of isometries (blue squares) is applied to the legs $i_1,\ldots,i_N$. Note after a layer of isometries is applied, the number of legs in the layer is reduced. Traversing upwards through the MERA, one reaches the top tensor where all upwards legs are eliminated, and the MERA represents a TN with legs $i_1,\ldots,i_N$ on the bottom. (c) property of a unitary tensor. (d) property of an isometric tensor. }%
    \label{fig:mera}%
\end{figure}

Here we give an example of a class of integrals that can be computed easily with arithmetic TN methods but for which quasi-MC is hard, because there is a theoretical guarantee that the integrand value cannot be efficiently evaluated without exponential cost in the number of variables $N$.

Consider the general functional form
\begin{align}
f(x_1,...,x_N)=\prod_{n=1}^Ng(x_n)^{i_n}(\mathrm{TN})_{i_1,...,i_N}
\label{eq:generalmera}
\end{align}
where $g(x_n)$ are single variable functions and $(\mathrm{TN})_{i_1,...,i_N}$ is some fixed tensor (possibly represented as a tensor network) for all $i_n=0,1$. We can choose the structure and the values of the $\mathrm{TN}_{i_1,...,i_N}$ and properties of $g(x_n)$ to allow for the easy computation of certain quantities. In our case, we choose the TN to be a 2D multi-scale entanglement renormalization ansatz~\cite{Evenbly2009_mera,Evenbly2014_mera} (see Fig.~\ref{fig:mera}). In this case, the TN in Eq.~(\ref{eq:generalmera}) is constructed from tensors satisfying certain unitary and isometric properties such that the evaluation of the trace $\sum_{\{ i\}} |\mathrm{TN}_{i_1, \ldots, i_N}|^2=1$. The unitary and isometric constraints on the tensors and how they lead to the trivial trace are shown in Figs.~\ref{fig:unitary},~\ref{fig:isometry}. Note that it is known to be classically hard (as a function of $N$) to compute the value $\mathrm{TN}_{i_1, \ldots, i_N}$ for a MERA composed of arbitrary unitaries and isometries~\cite{Ferris2012165147}.

To extend these properties to the more general functional form, we impose a normalization condition on the single variable functions 
\begin{align}
&\int_{-1}^1dx|g(x)|^2=2\\%\approx\sum_{p}(w)_p|g(x[p])|^2=2,\\
&\int_{-1}^1dxg(x)=0%\approx\sum_{p}(w)_pg(x[p])=0
\end{align}
Then, suppose we want to integrate the probability of the function over $\Omega=[-1,1]^N$
\begin{align}\label{eqn:mera_norm}
I
=\int_{\Omega}|f(\{x_{n}\})|^2\prod dx_{n}.
%\approx
%\sum_{p_n}(w_n)_{p_n}|f(\{x_n[p_n]\})|^2
\end{align}
Because of the additional constraints we have imposed on the single variable functions, 
contraction of the tensors associated with $I$ can be done efficiently for matching pairs, which then simplify to multiples of the identity, yielding $2^N$ as the value of the integral, even though  sampling in the basis of $\{ x_n \}$ is computationally hard. 
This result is clearly contrived because, in addition to requiring a restricted class of tensors and functions, we must also contract the tensor network in a certain order. Furthermore, the structure of the tensors permits sampling in a different basis than $\{ x_n \}$~\cite{Ferris2012165146,Ferris2012165147}. Nonetheless, this is a concrete example supported by complexity arguments where using the tensor network structure of the function circuit leads to an exponential improvement in cost versus quasi-Monte Carlo methods which assume the function is a blackbox.

\section{Conclusions}

In this work we introduced an arithmetic circuit tensor network representation of multivariable functions and demonstrated its power for high dimensional integration. Compared to existing techniques to represent functions with tensor networks, the ability to use the arithmetic circuit construction removes the need for extra computation to find the tensor representation, while the circuit structure suggests a tensor network connectivity natural to the function. In our examples of high dimensional integration, we find that the tensor network representation allow us to circumvent the curse of dimensionality in many cases, exchanging the exponential cost dependence on dimension for a cost dependence on other circuit characteristics, for example, the number of nonlinear circuit operations. In practice, we find superior performance to quasi-Monte Carlo integration across a range of dimensionalities and accuracies.

While the work here focused on integration as an example, we envision the arithmetic circuit tensor network construction to be powerful also in differential equation applications. Here, connections with existing tensor network techniques are intriguing, as are applications of these ideas to many other areas where tensor networks are currently employed, such as for many-body simulations.

\bibliography{export}
\end{document}